\documentclass[11pt]{amsart}
\newcommand{\preprint}[1]{}

\newcommand{\hide}[1]{}

\newcommand{\footnoteremember}[2]{\footnote{#2}\newcounter{#1}\setcounter{#1}{\value{footnote}}} 
\newcommand{\footnoterecall}[1]{\footnotemark[\value{#1}]} 

\usepackage{amssymb}
\usepackage{amsbsy}
\usepackage{amscd}
\usepackage{amsmath}
\usepackage{amsthm}
\usepackage{geometry}

\input xy
\xyoption{all}

\numberwithin{equation}{section}

\theoremstyle{plain}
\newtheorem{thm}{Theorem}[section]
\newtheorem{assumption}[thm]{Assumption}
\newtheorem{prop}[thm]{Proposition}

\newtheorem{problem}[thm]{Problem}
\newtheorem{conj}[thm]{Conjecture}
\newtheorem{cor}[thm]{Corollary}
\newtheorem{lem}[thm]{Lemma}

\theoremstyle{definition}
\newtheorem{defi}[thm]{Definition}

\theoremstyle{remark}
\newtheorem{example}[thm]{Example}
\newtheorem{question}[thm]{Question}
\newtheorem{rem}[thm]{Remark}
\newtheorem{caution}[thm]{Caution}

\newcommand{\A}{{\mathcal A}}

\newcommand{\BK}{{\mathcal B}{\mathcal K}}
\newcommand{\C}{{\mathcal C}}
\newcommand{\D}{{\mathcal D}}
\newcommand{\E}{{\mathcal E}}
\newcommand{\EE}{{\mathbb E}}

\newcommand{\FE}{{\mathcal F}{\mathcal E}}
\newcommand{\FR}{{\mathcal F}{\mathcal R}}
\newcommand{\FU}{{\mathcal F}{\mathcal U}}
\newcommand{\MV}{{\mathcal M}{\mathcal V}}

\newcommand{\HH}{{\mathbb H}}

\newcommand{\K}{{\mathcal K}}

\newcommand{\KT}{{\mathcal K}{\mathcal T}}
\newcommand{\LB}{{\mathcal L}}
\newcommand{\M}{{\mathcal M}}
\newcommand{\FM}{{\mathfrak M}}
\newcommand{\N}{{\mathcal N}}
\renewcommand{\P}{{\mathcal P}}

\newcommand{\PP}{{\mathbb P}}

\newcommand{\R}{{\mathcal R}}

\newcommand{\Spe}{{\mathcal S}pe}
\newcommand{\T}{{\mathcal T}}
\newcommand{\U}{{\mathcal U}}
\newcommand{\V}{{\mathcal V}}

\newcommand{\X}{{\mathcal X}}

\newcommand{\RealNumbers}{{\mathbb R}}
\newcommand{\Integers}{{\mathbb Z}}
\newcommand{\ComplexNumbers}{{\mathbb C}}
\newcommand{\RationalNumbers}{{\mathbb Q}}

\newcommand{\linsys}[1]{{\mid}#1{\mid}}

\newcommand{\LongRightArrowOf}[1]{\stackrel{#1}{\longrightarrow}}

\newcommand{\StructureSheaf}[1]{{\mathcal O}_{#1}}
\newcommand{\EndProof}{\hfill  $\Box$}
\newcommand{\restricted}[2]{#1_{\mid_{#2}}}

\newcommand{\rank}{{\rm rank}}

\renewcommand{\div}{{\rm div}}
\newcommand{\NS}{{\rm NS}}
\newcommand{\Pic}{{\rm Pic}}
\newcommand{\Pex}{{\mathcal P}ex}
\newcommand{\Eff}{{\mathcal E}{\rm ff}}
\newcommand{\Peff}{{\mathcal P}{\rm eff}}
\newcommand{\Sym}{{\rm Sym}}
\newcommand{\Ext}{{\rm Ext}}

\newcommand{\Hom}{{\rm Hom}}
\newcommand{\Aut}{{\rm Aut}}
\newcommand{\Bir}{{\rm Bir}}

\newcommand{\Orient}{{\rm Orient}}
\newcommand{\orient}{{\rm orient}}
\newcommand{\period}{p}

\begin{document}
\title[A survey of Torelli and monodromy results]
{A survey of Torelli and monodromy results for holomorphic-symplectic varieties}
\author{Eyal Markman}
\address{Department of mathematics and statistics, University of Massachusetts, Amherst MA 01003}
\email{markman@math.umass.edu}

\begin{abstract}
We survey recent results about the  Torelli question for
holomorphic-symplectic varieties. Following are the main topics. 
A Hodge theoretic Torelli theorem. 
A study of the subgroup $W_{Exc}$, of the isometry group of the weight $2$ Hodge structure,
generated by reflection with respect to exceptional divisors.
A description of the birational K\"{a}hler cone as a fundamental domain for the $W_{Exc}$ action on
the positive cone.
A proof of a weak version of Morrison's movable cone conjecture. 
A description of the moduli spaces of polarized holomorphic symplectic varieties
as monodromy quotients of period domains of type IV. 
\end{abstract}

\maketitle

\tableofcontents

%***************************************************************************
%
%***************************************************************************
\section{Introduction}

An {\em irreducible holomorphic symplectic manifold} is a simply connected
compact K\"{a}hler manifold such that $H^0(X,\Omega^2_X)$ is one-dimensional,
spanned by an everywhere non-degenerate holomorphic $2$-form
\cite{beauville}. 
There exists a unique non-degenerate symmetric integral and primitive 
bilinear pairing $(\bullet,\bullet)$ on $H^2(X,\Integers)$ 
of signature $(3,b_2(X)-3)$, with the following property.
There exists a positive rational number $\lambda_X$, such that the equality 
\[
(\alpha,\alpha)^n \ \ = \ \ \lambda_X\int_X\alpha^{2n}
\]
holds for all $\alpha\in H^2(X,\Integers)$, where $2n=\dim_\ComplexNumbers(X)$
\cite{beauville}. If $b_2(X)=6$, then 
we require\footnote{The condition is satisfied automatically by the assumption
that the signature is $(3,b_2(X)-3)$, if $b_2\neq 6$.} 
further that $(\alpha,\alpha)>0$, for every K\"{a}hler class $\alpha$.
The pairing is called the {\em Beauville-Bogomolov pairing} and
$(\alpha,\alpha)$ is called the {\em Beauville-Bogomolov degree} of the class $\alpha$.

Let $S$ be a $K3$ surface. Then the Hilbert scheme (or Douady space, in the K\"{a}hler case)
$S^{[n]}$, of length $n$ zero-dimensional subschemes of $S$,
is an irreducible holomorphic symplectic manifold.
If  $n\geq 2$, then $b_2(S^{[n]})=23$ \cite{beauville}.
If $X$ is deformation equivalent to $S^{[n]}$, we will say that $X$ is of
{\em $K3^{[n]}$-type.}

Let $T$ be a complex torus with an origin $0\in T$. Denote by $T^{(n)}$ the
$n$-th symmetric product. 
Let $T^{(n)}\rightarrow T$ be the addition morphism.
The composite morphism
\[
T^{[n+1]} \longrightarrow T^{(n+1)} \longrightarrow T
\]
is an isotrivial fibration. Each fiber is a $2n$-dimensional irreducible holomorphic 
symplectic manifold, called a {\em generalized Kummer variety}, and denoted by $K^{[n]}(T)$  
\cite{beauville}.
If $n\geq 2$, then $b_2\left(K^{[n]}(T)\right)=7$. 

O'Grady constructed two additional irreducible holomorphic symplectic manifolds,
a $10$-dimensional example $X$ with $b_2(X)=24$, and a $6$-dimensional example
$Y$ with $b_2(Y)=8$ \cite{ogrady-ten,ogrady-six,rapagnetta}.

We recommend Huybrechts' 
excellent survey of the subject of irreducible holomorphic symplectic manifolds
\cite{huybrechts-norway}. The aim of this note is to survey developments related
to the Torelli problem, obtained by various authors since Huybrechts' survey 
was written. The most important, undoubtedly, is Verbitsky's proof of
his version of the Global Torelli Theorem \cite{verbitsky,huybrechts-bourbaki}. 
%Several new results are obtained along the way.
%**************************************************************************************
%
%**************************************************************************************
\subsection{Torelli Theorems}
We hope to convince the reader that the concepts of monodromy and 
parallel-transport operators are essential for any discussion of the Torelli problem. 

\begin{defi}
\label{def-monodromy}
Let $X$, $X_1$, and $X_2$ be irreducible holomorphic symplectic manifolds.
\begin{enumerate}
\item
\label{def-item-parallel-transport-operator}
An isomorphism $f:H^*(X_1,\Integers)\rightarrow H^*(X_2,\Integers)$
is said to be a {\em parallel-transport operator}, 
if there exist a smooth and proper family\footnote{Note
that the family may depend on the isomorphism $f$.}
$\pi:\X\rightarrow B$ of irreducible holomorphic symplectic manifolds, 
over an analytic base $B$,  points $b_i\in B$, isomorphisms $\psi_i:X_i\rightarrow \X_{b_i}$,
$i=1,2$, and a continuous path $\gamma:[0,1]\rightarrow B$, satisfying $\gamma(0)=b_1$,
$\gamma(1)=b_2$, such that the parallel transport in the local system
$R{\pi_*}\Integers$ along $\gamma$ induces the homomorphism
$\psi_{2_*}\circ f\circ \psi_1^*:H^*(\X_{b_1},\Integers)\rightarrow H^*(\X_{b_2},\Integers)$. 
An isomorphism $g:H^k(X_1,\Integers)\rightarrow H^k(X_2,\Integers)$
is said to be a {\em parallel-transport operator}, if it is the $k$-th graded summand of 
a  parallel-transport operator $f$ as above.
\item
An automorphism $f:H^*(X,\Integers)\rightarrow H^*(X,\Integers)$ is said to be
a {\em monodromy operator}, if it is a parallel transport operator.
\item
The {\em monodromy group} $Mon(X)$ is the subgroup\footnoteremember{group-footnote}{If $f\in Mon(X)$ is associated to 
a family $\pi':\X'\rightarrow B'$ via an isomorphism $X\cong\X'_{b'}$, and
$g\in Mon(X)$ is associated to 
a family $\pi'':\X''\rightarrow B''$ via an isomorphism $X\cong\X''_{b''}$, then $fg$ is easily seen to be associated 
to the family $\pi:\X\rightarrow B$, obtained by ``gluing'' $\X'$ and $\X''$ via 
the isomorphism $\X'_{b'}\cong X\cong \X''_{b''}$ and connecting 
$B'$ and $B''$  at the points $b'$ and $b''$ to form the (reducible) base $B$.
} 
of $GL[H^*(X,\Integers)]$
consisting of all monodromy operators. We denote by $Mon^2(X)$ the image of $Mon(X)$ in
$O[H^2(X,\Integers)]$.
\item
Let $H_i$ be an ample line bundle on $X_i$, $i=1,2$. An isomorphism 
$f:H^2(X_1,\Integers)\rightarrow H^2(X_2,\Integers)$ 
is said to be a {\em polarized parallel-transport operator} from $(X_1,H_1)$ to
$(X_2,H_2)$, if 
there exists a family $\pi:\X\rightarrow B$, satisfying all the properties of part 
(\ref{def-item-parallel-transport-operator}), as well as a flat section $h$ of
$R^2{\pi_*}\Integers$, such that $h(b_i)=\psi_{i_*}(c_1(H_i))$, $i=1,2$,
and $h(b)$ is an ample class in $H^{1,1}(\X_b,\Integers)$, for
all $b\in B$.
\item
Given an ample line bundle $H$ on $X$, we
denote by $Mon(X,H)$ the subgroup of $Mon(X)$, consisting of
polarized parallel transport operators from $(X,H)$ to itself.
Elements of $Mon(X,H)$ will be called 
{\em polarized monodromy operators} of $(X,H)$.
\end{enumerate}
\end{defi}

Following is a necessary condition for an isometry
%When $X$ and $Y$ are $K3$ surfaces and 
$g:H^2(X,\Integers)\rightarrow H^2(Y,\Integers)$ to be a parallel transport operator.
%is an isometry, we can easily determine if $g$ is a parallel transport operator. 
Denote by $\widetilde{\C}_X\subset H^2(X,\RealNumbers)$ the cone
\[
\{\alpha\in H^2(X,\RealNumbers) \ : \ (\alpha,\alpha)>0\}.
\]
Then $H^2(\widetilde{\C}_X,\Integers)\cong\Integers$ and it comes with a canonical generator,
which we call the {\em orientation class} of $\widetilde{\C}_X$ 
(section \ref{sec-orientation}). 
Any isometry $g:H^2(X,\Integers)\rightarrow H^2(Y,\Integers)$ induces 
an isomorphism $\bar{g}:\widetilde{\C}_X\rightarrow \widetilde{\C}_Y$.
The isometry $g$ is said to be {\em orientation preserving} if $\bar{g}$ is.
A parallel transport operator  
$g:H^2(X,\Integers)\rightarrow H^2(Y,\Integers)$ is 
orientation preserving.
When $X$ and $Y$ are $K3$ surfaces, every orientation preserving isometry is 
a parallel transport operator. This is no longer the case for higher dimensional 
irreducible holomorphic symplectic varieties 
\cite{markman-constraints,namikawa-torelli}.
%Parallel transport operators will be discussed in sections 
%\ref{sec-parallel-transport-between-non-separated-marked-pairs} and \ref{sec-K3-type}.
%In particular, 
A necessary and sufficient criterion for an isometry to be
a parallel transport operator is provided in the $K3^{[n]}$-type case, for
all $n\geq 1$ (Theorem 
\ref{thm-necessary-and-sufficient-condition-to-be-a-parallel-transport-operator}).

A {\em marked pair} $(X,\eta)$ consists of an irreducible holomorphic symplectic manifold $X$
and an isometry $\eta:H^2(X,\Integers)\rightarrow \Lambda$ onto a fixed lattice
$\Lambda$. Let $\FM^0_\Lambda$ be a connected component of 
the moduli space of isomorphism classes of marked pairs 
(see section \ref{sec-global-torelli}). 
There exists a surjective period map $P_0:\FM_\Lambda^0\rightarrow \Omega_\Lambda$ onto
a period domain (\cite{huybrechts-basic-results}, Theorem 8.1). 
Each point $\period\in \Omega_\Lambda$ determines a weight $2$ Hodge structure 
on $\Lambda\otimes_\Integers\ComplexNumbers$, such that the marking $\eta$
is an isomorphism of Hodge structures. 
The positive cone $\C_X$ of $X$ is the connected component
of the cone $\{\alpha\in H^{1,1}(X,\RealNumbers) \ : \ (\alpha,\alpha)>0\}$,
containing the K\"{a}hler cone $\K_X$. Following is a concise version of the
{\em Global Torelli Theorem} (\cite{verbitsky}, or Theorem  \ref{thm-global-torelli} below).

\begin{thm} 
If $P_0(X,\eta)=P_0(\widetilde{X},\tilde{\eta})$, then $X$ and $\widetilde{X}$ 
are bimeromorphic.
A pair $(X,\eta)$ is the unique point in a fiber of $P_0$,
if and only if $\K_X=\C_X$.
This is the case, for example, 
if the sublattice $H^{1,1}(X,\Integers)$ is trivial, or of rank $1$, generated by an
element $\lambda$, with $(\lambda,\lambda)\geq 0$. 
\end{thm}

The following theorem combines the Global Torelli Theorem 
with results on the K\"{a}hler cone of 
irreducible holomorphic symplectic manifolds 
\cite{huybrechts-kahler-cone,boucksom-kahler-cone}.

\begin{thm} (A Hodge theoretic Torelli theorem)
\label{thm-Hodge-theoretic-Torelli}
Let $X$ and $Y$ be irreducible holomorphic symplectic manifolds, 
which are deformation equivalent.
\begin{enumerate}
\item 
\label{thm-item-bimeromorphic}
$X$ and $Y$ are bimeromorphic, if and only if there exists 
a parallel transport operator $f:H^2(X,\Integers)\rightarrow H^2(Y,\Integers)$,
which is an isomorphism of integral Hodge structures.
\item
\label{thm-item-isomorphic}
Let $f:H^2(X,\Integers)\rightarrow H^2(Y,\Integers)$ be 
a parallel transport operator,
which is an isomorphism of integral Hodge structures.
There exists an isomorphism $\tilde{f}:X\rightarrow Y$, such that
$f=\tilde{f}_*$, if and only if $f$ maps some K\"{a}hler class on $X$ 
to a K\"{a}hler class on $Y$. 
\end{enumerate}
\end{thm}

The theorem is proven in section \ref{sec-proof-of-hodge-theoretic-torelli}. 
It generalizes the Strong Torelli Theorem of Burns and Rapoport 
\cite{burns-rapoport} or (\cite{looijenga-peters}, Theorem 9.1). 

Given a bimeromorphic map $f:X\rightarrow Y$, of irreducible holomorphic symplectic
manifolds, denote by $f_*:H^2(X,\Integers)\rightarrow H^2(Y,\Integers)$
the homomorphism induced by the closure in $X\times Y$ of the graph of $f$.
The homomorphism $f_*$ is known to be an isometry 
(\cite{ogrady-weight-two}, Proposition 1.6.2).
Set $f^*:=(f^{-1})_*$. 

The {\em birational K\"{a}hler cone} $\BK_X$ of $X$ is the union of the cones $f^*\K_Y$,
as $f$ ranges through all bimeromorphic maps from $X$ to irreducible holomorphic
symplectic manifolds $Y$. Let $Mon^2_{Hdg}(X)$ be the subgroup of $Mon^2(X)$
preserving the Hodge structure. 
Results of Boucksom and Huybrechts, on the K\"{a}hler and birational K\"{a}hler cones, 
are surveyed in section \ref{sec-Kahler-cone}. We use them to define a
chamber decomposition of the positive cone $\C_X$, via $Mon^2_{Hdg}(X)$-translates
of cones of the form $f^*\K_Y$ 
(Lemma \ref{lemma-on-semi-chambers}). These chambers are said to be of
{\em K\"{a}hler type}.

Let $\FM^0_\Lambda$ be a  connected component  of the moduli space
of marked pairs. A detailed form of the Torelli theorem 
%\ref{thm-global-torelli}
provides a description of $\FM^0_\Lambda$ as a moduli space of 
Hodge theoretic data as follows.
A point $\period\in\Omega_\Lambda$ 
determines a Hodge structure on $\Lambda$, and so a real subspace 
$\Lambda^{1,1}(\period,\RealNumbers)$  in 
$\Lambda\otimes_\Integers\RealNumbers$, such that a marking
$\eta$ restricts to an isometry 
$H^{1,1}(X,\RealNumbers)\rightarrow \Lambda^{1,1}(\period,\RealNumbers)$, 
for every pair $(X,\eta)$ in the fiber $P_0^{-1}(\period)$. 
%The Torelli theorem yields
\begin{thm} (Theorem \ref{thm-modular-description-of-the-fiber-of-the-period-map})
The map $(X,\eta)\mapsto \eta(\K_X)$ establishes a 
one-to-one correspondence between 
points $(X,\eta)$ in the fiber $P_0^{-1}(\period)$ and chambers 
in the K\"{a}hler type
chamber decomposition of the positive cone in $\Lambda^{1,1}(\period,\RealNumbers)$.
%The chamber associated to $(X,\eta)$ is the image via $\eta$ of the K\"{a}hler cone of $X$. 
\end{thm}
%**************************************************************************************
%
%**************************************************************************************
\subsection{The fundamental exceptional chamber}
The next few results are easier to understand when compared to the following basic
fact about $K3$ surfaces. Let $S$ be a $K3$ surface and $\kappa_0$ a K\"{a}hler
class on $S$. The effective cone in $H^{1,1}(S,\Integers)$ 
is spanned by classes $\alpha$, such that $(\alpha,\alpha)\geq -2$, and 
$(\alpha,\kappa_0)>0$ (\cite{BHPV}, Ch. VIII Proposition 3.6). 
Set\footnote{$\Pex$ stands for prime exceptional classes, and $\Spe$ stands for stably 
prime exceptional classes, as will be explained below.}
\begin{eqnarray*}
\Spe & := & \{e\in H^{1,1}(S,\Integers) \ : \
(\kappa_0,e)>0, \ \mbox{and} \ (e,e)=-2\},
\\
\Pex & := & \{[C]\in H^{1,1}(S,\Integers) \ : \ C\subset S \ 
\mbox{is a smooth connected rational curve}
\}.
\end{eqnarray*}
Clearly, $\Pex$ is contained in $\Spe$.
Then the K\"{a}hler cone admits the following two characterizations
(\cite{BHPV}, Ch. VIII Proposition 3.7 and Corollary 3.8).
\begin{eqnarray}
\label{eq-Kahler-cone-of-a-K3-via-all-effective-minus-2}
\K_S & = & \{\kappa\in \C_S \ : \ 
(\kappa,e)>0, \ \mbox{for all} \ e\in \Spe\}.
\\
\label{eq-Kahler-cone-of-a-K3-via-prime-exceptional}
\K_S & = & \{\kappa\in \C_S \ : \ (\kappa,e)>0, \ \mbox{for all} \ e\in \Pex\}.
\end{eqnarray}

Equality (\ref{eq-Kahler-cone-of-a-K3-via-all-effective-minus-2}) is the simpler one,
depending only on the Hodge structure and the intersection pairing.
Equality (\ref{eq-Kahler-cone-of-a-K3-via-prime-exceptional})
expresses the fact that a class $e\in\Spe$ represents a smooth rational curve,
if and only if $\overline{\K}_S\cap e^\perp$ is a 
co-dimension one  face of the closure of $\K_S$ in $\C_S$.

Let $X$ be a projective irreducible holomorphic symplectic manifold.
A {\em prime exceptional divisor} on $X$ is a
reduced and irreducible effective divisor $E$ 
of negative Beauville-Bogomolov degree.
The {\em fundamental exceptional chamber} of the positive cone
is the set 
\begin{equation}
\label{eq-fundamental-exceptional-chamber}
\FE_X \ \ := \ \ 
\{\alpha\in\C_X \ : \ (\alpha,[E])>0, \ \mbox{for every prime exceptional divisor} \ E\}.
\end{equation}
When $X$ is a $K3$ surface, a prime exceptional divisor is simply 
a smooth rational curve. Furthermore, the cones $\K_X$, $\BK_X$, and
$\FE_X$ are equal. If $\dim(X)>2$, the cone $\BK_X$ need not be convex.
The following is thus a generalization of equality 
(\ref{eq-Kahler-cone-of-a-K3-via-prime-exceptional}) in the $K3$ surface case.

\begin{thm}
\label{thm-intro-numerical-characterization-of-BK-X}
(Theorem \ref{thm-numerical-characterization-of-BK-X} and Proposition \ref{prop-closure-of-FE-and-FU-are-equal})
$\FE_X$ is an open cone, which is the interior of a closed 
generalized convex polyhedron in $\C_X$
(Definition \ref{def-generalized-convex-polyhedron}).
The birational K\"{a}hler cone 
$\BK_X$ is a dense open subset of $\FE_X$.
% and the closure of both cones in $\C_X$ are equal.
\end{thm}

%The proof of Theorem \ref{thm-intro-numerical-characterization-of-BK-X}
%relies heavily on Boucksom's Zariski decomposition of
%divisor classes (\cite{boucksom-zariski-decomposition}, or Theorem
%\ref{thm-Zariski-decomposition} below).
Let $E$ be a prime exceptional divisor on a projective irreducible holomorphic symplectic manifold $X$.
In section \ref{sec-generators-for-Mon-2-Hdg} 
we recall that the reflection
\[
R_E \ : \ H^2(X,\Integers) \ \  \longrightarrow  \ \  H^2(X,\Integers),
\]
given by $R_E(\alpha):=  \alpha - \frac{2(\alpha,[E])}{([E],[E])}[E]$, 
is an element of $Mon^2_{Hdg}(X)$ 
(\cite{markman-prime-exceptional}, Corollary 3.6,
or Proposition \ref{prop-reflection-by-exceptional-divisor} below). 
Let $W_{Exc}(X)\subset Mon^2_{Hdg}(X)$ be the subgroup 
generated\footnote{Definition \ref{def-W-Exc} of $W_{Exc}$ is different. 
The two definitions will be shown to be equivalent
in Theorem \ref{thm-Mon-2-Hdg-is-a-semi-direct-product}.} 
by the reflections $R_E$, of all prime exceptional divisors in $X$.
In section \ref{sec-Mon-2-Hdg-is-semi-direct-product} we 
prove the following analogue of a well known result for $K3$ surfaces
(\cite{BHPV}, Ch. VIII, Proposition 3.9).

\begin{thm}
\label{thm-decomposition-of-a-Hodge-isometry-as-wg}
%(Theorem \ref{thm-Mon-2-Hdg-is-a-semi-direct-product}).
$W_{Exc}(X)$ is a normal subgroup of $Mon^2_{Hdg}(X)$.
Let $X_1$ and $X_2$ be projective irreducible holomorphic symplectic manifolds
and $f: H^2(X_1,\Integers)\rightarrow H^2(X_2,\Integers)$ a parallel-transport operator,
which preserves the weight $2$ Hodge structure. Then there exists a unique 
element $w\in W_{Exc}(X_2)$ 
% finite set of prime exceptional divisors $E_1$, \dots, $E_k$ in $X_2$ 
and a birational map $g:X_1\rightarrow X_2$, such that 
$f=w\circ g_*$.
The map $g$ is determined uniquely, up to composition with an automorphism
of $X_1$, which acts trivially on $H^2(X_1,\Integers)$.
\end{thm}

Let us emphasis the special case $X_1=X_2=X$ of the theorem. Denote by 
$Mon^2_{Bir}(X)\subset O[H^2(X,\Integers)]$  the subgroup of isometries
induced by birational maps from $X$ to itself.
Then $Mon^2_{Hdg}(X)$ is the semi-direct product of $W_{Exc}(X)$ and $Mon^2_{Bir}(X)$,
 by Theorem \ref{thm-Mon-2-Hdg-is-a-semi-direct-product} part \ref{thm-item-semi-direct-product}.
Theorem \ref{thm-decomposition-of-a-Hodge-isometry-as-wg} is proven in section
\ref{sec-Mon-2-Hdg-is-semi-direct-product}. The proof 
relies on a second $Mon^2_{Hdg}(X)$-equivariant 
chamber decomposition of the positive cone $\C_X$. 
We call these the {\em exceptional chambers} 
(Definition \ref{def-Kahler-type-chambers}). $W_{Exc}(X)$ acts 
simply-transitively on the set of exceptional chambers, one of which is
the fundamental exceptional chamber. 
The walls of a general exceptional chamber are hyperplanes orthogonal 
to classes of {\em stably prime-exceptional} line bundles. 
The latter are higher-dimensional analogues 
of effective line bundles of degree $-2$ on a $K3$ surface. 
Roughly, a line bundle $L$ 
on $X$ is stably prime-exceptional, if a generic small deformation 
$(X',L')$ of $(X,L)$ satisfies $L'\cong\StructureSheaf{X'}(E')$,
for a prime exceptional divisor $E'$ on $X'$ 
(Definition \ref{def-stably-prime-exceptional}).

Let $X$ be a projective irreducible holomorphic symplectic manifold.
Denote by $\Bir(X)$ the group of birational self-maps of $X$. 
The intersection of $\FE_X$ with the subspace $H^{1,1}(X,\Integers)\otimes_\Integers\RealNumbers$ 
is equal to the interior of the movable cone of $X$ (Definition \ref{def-movable-cone} and Lemma 
\ref{lemma-movable-cone-is-a-fundamental-domain}). 
We prove a weak version of Morrison's
movable cone conjecture, about the existence of a rational convex polyhedron, 
which is a fundamental domain for the
action of $\Bir(X)$ on the movable cone (Theorem \ref{thm-movable-conj-conj}). 
We use it to prove the following result.
When $X$ is a $K3$ surface, $\Bir(X)=\Aut(X)$. Hence the following is an analogue of 
a result of Looijenga and Sterk (\cite{sterk}, Proposition 2.6).

\begin{thm}
\label{thm-finite-number-of-Bir-X-orbits}
For every integer $d\neq 0$, the number of 
$Bir(X)$-orbits of complete linear systems, which contain an 
irreducible divisor of Beauville-Bogomolov degree $d$, is finite. 
For every positive integer $k$ 
there is only a finite number 
of $Bir(X)$-orbits of complete linear systems, which contain some  
irreducible divisor $D$ of Beauville-Bogomolov degree zero, such that the class $[D]$ is $k$ times 
a primitive class in $H^2(X,\Integers)$.
\end{thm}

Theorem \ref{thm-finite-number-of-Bir-X-orbits} is proven in section \ref{sec-finiteness}.
The proof follows an argument of  Looijenga and Sterk, adapted via 
an analogy between results on the ample cone of a projective $K3$ surface
and results on the movable cone of a projective  irreducible holomorphic symplectic manifold.

The following is an analogue of the characterization
of the K\"{a}hler cone of a $K3$ surface given in equation \nolinebreak
(\ref{eq-Kahler-cone-of-a-K3-via-all-effective-minus-2}).

\begin{prop}
\label{prop-intro-fundamental-exceptional-chamber-in-indeed-such}
(Proposition \ref{prop-fundamental-exceptional-chamber-is-indeed-an-exceptional-one})
The fundamental exceptional chamber $\FE_X$, defined in equation
(\ref{eq-fundamental-exceptional-chamber}), is equal to the set
\[
\{\alpha\in \C_X \ : \ (\alpha,\ell)>0, \ \mbox{for every stably prime exceptional
class} \ \ell\}.
\]
\end{prop}
The significance of
Proposition \ref{prop-intro-fundamental-exceptional-chamber-in-indeed-such}
stems from the fact that one has an explicit characterization
of the set of stably prime-exceptional classes, in terms of the weight $2$ Hodge structure and a
certain discrete monodromy invariant, 
at least in the $K3^{[n]}$-type case (Theorem 
\ref{thm-stably-prime-exceptional-K3-n}). 
Theorem \ref{thm-intro-numerical-characterization-of-BK-X} and 
Proposition \ref{prop-intro-fundamental-exceptional-chamber-in-indeed-such} thus 
yield an explicit description of the closure of the birational K\"{a}hler cone and of the movable cone.

%**************************************************************************************
%
%**************************************************************************************
\subsection{Torelli and monodromy in the polarized case}
In sections \ref{sec-polarized-monodromy} and
\ref{sec-monodromy-quotients-of-period-domains} we consider  
Torelli-type results for polarized irreducible holomorphic 
symplectic manifolds.
Another corollary of the Global Torelli Theorem is the following.

\begin{prop}
\label{prop-intro-Mon-X-H-is-a-stabilizer}
$Mon^2(X,H)$ is equal to the stabilizer of $c_1(H)$ in $Mon^2(X)$. 
\end{prop}

The above proposition is proven in section \ref{sec-polarized-monodromy}
(see Corollary \ref{cor-Mon-X-H-is-a-stabilizer}).

Coarse moduli spaces of polarized projective irreducible holomorphic symplectic
manifolds were constructed by Viehweg as quasi-projective  varieties \cite{viehweg}.
Given a polarized pair $(X,H)$ representing a point in such a coarse moduli space $\V$, 
the monodromy group $\Gamma:=Mon^2(X,H)$ is an arithmetic group, which acts on a
period domain $\D$ associated to $\V$. The 
quotient $\D/\Gamma$ is a quasi-projective variety \cite{baily-borel}. The following Theorem is 
a slight sharpening of Corollary 1.24 in \cite{verbitsky}.

\begin{thm}
(Theorem \ref{thm-viehweg-moduli-is-open-in-monodromy-quotient-of-period-space})
The period map $\V\rightarrow \D/\Gamma$ 
embeds each irreducible component $\V$, of the coarse moduli space 
of polarized irreducible holomorphic symplectic manifolds, 
as a Zariski  open subset of the quasi-projective
monodromy-quotient of the corresponding period domain.
\end{thm}

The above theorem provides a bridge between the powerful theory of modular forms, used to
study the quotient spaces $\D/\Gamma$, and the theory of projective holomorphic symplectic varieties.
The interested reader is referred to the excellent recent survey \cite{GHS-survey} for further reading on this topic.
%**************************************************************************************
%
%**************************************************************************************
\subsection{The $K3^{[n]}$-type}
In section \ref{sec-K3-type} we specialize to the case of 
varieties $X$ of $K3^{[n]}$-type and review the results of 
\cite{markman-monodromy-I,markman-constraints,markman-prime-exceptional}.
We introduce a Hodge theoretic {\em Torelli data}, consisting of 
the weight $2$ Hodge structure of $X$ and a certain discrete monodromy invariant
(Corollary \ref{cor-third-characterization-of-Mon-2}). 
We provide explicit computations, for many of the concepts introduced above, in terms of 
this Torelli data.
We enumerate the connected components of the moduli space of marked pairs of 
$K3^{[n]}$-deformation type (Corollary \ref{cor-enumeration-of-connected-components}).
We determine the monodromy group $Mon^2(X)$, 
as well as a necessary and sufficient condition for an isometry
$g:H^2(X,\Integers)\rightarrow H^2(Y,\Integers)$ to be a parallel transport operator 
(Theorems \ref{thm-mon-2-is-a-reflection-group} and  
\ref{thm-necessary-and-sufficient-condition-to-be-a-parallel-transport-operator}). 
We provide a numerical characterization of the set of 
stably prime-exceptional line bundles on $X$ 
(Theorem \ref{thm-stably-prime-exceptional-K3-n}). 
The latter, combined with the general 
Theorem \ref{thm-intro-numerical-characterization-of-BK-X}
and Proposition \ref{prop-intro-fundamental-exceptional-chamber-in-indeed-such},
determines the closure of the birational K\"{a}hler cone of $X$ in terms
of its Torelli data. 

In section \ref{sec-open-problems} we list a few open problems.

{\bf Acknowledgements:} I would like thank Klaus Hulek for encouraging  
me to write this survey, and for many insightful 
discussions and suggestions.
This note was greatly influenced by numerous conversations with Daniel Huybrechts
and by his foundational written  work. Significant improvements to an earlier version of this survey
are due to Daniel's detailed comments and suggestions, for which I am most grateful.
The note is the outcome of an extensive correspondence with 
Misha Verbitsky regarding his fundamental paper \cite{verbitsky}.
I am most grateful for his patience and for his numerous detailed answers.
Artie Prendergast-Smith kindly sent helpful comments to an earlier version of section
\ref{sec-finiteness}, for which I am grateful.
I would like to thank the two referees for their careful reading and their insightful comments. 
%****************************************************************
% 
%****************************************************************
\section{The Global Torelli Theorem}
\label{sec-global-torelli}
Fix a positive integer $b_2>3$ and an even lattice 
$\Lambda$ of signature $(3,b_2-\nolinebreak 3)$. 
%Assume, for simplicity, that $b_2\neq 6$.
Let $X$ be an irreducible holomorphic symplectic manifold, such that
$H^2(X,\Integers)$, endowed with its Beauville-Bogomolov pairing,
is isometric to $\Lambda$. A {\em marking} for $X$ is a choice of an isometry
$\eta:H^2(X,\Integers)\rightarrow \Lambda$. Two marked pairs
$(X_1,\eta_1)$, $(X_2,\eta_2)$ are isomorphic, if there exists 
an isomorphism $f:X_1\rightarrow X_2$, such that 
$\eta_1\circ f^*=\eta_2$.
There exists a coarse moduli space $\FM_\Lambda$ parametrizing
isomorphism classes of marked pairs \cite{huybrechts-basic-results}.
$\FM_\Lambda$ is a smooth complex manifold of dimension $b_2-2$,
but it is non-Hausdorff. 

The {\em period}, of the marked pair $(X,\eta)$, is the
line $\eta[H^{2,0}(X)]$ considered as a point in the projective space
$\PP[\Lambda\otimes_\Integers\ComplexNumbers]$. The period lies in the period domain 
\begin{equation}
\label{eq-period-domain}
\Omega_\Lambda \ := \ \{
\period \ : \ (\period,\period)=0 \ \ \ \mbox{and} \ \ \ 
(\period,\bar{\period}) > 0
\}.
\end{equation}
$\Omega_\Lambda$ is an open subset, in the classical topology, of the quadric in 
$\PP[\Lambda\otimes\ComplexNumbers]$ of isotropic lines \cite{beauville}. 
The period map 
\begin{eqnarray}
\label{eq-period-map}
P \ : \ {\mathfrak M}_\Lambda & \longrightarrow & \Omega_\Lambda,
\\
\nonumber
(X,\eta) & \mapsto & \eta[H^{2,0}(X)]
\end{eqnarray}
is a local isomorphism, by the Local Torelli Theorem \cite{beauville}. 

Given a point $\period\in \Omega_\Lambda$, set 
$\Lambda^{1,1}(\period):=\{\lambda\in \Lambda \ : \ (\lambda,\period)=0\}$. Note that 
$\Lambda^{1,1}(\period)$ is a sublattice of $\Lambda$ and 
$\Lambda^{1,1}(\period)=(0)$, if $\period$ does not belong to 
the countable union of hyperplane sections
$\cup_{\lambda\in \Lambda\setminus\{0\}}[\lambda^\perp\cap\Omega_\Lambda]$.
Given a marked pair $(X,\eta)$, we get the isomorphism
$H^{1,1}(X,\Integers)\cong\Lambda^{1,1}(P(X,\eta))$, via the restriction of $\eta$.

\begin{defi}
\label{def-positive-cone}
Let $X$ be an irreducible holomorphic symplectic manifold.
The cone 
$
\{\alpha\in H^{1,1}(X,\RealNumbers) \ : \ (\alpha,\alpha)>0\}
$
has two connected components. 
The {\em positive cone} $\C_X$
is the connected component containing the K\"{a}hler cone 
$\K_X$.
\end{defi}

%The two cones $\K_X$ and $\C_X$ are equal, when 
%$H^{1,1}(X,\Integers)$ is trivial, or if $H^{1,1}(X,\Integers)$ has rank $1$, 
%spanned by a class $\alpha$ with non-negative 
%Beauville-Bogomolov degree, by
%(\cite{huybrechts-basic-results}, Corollaries 5.7 and 7.2).

Two points $x$ and $y$ of a topological space $M$ are 
{\em inseparable}, if every pair of open subsets $U$, $V$, 
with $x\in U$ and $y\in V$, have a non-empty intersection $U\cap V$.
A point $x\in M$ is a {\em Hausdorff point}, 
if there does not exist any point $y\in [M\setminus\{x\}]$,
such that $x$ and $y$ are inseparable.

\begin{thm} 
\label{thm-global-torelli}
(The Global Torelli Theorem)
Fix a connected component $\FM^0_\Lambda$ of $\FM_\Lambda$.
\begin{enumerate}
\item
(\cite{huybrechts-basic-results}, Theorem 8.1)
The period map $P$ restricts to a surjective holomorphic map
$P_0:\FM^0_\Lambda\rightarrow \Omega_\Lambda$. 
\item 
\label{thm-item-injectivity}
(\cite{verbitsky}, Theorem 1.16)
The fiber $P_0^{-1}(\period)$ consists of pairwise 
inseparable points, for all $\period\in \Omega_\Lambda$.
\item
\label{thm-item-inseparable-are-bimeromorphic}
(\cite{huybrechts-basic-results}, Theorem 4.3)
Let  $(X_1,\eta_1)$ and $(X_2,\eta_2)$ be two inseparable points of $\FM_\Lambda$.
Then $X_1$ and $X_2$ are bimeromorphic.
\item
\label{thm-item-single-Hausdorff-point}
The marked pair $(X,\eta)$ is a Hausdorff point of $\FM_\Lambda$,
if and only if $\C_X=\K_X$.
%The fiber $P_0^{-1}\left[P_0(X,\eta)\right]$ consists of a single Hausdorff point, 
%if and only if $\C_X=\K_X$.
\item
\label{thm-item-cyclic-picard-and-projective-imply-K-X=C-X}
The fiber $P_0^{-1}(\period)$, $\period\in \Omega_\Lambda$, consists of a single Hausdorff point,
if $\Lambda^{1,1}(\period)$ is trivial, or if $\Lambda^{1,1}(\period)$ 
is of rank $1$, generated by a class $\alpha$ satisfying $(\alpha,\alpha)\geq 0$.
%$\C_X=\K_X$, if $H^{1,1}(X,\Integers)$ is trivial, or if $H^{1,1}(X,\Integers)$ 
%is of rank $1$, generated by a class $\alpha$ of
%non-negative  Beauville-Bogomolov degree.
%\item
%\label{thm-item-Hausdorff-points}
%The fiber $P_0^{-1}(\period)$ consists of a single Hausdorff point $(X,\eta)$,
%if and only if the K\"{a}hler cone $\K_X$ is equal to the positive cone $\C_X$.
%In particular, $P_0^{-1}(\period)$ consists of a single Hausdorff point, 
%if $\Lambda^{1,1}(\period)=(0)$, 
%or if $\Lambda^{1,1}(\period)$ is spanned by a class of positive degree.
\end{enumerate}
\end{thm}

\begin{proof}
Part (\ref{thm-item-single-Hausdorff-point}) of the theorem is due to Huybrechts and
Verbitsky. See Proposition \ref{prop-fiber} for a more general 
description of the fiber $P_0^{-1}\left[P_0(X,\eta)\right]$ in terms of 
the K\"{a}hler-type chamber
decomposition of the positive cone $\C_X$, and for 
further details about part (\ref{thm-item-single-Hausdorff-point}).

Part (\ref{thm-item-cyclic-picard-and-projective-imply-K-X=C-X}):
$\C_X=\K_X$, if $H^{1,1}(X,\Integers)$ is trivial, or if $H^{1,1}(X,\Integers)$ 
is of rank $1$, generated by a class $\alpha$ of
non-negative  Beauville-Bogomolov degree, by 
(\cite{huybrechts-basic-results}, Corollaries 5.7 and 7.2).
The statement of part (\ref{thm-item-cyclic-picard-and-projective-imply-K-X=C-X})
now follows from part (\ref{thm-item-single-Hausdorff-point}).
\end{proof}

\begin{rem}
Verbitsky states part (\ref{thm-item-injectivity}) of Theorem \ref{thm-global-torelli}
for a connected component of the Teichm\"{u}ller space, but  Theorem 1.16
in \cite{verbitsky} is a consequence of the two more general Theorems
4.22 and 6.14 in \cite{verbitsky}, and both the Teichm\"{u}ller space and
the moduli space of marked pairs $\FM_\Lambda$ satisfy the hypothesis of these theorems.
A complete proof of part (\ref{thm-item-injectivity}) of Theorem \ref{thm-global-torelli}
can be found in Huybrechts excellent Bourbaki seminar paper \cite{huybrechts-bourbaki}.
\end{rem}
%****************************************************************
% 
%****************************************************************
\section{The Hodge theoretic Torelli Theorem}
In section \ref{sec-parallel-transport-between-non-separated-marked-pairs} we 
review two theorems of Huybrechts, which relate bimeromorphic maps and 
parallel-transport operators. 
The Hodge theoretic Torelli Theorem 
\ref{thm-Hodge-theoretic-Torelli} is proven in section \ref{sec-proof-of-hodge-theoretic-torelli}.

%****************************************************************
% 
%****************************************************************
\subsection{Parallel transport operators between inseparable marked pairs}
\label{sec-parallel-transport-between-non-separated-marked-pairs}
Let $X_1$ and $X_2$ be two irreducible holomorphic symplectic manifolds
of dimension $2n$. Denote by $\pi_i$ the projection from $X_1\times X_2$ 
onto $X_i$, $i=1,2$. 
Given a correspondence $Z$ in $X_1\times X_2$, of pure complex co-dimension $2n+d$, 
denote by $[Z]$ the cohomology class
Poincar\'{e} dual to $Z$ and by $[Z]_*:H^*(X_1)\rightarrow H^{*+2d}(X_2)$
the homomorphism defined by $[Z]_*\alpha:=\pi_{2_*}\left(\pi_1^*(\alpha)\cup [Z]\right)$.
The following are two fundamental results of Huybrechts.

Assume that $X_1$ and $X_2$ are bimeromorphic.
Denote the graph of a bimeromorphic map by 
$Z\subset X_1\times X_2$. 

\begin{thm}
\label{thm-bimeromorphic-implies-inseparable}
(\cite{huybrechts-kahler-cone}, Corollary 2.7)
There exists an effective cycle 
$\Gamma:=Z+\sum Y_j$ in $X_1\times X_2$, 
of pure dimension $2n$, with the following properties.
\begin{enumerate}
\item
The correspondence $[\Gamma]_*:H^*(X_1,\Integers)\rightarrow H^*(X_2,\Integers)$
is a parallel-transport operator.
\item
The image $\pi_i(Y_j)$ has codimension $\geq 2$ in $X_i$, for all $j$.
In particular, the correspondences $[\Gamma]_*$ and $[Z]_*$ 
coincide on $H^2(X_1,\Integers)$.
\end{enumerate}
\end{thm}

Let $(X_1,\eta_1)$, $(X_2,\eta_2)$ be two marked pairs corresponding to 
inseparable points of ${\mathfrak M}_\Lambda$. 

\begin{thm} 
\label{thm-non-separated-implies-bimeromorphic}
(\cite{huybrechts-basic-results}, Theorem 4.3 and its proof)
There exists an effective cycle
$\Gamma:= Z + \sum_j Y_j$ in $X_1\times X_2$,
of pure dimension $2n$, 
satisfying the following conditions.
\begin{enumerate}
\item 
$Z$ is the graph of a bimeromorphic map from $X_1$ to $X_2$.
\item
\label{thm-item-composition-of-markings-is-a-parallel-transport-op}
The correspondence $[\Gamma]_*:H^*(X_1,\Integers)\rightarrow H^*(X_2,\Integers)$
is a parallel-transport operator. Furthermore, 
the composition 
\[
\eta_2^{-1}\circ\eta_1:H^2(X_1,\Integers)\rightarrow H^2(X_2,\Integers)
\]
is equal to the restriction of $[\Gamma]_*$.
\item
(\cite{huybrechts-kahler-cone}, Theorem 2.5 and its proof)
The codimensions of $\pi_1(Y_j)$ in $X_1$ and of
$\pi_2(Y_j)$ in $X_2$ are equal and positive.
\item
\label{thm-item-uniruled-divisor}
If $\pi_i(Y_j)$ has codimension $1$, then it is supported by a uniruled divisor.
\end{enumerate}
\end{thm}

%\begin{question}
%Can one further show that the divisors appearing in part \ref{thm-item-uniruled-divisor}
%of Theorem \ref{thm-non-separated-implies-bimeromorphic}
%all have negative Beauville-Bogomolov degree?
%\end{question}

The statement that the isomorphisms $[\Gamma]_*$ in Theorems
\ref{thm-bimeromorphic-implies-inseparable} and \ref{thm-non-separated-implies-bimeromorphic}
are parallel transport operators is implicit in Huybrechts proofs, so we clarify that point next.
In each of the proofs Huybrechts shows that there exist two smooth and proper families
$\X\rightarrow B$ and $\X'\rightarrow B$, over the same one-dimensional disk $B$, a point $b_0$ in $B$,
isomorphisms $X_1\cong \X_{b_0}$ and $X_2\cong \X'_{b_0}$, and an isomorphism
$\tilde{f}:\restricted{\X}{B\setminus\{b_0\}}\rightarrow \restricted{\X'}{B\setminus\{b_0\}}$, compatible with projections to 
$B$.
The cycle $\Gamma\subset X_1\times X_2$ is the fiber over $b_0$ of the closure in $\X\times_B\X'$
of the graph of $\tilde{f}$. Choose a point $b_1$ in $B\setminus \{b_0\}$ and let 
$\gamma$ be a continuous path in $B$ from $b_0$ to $b_1$. Let 
$g_1:H^*(\X_{b_0},\Integers)\rightarrow H^*(\X_{b_1},\Integers)$ and
$g_2:H^*(\X'_{b_0},\Integers)\rightarrow H^*(\X'_{b_1},\Integers)$
be the two parallel transport operators along $\gamma$. Then 
the isomorphism $g_2^{-1}\circ g_1: H^*(\X_{b_0},\Integers)\rightarrow H^*(\X'_{b_0},\Integers)$ 
is induced by the correspondence $[\Gamma]_*$. Furthermore, $g_2^{-1}\circ g_1$ is a parallel transport operator,
being a composition of such operators (parallet transport operators form a groupoid, by
an argument similar to that used in
footnote \footnoterecall{group-footnote}). 

The reader may wonder why the image in $X_i$ of a component $Y_j$ of $\Gamma$ has codimension
 $\geq 2$ in Theorem \ref{thm-bimeromorphic-implies-inseparable}, while 
 the codimension is only $\geq 1$ in 
Theorem \ref{thm-non-separated-implies-bimeromorphic}.
The reason is that in the proof of Theorem \ref{thm-non-separated-implies-bimeromorphic} 
one does not have control on the choice of the above mentioned families $\X$ and $\X'$,
beyond the condition that $\eta_2\circ[\Gamma]_*=\eta_1$.
In the proof of Theorem \ref{thm-bimeromorphic-implies-inseparable},
given a bimeromorphic map $f:X_1\rightarrow X_2$, Huybrechts constructs the above two families 
$\X$ and $\X'$ in such a way that the following two
properties hold. (1) The bimeromorphic map
$\tilde{f}$ from $\X$ to $\X'$ restricts to the
bimeromorphic map $f$ between the  fibers $X_1$ and $X_2$ over $b_0$.
(2) $[\Gamma]_*$  restricts to the isometry  $f_*:H^2(X_1,\Integers)\rightarrow H^2(X_2,\Integers)$
(see Theorem 2.5 in \cite{huybrechts-kahler-cone} and its proof).

%****************************************************************
% 
%****************************************************************
\subsection{Proof of the Hodge theoretic Torelli Theorem \ref{thm-Hodge-theoretic-Torelli}}
\label{sec-proof-of-hodge-theoretic-torelli}
\hspace{1ex}

{\bf Proof of part \ref{thm-item-bimeromorphic}:}
If $X$ and $Y$ are bimeromorphic, then there exists a parallel-transport operator
$f:H^2(X,\Integers)\rightarrow H^2(Y,\Integers)$, which is an isomorphism
of Hodge structures, by Theorem \ref{thm-bimeromorphic-implies-inseparable}. 
Conversely, assume that such $f$ is given.
Let $\eta_Y:H^2(Y,\Integers)\rightarrow \Lambda$ be a marking.
Set $\eta_X:=\eta_Y\circ f$. The assumption that $f$ is a parallel transport operator 
implies that $(X,\eta_X)$ and $(Y,\eta_Y)$ belong to the same connected component 
$\FM_\Lambda^0$ of $\FM_\Lambda$. 
Both have the same period
\[
P(X,\eta_X)=\eta_X(H^{2,0}(X))=\eta_Y(f(H^{2,0}(X)))=\eta_Y(H^{2,0}(Y))=P(Y,\eta_Y),
\]
where the third equality follows from the assumption that $f$ is an isomorphism of 
Hodge structures.
Hence, $(X,\eta_X)$ and $(Y,\eta_Y)$  are inseparable points of $\FM_\Lambda^0$,
by Theorem \ref{thm-global-torelli} part
\ref{thm-item-injectivity}. 
$X$ and $Y$ are thus bimeromorphic, by
Theorem \ref{thm-global-torelli} part
\ref{thm-item-inseparable-are-bimeromorphic}.

{\bf Proof of part \ref{thm-item-isomorphic}:}
Let $\eta_X$ and $\eta_Y$ be the markings constructed in the 
proof of part \ref{thm-item-bimeromorphic}. Note that
$f=\eta_Y^{-1}\circ \eta_X$. 
There exists an {\em effective} correspondence
$\Gamma=Z+\sum_{i=1}^N W_i$ of pure dimension $2n$ in $X\times Y$,
such that $Z$ is the graph of a bimeromorphic map,
$W_i$ is irreducible, but not necessarily reduced, 
the images of the projections $W_i\rightarrow X$,
$W_i\rightarrow Y$ have positive co-dimensions, and 
$[\Gamma]_*:H^*(X,\Integers)\rightarrow H^*(Y,\Integers)$
is a parallel transport operator, which is equal to $f$ in degree $2$, by 
Theorem \ref{thm-non-separated-implies-bimeromorphic}
%Such a correspondence $\Gamma$ is constructed in the proof of 
%(\cite{huybrechts-basic-results}, Theorem 4.3), 
and the assumption that
the two points $(X,\eta_X)$ and $(Y,\eta_Y)$ are inseparable.

Assume that $\alpha\in \K_X$ is 
a K\"{a}hler class, such that $f(\alpha)$ is a
K\"{a}hler class. 
%Such a class $\alpha$ exists, by assumption.
The relationship between $f$ and $\Gamma$ yields:
\[
f(\alpha) \ = \ [\Gamma]_*(\alpha) \ = \ 
[Z]_*(\alpha) + \sum_{i=1}^N [W_i]_*(\alpha).
\]
%The class $[Z]_*(\alpha)$ belongs to the fundamental uniruled chamber of $Y$,
%by Theorem \ref{thm-kahler-cone}.
Each class $[W_i]_*(\alpha)$ is either zero or a multiple $c_i[D_i]$ 
of the class of a prime divisor $D_i$, where $c_i$ is a positive\footnote{The 
coefficient $c_i$ is positive since $\Gamma$ is effective and 
$\alpha$ is a K\"{a}hler class.} 
real number.

We prove next that $[W_i]_*(\alpha)=0$, for $1\leq i\leq N$.
Write $f(\alpha)=[Z]_*(\alpha)+\sum_{i=1}^Nc_i[D_i]$, where $c_i$ are all
positive real numbers, and $D_i$ is either a prime divisor, or zero.
Set $D:=\sum_{i=1}^Nc_iD_i$. 
We need to show that all $D_i$ are equal to zero. The Beauville-Bogomolov degree of $\alpha$ 
satisfies
\[
(\alpha,\alpha)=(f(\alpha),f(\alpha))=
([Z]_*\alpha,[Z]_*\alpha)+2\sum_{i=1}^Nc_i([Z]_*\alpha,[D_i])+([D],[D]).
\]
The homomorphism $[Z]_*$, induced by the graph of the bimeromorphic map, 
is an isometry, by \cite{ogrady-weight-two}, Proposition 1.6.2 (also by 
the stronger Theorem
\ref{thm-bimeromorphic-implies-inseparable}).  
Furthermore, if $D_i$ is non-zero, then $D_i$ is the strict transform of 
a prime divisor $D'_i$ on $X$, such that $[Z]_*([D'_i])=[D_i]$. 
Set $D':=\sum_{i=1}^Nc_iD'_i$. 
We get the equalities
\begin{eqnarray}
\label{eq-D-D}
([D],[D]) & = &  -2(\alpha,[D']),
%\sum_{i=1}^Nc_i(\alpha,[D_i']),
\\
\label{eq-D-prime}
\hspace{0ex}[D]&=&[Z]_*[D'],
\end{eqnarray}
and 
\[
%0\leq 
([D],f(\alpha))=([D],[Z]_*\alpha)+([D],[D])
\stackrel{(\ref{eq-D-prime})}{=}
([D'],\alpha)+([D],[D])
\stackrel{(\ref{eq-D-D})}{=}
-(\alpha,[D']).
%-\sum_{i=1}^Nc_i(\alpha,[D_i'])
%\leq 0.
\]
Now $(\alpha,[D_i'])$ is zero, if $D_i=0$, and positive, if $D_i\neq 0$, since $\alpha$
is a K\"{a}hler class. Hence, the right hand side above is $\leq 0$.
The left hand side is $\geq 0$, due to the assumption that 
the class $f(\alpha)$ is a K\"{a}hler class.
%The right inequality is due to the fact that $c_i$ are all positive, $\alpha$ is a K\"{a}hler class,
%and $D_i'$ are all effective.
Hence, $D'_i=0$, for all $i$. We conclude  that
$[W_i]_*(\alpha)=0$, for $1\leq i\leq N$, as claimed.
 
The equality $[Z]_*(\alpha)=f(\alpha)$ was proven above. Consequently, 
$Z$ is the graph of a bimeromorphic map, which maps a 
K\"{a}hler class to a K\"{a}hler class. Hence, $Z$ is the graph of 
an isomorphism, by a theorem of Fujiki \cite{fujiki}.
\EndProof

%****************************************************************
% 
%****************************************************************
\section{Orientation}
\label{sec-orientation}
Let $\Omega_\Lambda$ be the period domain (\ref{eq-period-domain}).
Following are two examples, in
which spaces arise with two connected components.
\begin{enumerate}
\item
Fix a primitive class $h\in \Lambda$, with $(h,h)>0$. The hyperplane section
\[
\Omega_{h^\perp} \ \ := \ \ \Omega_\Lambda\cap h^\perp
\]
has two connected components. 
\item
Let $\period\in \Omega_\Lambda$.
Set $\Lambda_\RealNumbers:=\Lambda\otimes_\Integers\RealNumbers$ and
$\Lambda^{1,1}(\period,\RealNumbers):=
\{\lambda\in \Lambda_\RealNumbers : (\lambda,\period)=0\}$.
Then the cone
$\C'_\period:=\{\lambda\in \Lambda^{1,1}(\period,\RealNumbers) \ : \ (\lambda,\lambda)>0\}$
has two connected components. 
\end{enumerate}

We recall in this section that a connected component $\FM^0_\Lambda$,
of the moduli space of marked pairs, determines a choice of a 
component of $\Omega_{h^\perp}$ and of $\C'_\period$, for all 
$h\in \Lambda$, with $(h,h)>0$, and for all $\period \in \Omega_\Lambda$. Let us first
relate the choice of one
of the two components in the two examples above. The relation
can be explained in terms of the following larger cone.
Set 
%$\Lambda_\RealNumbers:=\Lambda\otimes_\Integers\RealNumbers$ and
\[
\widetilde{\C}_\Lambda:=\{\lambda\in  \Lambda_\RealNumbers \ : \ (\lambda,\lambda)>0\}.
\]
A subspace $W\subset \Lambda_\RealNumbers$ is said to be {\em positive},
if the pairing of $\Lambda_\RealNumbers$ restricts to $W$ as a positive definite pairing.

\begin{lem}
\begin{enumerate}
\item
\label{lemma-item-orientation-character}
$H^2(\widetilde{\C}_\Lambda,\Integers)$
is a free abelian group of  rank $1$. 
\item
\label{lemma-item-orientation-character-is-spinnor-norm}
Let $e\in \Lambda$ be an element with $(e,e)\neq 0$ and 
$R_e:\Lambda_\RealNumbers\rightarrow \Lambda_\RealNumbers$ the reflection given by
$R_e(\lambda)=\lambda-\frac{2(e,\lambda)}{(e,e)}e$. Then $R_e$ acts on $H^2(\widetilde{\C}_\Lambda,\Integers)$
by $-1$, if $(e,e)>0$, and trivially if $(e,e)<0$.
\item
\label{lemma-item-deformation-retract}
Let $W$ be a positive three dimensional subspace of $\Lambda_\RealNumbers$.
Then $W\setminus \{0\}$ is a deformation retract of $\widetilde{\C}_\Lambda$. 
\end{enumerate}
\end{lem}

\begin{proof}
(\ref{lemma-item-deformation-retract})
Set $I:=[0,1]$. 
We need to construct a continuous map 
$F:\widetilde{\C}_\Lambda\times I\rightarrow \widetilde{\C}_\Lambda$
satisfying 
\begin{eqnarray*}
F(\lambda,0) & = & \lambda, \ \ \ \ \ \ \ \ \ \ \ \mbox{for all} \ \lambda \in \widetilde{\C}_\Lambda,
\\
F(\lambda,1) & \in & W\setminus \{0\}, \ \ \mbox{for all} \ \lambda \in \widetilde{\C}_\Lambda,
\\
F(w,t) & = & w,   \ \ \ \ \ \ \ \ \  \ \mbox{for all} \ w\in W\setminus\{0\}.
\end{eqnarray*}
Choose a basis $\{e_1, e_2, e_3, \dots, e_{b_2}\}$ of $\Lambda_\RealNumbers$, so that 
$\{e_1, e_2, e_3\}$ is a basis of $W$, and for $\lambda=\sum_{i=1}^{b_2}x_ie_i$, we have
$(\lambda,\lambda)=x_1^2+x_2^2+x_3^2-\sum_{i=4}^{b_2}x_i^2$.
Then $\widetilde{\C}_\Lambda$ consists of $\lambda$ satisfying
$
x_1^2+x_2^2+x_3^2>\sum_{i=4}^{b_2}x_i^2.
$
Set $F\left(\sum_{i=1}^{b_2} x_ie_i,t\right)  =  \sum_{i=1}^{3} x_ie_i+(1-t)\sum_{i=4}^{b_2}x_ie_i$.
Then $F$ has the above properties of a deformation retract of $\widetilde{\C}_\Lambda$ onto $W\setminus\{0\}$.

Part (\ref{lemma-item-orientation-character}) follows immediately from part 
(\ref{lemma-item-deformation-retract}).

(\ref{lemma-item-orientation-character-is-spinnor-norm})
If $(e,e)>0$, we can choose a positive $3$ dimensional subspace $W$ containing $e$, and
if $(e,e)<0$ we can choose $W$ to be orthogonal to $e$. 
Then $W\setminus\{0\}$ is $R_e$ invariant and $R_e$ acts as stated on
$H^2(W\setminus\{0\},\Integers)$, hence also on $H^2(\widetilde{\C}_\Lambda,\Integers)$,
by part (\ref{lemma-item-deformation-retract}).
\end{proof}

%There is a canonical isomorphism 
%$H^2(W\setminus\{0\},\Integers)\cong H^2(\widetilde{\C}_\Lambda,\Integers)$,
%by Lemma \ref{lemma-deformation-retract}.
%We conclude that $H^2(\widetilde{\C}_\Lambda,\Integers)$
%is a free abelian group of  rank $1$. 
The character $H^2(\widetilde{\C}_\Lambda,\Integers)$
of $O(\Lambda)$ is known as the {\em spinor norm}.

A point $\period\in \Omega_{h^\perp}$ determines the 
three dimensional positive definite subspace
$W_\period:=\mbox{Re}(\period)\oplus\mbox{Im}(\period)\oplus \mbox{span}_\RealNumbers\{h\}$,
which comes with an orientation associated to the basis 
$\{\mbox{Re}(\sigma),\mbox{Im}(\sigma),h\}$, for some choice of a non-zero 
element $\sigma\in \period\subset \Lambda_\ComplexNumbers$. 
The orientation of the basis is independent of the choice of $\sigma$. 
Consequently, an element $\period\in \Omega_{h^\perp}$ determines
a generator of $H^2(\widetilde{\C}_\Lambda,\Integers)$.
The two components of $\Omega_{h^\perp}$ are distinguished by the two generators of
the rank $1$ free abelian group $H^2(\widetilde{\C}_\Lambda,\Integers)$.
We refer to each of the two generators as an {\em orientation class} of
the cone $\widetilde{\C}_\Lambda$.

A point $\lambda\in \C'_\period$ determines an orientation of $\widetilde{\C}_\Lambda$
as follows. Choose a class $\sigma\in \period$. Again we get the 
three dimensional positive definite subspace
$W_\lambda:=
\mbox{Re}(\period)\oplus\mbox{Im}(\period)\oplus \mbox{span}_\RealNumbers\{\lambda\}$,
which comes with an orientation associated to the basis 
$\{\mbox{Re}(\sigma),\mbox{Im}(\sigma),\lambda\}$. 
Consequently, $\lambda$ determines an orientation 
of $\widetilde{\C}_\Lambda$. The orientation remains the same as $\lambda$ varies in a
connected component of $\C'_\period$. Hence, a connected component of 
$\C'_\period$ determines an orientation of $\widetilde{\C}_\Lambda$.

Let $X$ be an irreducible holomorphic symplectic manifold. 
Recall that the positive cone
$\C_X\subset H^{1,1}(X,\RealNumbers)$ is the distinguished 
connected component of the cone
$\C'_X:=\{\lambda\in H^{1,1}(X,\RealNumbers) \ : \ (\lambda,\lambda)>0\}$,
which contains the K\"{a}hler cone 
(Definition \ref{def-positive-cone}). Denote by 
$\widetilde{\C}_X$ the positive cone in $H^2(X,\RealNumbers)$. 
We conclude that $\widetilde{\C}_X$ comes with a distinguished orientation.
%orientation class, which is determined as follows.
%Let $\sigma$ be a non-zero holomorphic $2$-form and $\alpha$ a 
%K\"{a}hler class or, more generally, a class in $\C_X$. 
%Then $\{\mbox{Re}(\sigma),\mbox{Im}(\sigma),\alpha\}$  
%is a basis of a positive three dimensional subspace $W\subset H^2(X,\RealNumbers)$.
%This basis determines an orientation of $W$, hence a generator of $H^2(\C_W,\Integers)$,
%hence a generator of $H^2(\widetilde{\C}_X,\Integers)$.

Let $\FM_\Lambda^0$ be a connected component of the moduli space of marked pairs
and $P_0:\FM_\Lambda^0\rightarrow \Omega_\Lambda$ the period map.
A marked pair $(X,\eta)$ in $\FM_\Lambda^0$ determines an orientation 
of $\widetilde{\C}_\Lambda$, via the isomorphism $\widetilde{C}_X\cong \widetilde{\C}_\Lambda$ induced by
the marking $\eta$. This orientation of $\widetilde{\C}_\Lambda$ is constant
throughout the connected component $\FM_\Lambda^0$.
In particular, for each class $h\in \Lambda$, with $(h,h)>0$, we get 
a choice of a connected component 
\begin{equation}
\label{eq-Omega-h-perp-plus}
\Omega_{h^\perp}^+
\end{equation}
of $\Omega_{h^\perp}$, compatible with the orientation of $\widetilde{\C}_\Lambda$ induced by 
$\FM_\Lambda^0$. 
%Denote the other component of $\Omega_{h^\perp}$ by $\Omega_{h^\perp}^-$.

Let $\Orient(\Lambda)$ be the set of two orientations of the positive cone $\widetilde{\C}_\Lambda$.
Let 
\begin{equation}
\label{eq-orient}
\orient \ : \ \FM_\Lambda \ \ \ \rightarrow \ \ \ \Orient(\Lambda)
\end{equation}
be the natural map constructed above.

%****************************************************************
% 
%****************************************************************
\section{A modular description of each fiber of the period map}
\label{sec-Kahler-cone}

We provide a modular description of the fiber of the period map 
$\FM_\Lambda^0\rightarrow \Omega_\Lambda$ from a connected component
$\FM_\Lambda^0$ of the moduli space of marked pairs
(Theorem \ref{thm-modular-description-of-the-fiber-of-the-period-map}).
Throughout this section $X$ is an irreducible holomorphic symplectic manifold,
which need not be projective.

%*********************************************************************
%
%*********************************************************************
\subsection{Exceptional divisors}

A reduced and irreducible effective divisor $D\subset X$
will be called a {\em prime} divisor.

\begin{defi}
\label{def-exceptional-divisor}
\begin{enumerate}
\item
A set $\{E_1, \dots, E_r\}$ of prime divisors is {\em exceptional}, if and only if
its Gram matrix $\left(([E_i],[E_j])\right)_{ij}$ is negative definite.
\item
An effective divisor $E$ is {\em exceptional}, if the support of $E$
is an exceptional set of prime divisors. 
\end{enumerate}
\end{defi}

\begin{defi}
\label{def-fundamental-exceptional-chamber}
The {\em fundamental exceptional chamber} 
$\FE_X$ is the cone of classes $\alpha$, such that 
$\alpha\in \C_X$, and $(\alpha,[E])>0$, for every prime exceptional divisor $E$.
\end{defi}

%****************************************************************
% 
%****************************************************************
\subsubsection{The fundamental exceptional chamber versus the birational K\"{a}hler cone}

%*********************
% Hide
%*********************
\hide{
We show next that 
for each uniruled divisor $D$, there exists an effective rational $1$-cycle $C$,
such that $C^\perp=D^\perp$ (Lemma \ref{lemma-intersection-with-fiber}).
Consequently, each rational chamber is contained in a unique uniruled chamber.
Let $X$ be an irreducible holomorphic symplectic manifold and 
$E\subset X$ a reduced, irreducible, and uniruled divisor. 
Denote by $f\in H_2(X,\Integers)$ the class of the generic rational curve in a family which 
dominates $E$.

\begin{lem}
\label{lemma-intersection-with-fiber}
There exists a positive constant $\lambda_E$, such that 
\[
([E],\alpha) \ \ \ = \ \ \ \lambda_E (f,\alpha), \ \ \ \forall \alpha\in H^{1,1}(X,\RealNumbers).
\]
\end{lem}

\begin{proof}
There exists\footnote{(???) True if $X$ is projective, in which case the Hilbert scheme
is projective, but when $X$ is K\"{a}hler, it is not clear
that the Douady space contains such a subvariety.
Nevertheless, the Douady space is K\"{a}hler 
\cite{fujiki-schumacher1,fujiki-schumacher2} and we can replace $U$ by an
open subset of the Douady space and replace
$\int_U(\bullet)$ below by $\int_U\omega^c(\bullet)$,
where $\omega$ is a K\"{a}hler class on $U$ and $c=\dim(U)+2-2n$ (???).}
a smooth and proper family $p:\C\rightarrow U$,
of rational curves over a smooth complex manifold $U$ of dimension $2n-2$,
and a morphism $ev:\C\rightarrow X$, such that 
$ev(\C)$ is a dense open subset of $E$, since $E$ is uniruled. 
Let $\sigma$ be a non-zero holomorphic $2$-form on $X$.
Then $ev^*(\sigma)=p^*(\tau)$, for some holomorphic $2$-form on $U$,
since $p$ is a $\PP^1$-bundle. Let $d$ be the degree of the
morphism $ev:\C\rightarrow E$.
There exists a positive integer $\lambda_X$, such that 
$(\alpha,\beta)=\lambda_X\int_X\alpha\beta(\sigma\bar{\sigma})^{n-1}$,
for all $\alpha$, $\beta\in H^{1,1}(X,\RealNumbers)$ \cite{beauville}.
We have
\[
(\alpha,[E])=\lambda_X\int _E \alpha(\sigma\bar{\sigma})^{n-1}=
\frac{\lambda_X}{d}\int_\C ev^*[\alpha(\sigma\bar{\sigma})^{n-1}]=
(\alpha,f)\frac{\lambda_X}{d}\int_U(\tau\bar{\tau})^{n-1}.
\]
The constant $\frac{\lambda_X}{d}\int_U(\tau\bar{\tau})^{n-1}$
is independent of $\alpha$ and it is positive, as both
$(\alpha,[E])$ and $(\alpha,f)$ are positive, when $\alpha$ is a K\"{a}hler class.
\end{proof}
%*********************
% End Hide
%*********************
}

Huybrechts and Boucksom stated an important result (Theorem \ref{thm-kahler-cone} below) 
in terms of another chamber, which we introduce next.

\begin{defi}
\label{def-chambers}
(\cite{boucksom-zariski-decomposition}, Section 4.2.2)
\begin{enumerate}
\item
A {\em rational effective $1$-cycle} $C$ is a linear combination, with positive 
integral coefficients, of irreducible rational curves on $X$. 
%The hyperplane $C^\perp\subset H^2(X,\RealNumbers)$
%consists of classes $\alpha$, such that $\int_C\alpha=0$.
%\item
%A {\em rational chamber} of $\C_X$ is a connected component of the complement 
%$\C_X\setminus \cup\{C^\perp \ : \ C \ \mbox{is a rational effective} \ 1\mbox{-cycle in} \ X\}$.
%\item
%The {\em fundamental rational chamber} is the subset of $\C_X$ consisting of classes 
%$\alpha\in \C_X$, such that $\int_C\alpha>0$, for every non-zero 
%rational effective $1$-cycle $C$.
\item
A {\em uniruled divisor} $D$ is an effective divisor each of which irreducible components 
$D_i$ is covered by rational curves. 
%of class $c\in H_2(X,\Integers)$, such that $([D_i],c)<0$ (???).
%\item
%A {\em uniruled chamber} of $\C_X$ is a connected component of the complement 
%$\C_X\setminus \cup\{D^\perp \ : \ D \ \mbox{is a uniruled divisor in} \ X\}$,
%where $D^\perp$ is the orthogonal complement to the class of $D$ with respect to
%the Beauville-Bogomolov pairing.
\item
The {\em fundamental uniruled chamber}  $\FU_X$
is the subset of $\C_X$ consisting of classes 
$\alpha\in \C_X$, such that $(\alpha,D)>0$, for every non-zero 
uniruled divisor $D$.
\item
The {\em birational K\"{a}hler cone} $\BK_X$ of $X$ is the 
union of $f^*\K_Y$, as $f$ ranges over all bimeromorphic maps $f:X\rightarrow Y$
to an irreducible holomorphic symplectic manifold $Y$. 
\end{enumerate}
\end{defi}

Note that the birational K\"{a}hler cone is not convex in general.

%Then 
%the closure of $\FU_X$ is equal to the closure of the birational K\"{a}hler cone
%$\BK_X$ (Remark \ref{rem-birational-Kahler-cone}). 
%The following theorem is due to Huybrechts and Boucksom. 

\begin{thm}
\label{thm-kahler-cone}
(\cite{huybrechts-kahler-cone} and \cite{boucksom-zariski-decomposition}, Theorem 4.3)
\begin{enumerate}
%\item
%The positive cone $\C_X$ is contained in the pseudo-effective cone. 
\item
\label{thm-item-kahler-cone-is-fundamental-rational-chamber}
The K\"{a}hler cone $\K_X$ is equal to the 
%fundamental rational chamber.
subset of $\C_X$ consisting of classes 
$\alpha\in \C_X$, such that $\int_C\alpha>0$, for every non-zero 
rational effective $1$-cycle $C$.
\item
\label{thm-item-fundamental-uniruled-chamber-is-closure-of-birl-kahler-cone}
Let $\alpha\in \C_X$ be a class, such that $\int_C\alpha\neq 0$, for every
rational $1$-cycle.
%which belongs to a rational chamber. 
Then $\alpha$ belongs to $\FU_X$, if and only if $\alpha$
belongs to the birational K\"{a}hler cone $\BK_X$.
%there exists an irreducible holomorphic symplectic manifold $Y$, 
%and a bimeromorphic map $f:X\rightarrow Y,$ such that $f_*(\alpha)$ is a 
%K\"{a}hler class of $Y$. 
\item
(\cite{boucksom-zariski-decomposition}, Theorem 4.3 part ii, and 
\cite{huybrechts-basic-results}, Corollary 5.2)
Let $\alpha\in \C_X$ be a class, which does not belong to $\FU_X$. 
Assume that $\int_C\alpha\neq 0$, for every
rational $1$-cycle.
%belongs to a rational chamber, which is not contained in the fundamental uniruled chamber. 
Then there exists 
an irreducible holomorphic symplectic manifold $Y$, and a bimeromorphic map 
$f:X\rightarrow Y$, such that $f_*(\alpha)=\beta+D'$, where $\beta$ is a K\"{a}hler 
class on $Y$ and $D'$ is a non-zero linear combination
%\footnote{In Boucksom's language, the sum $\beta+D'$ is the Zariski decomposition of
%$f_*(\alpha)$ and it is uniquely determined by the class $f_*(\alpha)$.} 
of finitely many 
uniruled reduced and irreducible divisors with positive real coefficients. 
\end{enumerate}
\end{thm}

\begin{rem}
\label{rem-birational-Kahler-cone}
Let $X$ be an irreducible holomorphic symplectic manifold. 
%The union $f^*\K_Y$, as $f$ ranges over all bimeromorphic maps $f:X\rightarrow Y$
%to an irreducible holomorphic symplectic manifold $Y$, is called the
%{\em birational K\"{a}hler cone} $\BK_X$ of $X$. 
Part (\ref{thm-item-fundamental-uniruled-chamber-is-closure-of-birl-kahler-cone})
of the theorem asserts that if a class $\alpha$ satisfies the assumptions
stated, then $\alpha$ is contained in $\FU_X$, if and only if it is contained in $\BK_X$.
The `only if' direction of part  
(\ref{thm-item-fundamental-uniruled-chamber-is-closure-of-birl-kahler-cone})
is stated in (\cite{boucksom-zariski-decomposition}, Theorem 4.3). 
The `if' part is the obvious direction. 
Indeed, let $f:X\rightarrow Y$ be a birational map,
such that $f_*(\alpha)$ is a K\"{a}hler class on $Y$. Let $D$ be an effective uniruled 
reduced and irreducible divisor in $X$, and $D'$ its strict transform in $Y$. We have
$([D],\alpha)=([D'],f_*(\alpha))>0$. Hence, $\alpha$ is in the fundamental uniruled 
chamber.
\end{rem}

%***********************
% Hide
%***********************
\hide{
\begin{caution}
We could not exclude the possibility that in part 
(\ref{thm-item-fundamental-uniruled-chamber-is-closure-of-birl-kahler-cone}) of 
Theorem \ref{thm-kahler-cone} the set $f^*\K_Y$ intersects more than one rational chamber 
in $\C_X$. In other words, we could not exclude the possibility
that the homomorphism $f_*:H_2(X,\Integers)\rightarrow H_2(Y,\Integers)$ 
maps the class of a rational effective $1$-cycle $C$ to a class $f_*[C]$,
such that neither $f_*[C]$, nor $-f_*[C]$ is effective. For that reason, we will not use
below the rational chambers other than $\K_X$.
\end{caution}
%***********************
% End hide
%***********************
}

Let $\overline{\BK}_X$ be the closure of the birational K\"{a}hler cone $\BK_X$ in $\C_X$.

\begin{prop}
\label{prop-closure-of-FE-and-FU-are-equal}
%\begin{enumerate}
%\item
%\label{cor-item-fu-same-as-Sigma}
The following inclusions and equality hold:
\[
\BK_X\subset \FU_X= \FE_X\subset \overline{\BK}_X.
\]
%\item
%\label{cor-item-difference-between-two-fundamental-chambers}
%The difference $\left[\FE_X\setminus \FU_X\right]$ is contained in the following union
%of hyperplanes.
%\[  
%\bigcup\left\{[D]^\perp  \ : \   D \ 
%\mbox{is a uniruled  divisor but not a prime exceptional divisor}
%\right\}.
%\]
%\end{enumerate}
\end{prop}

\begin{proof} 
%(Of Proposition \ref{prop-closure-of-FE-and-FU-are-equal})
%Part (\ref{cor-item-fu-same-as-Sigma}):
An exceptional divisor is uniruled, by
(\cite{boucksom-zariski-decomposition}, Proposition 4.7). 
%This follows also from Proposition \ref{prop-reflection-by-exceptional-divisor} below.
The inclusion $\FU_X\subset \FE_X$ follows.
We prove next the inclusion $\FE_X\subset\FU_X$. Let $\alpha$ be a class in $\FE_X$
and $D$ a prime uniruled divisor. If $[D]$ belongs to the closure 
$\overline{\C}_X$ of the positive cone, then $(\alpha,[D])>0$, since $\alpha$ belongs to $\C_X$.
Otherwise, $[D]$ is a prime exceptional divisor, and so $(\alpha,[D])>0$.
The inclusion $\FE_X\subset \FU_X$ follows.

The inclusion $\BK_X\subset  \FU_X$ follows from the `if' direction of
Theorem
\ref{thm-kahler-cone} part 
\ref{thm-item-fundamental-uniruled-chamber-is-closure-of-birl-kahler-cone},
and 
the inclusion  $\FE_X\subset \overline{\BK}_X$ 
follows from the `only if' direction.
%is due to (\cite{huybrechts-kahler-cone}, Proposition 4.2), which states that $\alpha\in\C_X$
%belongs to $\overline{\BK}_X$, if and only if $(\alpha,[D])\geq 0$, for every
%uniruled prime divisor $D$. Indeed, let $\alpha$ belong to $\FE_X$.
%If $([D],[D])\geq 0$, then
%$(\alpha,[D])\geq 0$, since the closure of the positive cone is self-dual with
%respect to the Beauville-Bogomolov pairing. 
%If $([D],[D])< 0$, then
%$(\alpha,[D])> 0$, by definition of $\FE_X$.
%******************
%  Hide
%******************
\hide{
The latter also implies that 
$\overline{\BK}_X$ is equal to the closure of $\FU_X$ in $\C_X$.
It suffices to prove the inclusion $\FE_X\subset \overline{\FU}_X$.
Let $\alpha$ be a class in $\FE_X$. Then $\alpha\in\C_X$, by definition. Let
$\alpha=P(\alpha)+N(\alpha)$ be the Zariski decomposition 
of Theorem \ref{thm-Zariski-decomposition}.
We need to show that $N(\alpha)=0$. Assume otherwise. Then
$(\alpha,N(\alpha))>0$, by definition of $\FE_X$. On the other hand, we have
\[
(\alpha,N(\alpha))=(P(\alpha)+N(\alpha),N(\alpha))=(N(\alpha),N(\alpha))<0.
\]
A contradiction. Hence, $N(\alpha)=0$ and $\alpha=P(\alpha)$ belongs to 
$\overline{\FU}_X$.
%******************
%  End Hide
%******************
}
\end{proof}
%Part (\ref{cor-item-difference-between-two-fundamental-chambers}):
%Let $\alpha$ be a class in $\left[\FE_X\setminus \FU_X\right]$. 
%Then there exists a uniruled divisor $D$, such that $(\alpha,[D])\leq 0$. 
%The inclusion $\FE_X\subset \overline{\FU}_X$, which
%follows from part (\ref{cor-item-fu-same-as-Sigma}), implies that
%$(\alpha,[D])\geq 0$. Hence, $(\alpha,[D])=0.$
%The divisor $D$ can not be exceptional, since $\alpha$ belongs to $\FE_X$.
%
%Part (\ref{cor-item-semi-rational-chambers-in-Sigma}):
%Let $\overline{\FE}_X$ be the closure of $\FE_X$ in $\C_X$. 
%The boundary $\partial\FE_X:=\overline{\FE}_X\setminus\FE_X$
%is contained in $[E]^\perp$, for some prime exceptional divisor $E$,
%by definition of $\FE_X$. In particular, $\partial\FE_X$ is either empty,
%or $\dim\left(\partial\FE_X\right)<h^{1,1}(X)$. Let $ch$ be a K\"{a}hler-type
%chamber, which intersects $\C_X\setminus \FE_X$. Then $ch$
%intersects $\C_X\setminus \overline{\FE}_X$, since $ch$ is an open subset of 
%$\C_X$ of dimension $h^{1,1}(X)$. Hence, $ch$ intersects 
%$\C_X\setminus \overline{\FU}_X$. Thus, $ch$ is contained in 
%$\C_X\setminus \overline{\FU}_X$, by Lemma \ref{lemma-on-semi-chambers}.
%
%Part (\ref{cor-item-seri-rational-chambers-in-Sigma-are-rational}):
%Let $ch$ be a K\"{a}hler-type chamber, which is contained in $\FE_X$.
%Then $ch$ intersects $\FU_X$, by Part (\ref{cor-item-fu-same-as-Sigma}).
%Hence, $ch$ is contained in $\BK_X$, by Lemma \ref{lemma-on-semi-chambers}.
%\end{proof}

The notation $\FE_X$  will replace $\FU_X$ from now on, in view
of Proposition \ref{prop-closure-of-FE-and-FU-are-equal}.
A class $\alpha\in \C_X$ is said to be {\em very general}, if 
$\alpha^\perp\cap H^{1,1}(X,\Integers)=0$.

\begin{cor}
\label{cor-stabilizer-of-birational-Kahler-cone}
Let $X_1$ and $X_2$ be irreducible holomorphic symplectic manifolds, 
$g:H^2(X_1,\Integers)\rightarrow H^2(X_2,\Integers)$ a parallel transport operator,
which is an isomorphism of Hodge structures, and 
$\alpha_1\in \FE_{X_1}$ a very general class. 
Then $g(\alpha_1)$ belongs to $\FE_{X_2}$, if and only if 
there exists a bimeromorphic map $f:X_1\rightarrow X_2$, such that $g=f_*$. 
\end{cor}

\begin{proof}
The `if' part is clear, since $f_*$ induces a bijection between the sets of exceptional 
divisors on $X_i$, $i=1,2$.
Set $\alpha_2:=g(\alpha_1)$. There exist irreducible holomorphic symplectic manifolds $Y_i$
and bimeromorphic maps $f_i:X_i\rightarrow Y_i$, 
such that $f_{i_*}(\alpha_i)$ is a K\"{a}hler class on $Y_i$,
by part (\ref{thm-item-fundamental-uniruled-chamber-is-closure-of-birl-kahler-cone}) 
of Theorem \ref{thm-kahler-cone}.
The homomorphisms $f_{i_*}:H^2(X_i,\Integers)\rightarrow H^2(Y_i,\Integers)$
are parallel transport operators, by Theorem
\ref{thm-bimeromorphic-implies-inseparable}.
Thus $(f_2^{-1})^*\circ g\circ f_1^*:H^2(Y_1,\Integers)\rightarrow H^2(Y_2,\Integers)$
is a parallel transport operator and a Hodge-isometry,
mapping the K\"{a}hler class $f_{1_*}(\alpha_1)$ to the 
the K\"{a}hler class $f_{2_*}(\alpha_2)$. Hence, 
there exists an isomorphism $h:Y_1\rightarrow Y_2$, such that 
$h_*=(f_2^{-1})^*\circ g\circ f_1^*,$ 
by Theorem \ref{thm-Hodge-theoretic-Torelli}.
Thus, $g=[(f_2)^{-1}hf_1]_*$.
\end{proof}

%*********************************************************************
%
%*********************************************************************
\subsubsection{The divisorial Zariski decomposition}
The following fundamental result of Bouksom will be needed in section 
\ref{sec-stably-prime-exceptional-line-bundles}.
The {\em effective cone} of $X$ is the cone 
in $H^{1,1}(X,\Integers)\otimes_\Integers\RealNumbers$ generated
by the classes of effective divisors. The
{\em algebraic pseudo-effective cone} $\Peff_X$ 
is the closure of  the effective cone.
Boucksom defines a larger transcendental analogue, a cone in 
$H^{1,1}(X,\RealNumbers)$, which he calls the 
{\em pseudo-effective cone} (\cite{boucksom-zariski-decomposition}, section 2.3).
We will not need the precise definition, but only the fact that the
pseudo-effective cone contains $\C_X$ 
(\cite{boucksom-zariski-decomposition}, Theorem 4.3  part (i)).
The sum $\C_X+\Peff_X$ is thus a sub-cone 
of Boucksom's pseudo-effective cone in $H^{1,1}(X,\RealNumbers)$.
Denote by $\overline{\FE}_X$ the closure of the fundamental exceptional chamber
in $H^{1,1}(X,\RealNumbers)$. 
%Denote by $\Peff_X\subset H^{1,1}(X,\Integers)\otimes_\Integers\RealNumbers$ the 
%closure of the cone generated by effective line bundles. 
%The following is a fundamental result of Boucksom.

\begin{thm}
\label{thm-Zariski-decomposition}
\begin{enumerate}
\item
(\cite{boucksom-zariski-decomposition}, Theorem 4.3 part (i), 
Proposition 4.4, and Theorem 4.8).
Let $X$ be an irreducible holomorphic symplectic manifold and
$\alpha$ a class in $\C_X+\Peff_X$. 
Then there exists a unique decomposition
\[
\alpha \ \ = \ \ P(\alpha) + N(\alpha),
\]
where $(P(\alpha),N(\alpha))=0$, $P(\alpha)$ belongs to 
$\overline{\FE}_X$, and $N(\alpha)$
is an exceptional $\RealNumbers$-divisor.
\item
\label{cor-item-fixed-part}
(\cite{boucksom-zariski-decomposition}, Corollary 4.11).
Let $L$ be a line bundle with $c_1(L)\in \C_X+\Peff_X$. 
Set $\alpha:=c_1(L)$.
%belongs to $\C_X\cap H^{1,1}(X,\RationalNumbers)$,
Then the classes $P(\alpha)$ and $N(\alpha)$ correspond to 
$\RationalNumbers$-divisors classes, which we denote by $P(\alpha)$ and $N(\alpha)$
as well. Furthermore, the homomorphism
\[
H^0\left(X,\StructureSheaf{X}(kP(\alpha))\right) \rightarrow
H^0(X,L^k)
\]
is surjective, for every non-negative integer $k$, such that $kP(\alpha)$
is an integral class.
\end{enumerate}
\end{thm}

\begin{rem}
%Theorem 4.8, and  Corollary 4.11 in \cite{boucksom-zariski-decomposition}
%are stated for a class $\alpha$ in the pseudo-effective cone (in the
%sense of \cite{boucksom-zariski-decomposition}, section 2.3).
%The cone $\C_X+\Peff_X$ is contained in the pseudo-effective cone, by  
%(\cite{boucksom-zariski-decomposition}, Theorem 4.3 part (i)).
The class $P(\alpha)$ is stated as a class in the {\em modified nef cone} in 
(\cite{boucksom-zariski-decomposition}, Theorem 4.8), but the modified nef cone
is equal to the closure of the birational K\"{a}hler cone, 
by (\cite{boucksom-zariski-decomposition}, 
Proposition 4.4), and hence also to $\overline{\FE}_X$.
\end{rem}

Part (\ref{cor-item-fixed-part}) of the above Theorem implies that
%if $L$ is a line bundle of non-negative Beauville-Bogomolov degree, then
the exceptional divisor $N(kc_1(L))$ is the fixed part of the linear system
$\linsys{L^k}$. In particular, if $c_1(L)=N(c_1(L))$, then the linear system
$\linsys{L^k}$ is either empty, or consists of a single exceptional divisor.
Exceptional divisors are thus rigid.

%****************************************************************
% 
%****************************************************************
\subsection{A K\"{a}hler-type chamber decomposition of the positive cone}
Let $X$ be an irreducible holomorphic symplectic manifold.
Denote the subgroup of $Mon^2(X)$ preserving the weight $2$ Hodge structure
by $Mon^2_{Hdg}(X)$.  
Note that the positive cone $\C_X$ is invariant under $Mon^2_{Hdg}(X)$,
since the orientation class of $\widetilde{\C}_X$ is invariant under the whole
monodromy group $Mon^2(X)$ (see section \ref{sec-orientation}).

\begin{defi}
\label{def-Kahler-type-chambers}
\begin{enumerate}
\item
An {\em exceptional chamber}
%\footnote{The prefix `semi' is added to distinguish these from
%uniruled chambers defined in \cite{boucksom-zariski-decomposition}. Uniruled chambers,
%other than the fundamental uniruled chamber, will not
%be mentioned elsewhere in this paper. The term 
%semi-uniruled is superfluous when $X$ is projective, since we 
%will introduce later the exceptional chambers 
%(Definition \ref{def-exceptional-chamber}) and prove that a
%semi-uniruled chamber is the same as an exceptional chamber
%(Theorem \ref{thm-numerical-characterization-of-BK-X}).} 
of the positive cone $\C_X$ is a subset of the form
$g[\FE_X]$, $g\in Mon^2_{Hdg}(X)$.
\item
A {\em K\"{a}hler-type chamber} of the positive cone $\C_X$  is a subset of the form
$g[f^*(\K_Y)]$, where $g\in Mon^2_{Hdg}(X)$, and $f:X\rightarrow Y$ is a bimeromorphic 
map to an irreducible holomorphic symplectic manifold $Y$.
\end{enumerate}
\end{defi}

Let  $Mon^2_{Bir}(X)\subset Mon^2_{Hdg}(X)$ be the subgroup of monodromy
operators induced by bimeromorphic maps from $X$ to itself
(see Theorem \ref{thm-bimeromorphic-implies-inseparable}).

\begin{lem}
\label{lemma-on-semi-chambers}
\begin{enumerate}
\item
\label{lemma-item-very-general-class-is-in-semi-uniruled-chamber}
Every very general class $\alpha\in \C_X$ belongs to some K\"{a}hler-type chamber.
\item
\label{lemma-item-semi-rational-chamber-contained-in-semi-uniruled-one}
Every K\"{a}hler-type chamber is contained in some exceptional chamber.
\item
\label{lemma-item-intersecting-semi-rational-chambers-are-equal}
If two K\"{a}hler-type chambers intersect, then they are equal.
\item
\label{lemma-item-overlapping-semi-uniruled-chambers}
If two exceptional chambers $g_1[\FE_X]$ and $g_2[\FE_X]$ contain
a common very general class $\alpha$, then they are equal.
\item
\label{lemma-item-transitive-action}
$Mon^2_{Hdg}(X)$ acts transitively on the set of exceptional chambers.
\item
\label{lemma-item-stabilizer-of-fundamental-uniruled-chamber}
The subgroup of $Mon^2_{Hdg}(X)$ stabilizing $\FE_X$ 
is equal to $Mon^2_{Bir}(X)$.
\end{enumerate}
\end{lem}

\begin{proof}
Part (\ref{lemma-item-very-general-class-is-in-semi-uniruled-chamber}):
There exists an irreducible holomorphic symplectic manifold $\widetilde{X}$
and a correspondence $\Gamma:=Z+\sum_i Y_i$ in $X\times \widetilde{X}$,
such that $Z$ is the graph of a bimeromorphic map $f:X\rightarrow \widetilde{X}$,
the restriction $g:H^2(X,\Integers)\rightarrow H^2(\widetilde{X},\Integers)$ of 
$[\Gamma]_*$ is a parallel transport operator, and $g(\alpha)$ is a 
K\"{a}hler class of $\widetilde{X}$, by 
(\cite{huybrechts-basic-results}, Corollary 5.2). Set $h:=f^*\circ g$.
Then $h$ belongs to $Mon^2_{Hdg}(X)$, by
Theorem \ref{thm-bimeromorphic-implies-inseparable}, 
$h(\alpha)=(f^*\circ g)(\alpha)$ belongs to $f^*\K_{\widetilde{X}}$, and $f^*\K_{\widetilde{X}}$
is a K\"{a}hler-type chamber, by Definition \ref{def-chambers}. 
Consequently, $h^{-1}(f^*\K_{\widetilde{X}})$ is a K\"{a}hler-type chamber 
containing $\alpha$.

Part (\ref{lemma-item-semi-rational-chamber-contained-in-semi-uniruled-one}):
Let $Ch$ be the K\"{a}hler-type chamber $g[f^*(\K_Y)]$, where $f$, $g$, and $Y$
are as in Definition \ref{def-Kahler-type-chambers}. Then $f^*(\FE_Y)=\FE_X$, 
by Corollary \ref{cor-stabilizer-of-birational-Kahler-cone}, and so
$Ch$ is contained in the exceptional chamber $g[\FE_X]$.

Part (\ref{lemma-item-intersecting-semi-rational-chambers-are-equal}):
Let $Y_i$ be irreducible holomorphic symplectic manifolds, $f_i:X\rightarrow Y_i$
bimeromorphic maps, $g_i\in Mon^2_{Hdg}(X)$, $i=1,2$, and 
$\alpha$ a class in $g_1[f_1^*(\K_{Y_1})]\cap g_2[f_2^*(\K_{Y_2})]$. 
%The K\"{a}hler cones $\K_{Y_i}$ are open, so we may assume that $\alpha$ is very general.
%Set $g:=g_1^{-1}g_2$. Then 
%$g_1^{-1}(\alpha)$ belongs to $[f_1^*(\K_{Y_1})]\cap g[f_2^*(\K_{Y_2})]$. 
%Thus, $g[\FU_X]=\FU_X$, by part (\ref{lemma-item-overlapping-semi-uniruled-chambers}).
%Hence, there exists a bimeromorphic map $h:X\rightarrow X$, such that 
%$g=h^*$, by part (\ref{lemma-item-stabilizer-of-fundamental-uniruled-chamber}).
The composition 
$\phi:=f_{2_*}\circ g_2^{-1}\circ g_1\circ f_1^*:H^2(Y_1,\Integers)
\rightarrow H^2(Y_2,\Integers)$ is a parallel-transport operator,
which maps the K\"{a}hler class $f_{1_*}(g_1^{-1}(\alpha))$ to the 
K\"{a}hler class $f_{2_*}(g_2^{-1}(\alpha))$.
Hence, $\phi$ is induced by an isomorphism $\tilde{\phi}:Y_1\rightarrow Y_2$, 
by Theorem \ref{thm-Hodge-theoretic-Torelli}. We get the equality 
$g_1^{-1}g_2f_2^*(\K_{Y_2})=f_1^*\tilde{\phi}^*(\K_{Y_2})=f_1^*(\K_{Y_1})$.
Consequently, $g_1[f_1^*(\K_{Y_1})]= g_2[f_2^*(\K_{Y_2})]$.

Part (\ref{lemma-item-overlapping-semi-uniruled-chambers}):
Set $g:=g_2^{-1}g_1$ and $\beta:=g_2^{-1}(\alpha)$. Then $\beta$ belongs to the intersection 
$g[\FE_X]\cap\FE_X$. So $g^{-1}(\beta)$ and
$\beta$ both belong to $\FE_X$ and $g$ maps the former to the latter. Hence, 
$g$ is induced by a birational map from $X$ to itself, 
by Corollary \ref{cor-stabilizer-of-birational-Kahler-cone}.
Thus, $g[\FE_X]=\FE_X$ and so $g_1[\FE_X]=g_2[\FE_X]$.

Part (\ref{lemma-item-transitive-action}):
%The action is well defined, by part
%(\ref{lemma-item-overlapping-semi-uniruled-chambers}). 
The action is transitive, by definition. 

%Completion of the proof of part (\ref{lemma-item-overlapping-semi-uniruled-chambers}):
%The set of very general classes is invariant under $Mon^2_{Hdg}(X)$.
%We may thus assume that $g_2=id$, by part (\ref{lemma-item-transitive-action}).

Part (\ref{lemma-item-stabilizer-of-fundamental-uniruled-chamber})
is an immediate consequence of Corollary \ref{cor-stabilizer-of-birational-Kahler-cone}.
\end{proof}

\begin{lem}
\label{lemma-hodge-parallel-transport-map-semi-chambers-to-such}
Let $X_1$ and $X_2$ be irreducible holomorphic symplectic manifolds and
$g:H^2(X_1,\Integers)\rightarrow H^2(X_2,\Integers)$ a parallel transport operator,
which is an isomorphism of Hodge structures. 
\begin{enumerate}
\item
\label{lemma-item-g-maps-semi-uniruled-chambers-to-such}
$g$ maps each exceptional chamber  in $\C_{X_1}$
onto an exceptional chamber  in $\C_{X_2}$.
\item
\label{lemma-item-g-maps-Kahler-type-chambers-to-such}
$g$ maps each K\"{a}hler-type chamber in $\C_{X_1}$ onto a K\"{a}hler-type chamber in
$\C_{X_2}$. 
\end{enumerate}
\end{lem}

\begin{proof}
There exists a bimeromorphic map $h:X_1\rightarrow X_2$, 
by Theorem \ref{thm-Hodge-theoretic-Torelli}. 
The homomorphism $h_*:H^2(X_1,\Integers)\rightarrow H^2(X_2,\Integers)$
is a parallel transport operator, and an isomorphism of Hodge structures,
by Theorem \ref{thm-bimeromorphic-implies-inseparable}.
%Hence, $f:=h^*\circ g$ belongs to $Mon^2_{Hdg}(X_1)$.

Part (\ref{lemma-item-g-maps-semi-uniruled-chambers-to-such}):
Let $f$ be an element of $Mon^2_{Hdg}(X_1)$. 
We need to show that $g(f[\FE_{X_1}])$ is an exceptional chamber in $\C_{X_2}$. 
Indeed, we have the equalities
\[
g(f[\FE_{X_1}])=
%(fh_*\tilde{f})[\FE_{X_1}]=
(gfh^*)\left\{h_*[\FE_{X_1}]\right\}=
(gfh^*)[\FE_{X_2}],
\]
and $gfh^*$ belongs to $Mon^2_{Hdg}(X_2)$.

Part (\ref{lemma-item-g-maps-Kahler-type-chambers-to-such}):
Any K\"{a}hler-type chamber of $\C_{X_1}$ is of the form 
$f[\tilde{h}^*(\K_{Y_1})]$, where $\tilde{h}:X_1\rightarrow Y_1$ is a 
bimeromorphic map to an irreducible holomorphic symplectic manifold $Y_1$,
and $f$ is an element of $Mon^2_{Hdg}(X_1)$. 
We have the equality
\[
gf[\tilde{h}^*(\K_{Y_1})] \ \ = \ \ 
(gfh^*)\left\{(h\tilde{h}^{-1})_*(\K_{Y_1})\right\},
\]
$(h\tilde{h}^{-1})_*(\K_{Y_1})$ is a K\"{a}hler-type chamber of $X_2$ and 
$gfh^*$ belongs to $Mon^2_{Hdg}(X_2)$, 
by Theorem \ref{thm-bimeromorphic-implies-inseparable}.
Thus $gf[\tilde{h}^*(\K_{Y_1})]$ is a K\"{a}hler-type chamber of $X_2$.
\end{proof}

\begin{cor}
\label{cor-inseparable-pairs-and-semi-rational-chambers}
Let $(X_1,\eta_1)$, $(X_2,\eta_2)$ be two inseparable points in ${\mathfrak M}^0_\Lambda$. 
\begin{enumerate}
\item
\label{cor-item-inseparable-pairs-and-semi-uniruled-chambers}
The composition $\eta_2^{-1}\circ \eta_1$ 
maps each K\"{a}hler-type chamber in $\C_{X_1}$ onto a K\"{a}hler-type chamber in
$\C_{X_2}$. 
Similarly, $\eta_2^{-1}\circ \eta_1$ maps each exceptional chamber  in $\C_{X_1}$
onto an exceptional chamber  in $\C_{X_2}$.
\item
\label{cor-item-inseparable-pairs-and-semi-rational-chambers}
$(\eta_2^{-1}\circ \eta_1)(\FE_{X_1})=\FE_{X_2}$, if and only if 
there exists a bimeromorphic map 
$f$ from $X_1$ to $X_2$, such that $\eta_2^{-1}\circ \eta_1=f_*$.
\end{enumerate}
\end{cor}

\begin{proof}
The composition $\eta_2^{-1}\circ \eta_1$ is a parallel-transport operator, 
and a Hodge-isometry, by Theorem 
\ref{thm-non-separated-implies-bimeromorphic} part 
\ref{thm-item-composition-of-markings-is-a-parallel-transport-op}.
Part (\ref{cor-item-inseparable-pairs-and-semi-uniruled-chambers})
follows from Lemma 
\ref{lemma-hodge-parallel-transport-map-semi-chambers-to-such}.
Part (\ref{cor-item-inseparable-pairs-and-semi-rational-chambers})
follows from Corollary \ref{cor-stabilizer-of-birational-Kahler-cone}.
\end{proof}

%*********************************************************************
%
%*********************************************************************
\subsection{$\FM_\Lambda$ as the moduli space of K\"{a}hler-type chambers}
Consider the period map $P_0:\FM_\Lambda^0\rightarrow \Omega_\Lambda$
from the connected component $\FM_\Lambda^0$ containing the isomorphism 
class of the marked pair $(X,\eta)$.
Denote by $\KT(X)$ the set of K\"{a}hler-type chambers in $\C_X$.
Let 
\begin{equation}
\label{eq-from-fiber-to-SR}
\rho \ : \ P_0^{-1}\left[P_0(X,\eta)\right] \ \ \ \longrightarrow \ \ \ \KT(X)
\end{equation}
be the map given by $\rho(\widetilde{X},\tilde{\eta})=(\eta^{-1}\tilde{\eta})(\K_{\widetilde{X}}).$
%The homomorphism $\eta^{-1}\tilde{\eta}$ is a parallel transport operator, 
%and an isomorphism of Hodge structures, by
%Theorem \ref{thm-non-separated-implies-bimeromorphic} part
%\ref{thm-item-composition-of-markings-is-a-parallel-transport-op}.
The map $\rho$ is well defined, by 
Corollary \ref{cor-inseparable-pairs-and-semi-rational-chambers}.
$Mon^2_{Hdg}(X)$ acts on $\KT(X)$, by Lemma 
\ref{lemma-hodge-parallel-transport-map-semi-chambers-to-such}.

Note that each period $P(X,\eta)\in\Omega_\Lambda$
is invariant under the subgroup
\begin{equation}
\label{eq-Mon-2-Hdg-conjugated-by-eta}
Mon^2_{Hdg}(X)^\eta \ \ := \ \ \{\eta g \eta^{-1} \ : \ g \in Mon^2_{Hdg}(X)\}
\end{equation}
of $O(\Lambda)$.
Consequently, $Mon^2_{Hdg}(X)$ acts on the fiber $P_0^{-1}[P_0(X,\eta)]$
of the period map by
\[
g(\widetilde{X},\tilde{\eta}) \ \ := \ \ (\widetilde{X},\eta g \eta^{-1}\tilde{\eta}).
\]

\begin{prop}
\label{prop-fiber}
\begin{enumerate}
\item
\label{prop-item-rho-is-bijective}
The map $\rho$ is a $Mon^2_{Hdg}(X)$-equivariant bijection. 
\item
\label{prop-item-single-Hausdorff-point}
The marked pair $(X,\eta)$ is a Hausdorff point of $\FM_\Lambda$,
if and only if $\C_X=\K_X$.
\item
\label{prop-item-cyclic-picard-and-projective-imply-K-X=C-X}
(\cite{huybrechts-basic-results}, Corollaries 5.7 and 7.2)
$\C_X=\K_X$, if $H^{1,1}(X,\Integers)$ is trivial, or if $H^{1,1}(X,\Integers)$ 
is of rank $1$, generated by a class $\alpha$ of
non-negative  Beauville-Bogomolov degree.
\end{enumerate}
\end{prop}

\begin{proof}
Part (\ref{prop-item-rho-is-bijective}):
Assume that $\rho(X_1,\eta_1)=\rho(X_2,\eta_2)$.
Then $\eta_2^{-1}\eta_1(\K_{X_1})=\K_{X_2}.$ Hence, 
$\eta_2^{-1}\eta_1=f_*$, for an isomorphism $f:X_1\rightarrow X_2$, 
by Theorem \ref{thm-Hodge-theoretic-Torelli}.
Thus, $(X_1,\eta_1)$ and $(X_2,\eta_2)$ are isomorphic, and $\rho$ is injective.

Given a K\"{a}hler-type chamber $Ch$ in $\C_X$ 
and a very general class $\alpha$ in $Ch$, 
there exists an element $g\in Mon^2_{Hdg}(X)$, such that 
$g(\alpha)$ belongs to $\FE_X$, by Lemma 
\ref{lemma-on-semi-chambers} part \ref{lemma-item-transitive-action}.
There exists an irreducible holomorphic symplectic manifold $Y$
and a bimeromorphic map $h:X\rightarrow Y$, such that $h_*(g(\alpha))$
belongs to $\K_Y$, by Theorem \ref{thm-kahler-cone} part 
\ref{thm-item-fundamental-uniruled-chamber-is-closure-of-birl-kahler-cone}.
Thus, $(h_*\circ g)(Ch)=\K_Y$, by Lemma 
\ref{lemma-hodge-parallel-transport-map-semi-chambers-to-such}.
We conclude that $\rho(Y,\eta\circ g^{-1}\circ h^*)=g^{-1}h^*(\K_Y)=Ch$
and $\rho$ is surjective.

Part (\ref{prop-item-single-Hausdorff-point}) follows from part 
(\ref{prop-item-rho-is-bijective}).
%Part (\ref{prop-item-cyclic-picard-and-projective-imply-K-X=C-X}) is proven in
\end{proof}

Fix a connected component $\FM_\Lambda^0$ of the moduli space of marked pairs.
We get the following modular description of the fiber $P_0^{-1}(\period)$
in terms of the period $\period$. Set $\Lambda^{1,1}(\period,\RealNumbers):=
\{\lambda\in \Lambda\otimes_\Integers\RealNumbers \ : \ (\lambda,\period)=0\}$.
Let $\C_\period$ be the connected component, of the  cone $\C'_\period$ in
$\Lambda^{1,1}(\period,\RealNumbers)$, which is compatible with the 
orientation of the positive cone $\widetilde{\C}_\Lambda$ 
determined by $\FM_\Lambda^0$ (see section \ref{sec-orientation}).

\begin{defi}
A {\em K\"{a}hler-type chamber} of $\C_\period$ is a subset of the form 
$\eta(Ch)\subset\C_\period$, where $(X,\eta)$ is 
a marked pair $\FM_\Lambda^0$ and $Ch\subset \C_X$ is a K\"{a}hler-type chamber of $X$.
\end{defi}

Denote by $\KT(\period)$ the set of K\"{a}hler-type chambers in $\C_\period$. 
The map 
\[
\eta \ : \ \KT(X) \ \ \ \longrightarrow \ \ \ \KT(\period),
\]
sending a K\"{a}hler-type chamber $Ch\in \KT(X)$ to $\eta(Ch)$,
is a bijection, for every marked pair
$(X,\eta)$ in the fiber $P_0^{-1}(\period)$, by 
Corollary \ref{cor-inseparable-pairs-and-semi-rational-chambers} and
Proposition \ref{prop-fiber}. 
$Mon^2_{Hdg}(X)^\eta$, given in equation (\ref{eq-Mon-2-Hdg-conjugated-by-eta}),
is the same  subgroup of $O(\Lambda)$,
for all $(X,\eta)\in P_0^{-1}(\period)$, and we denote it by $Mon^2_{Hdg}(\period)$.
The following statement is an immediate corollary of Proposition \ref{prop-fiber}.

\begin{thm}
\label{thm-modular-description-of-the-fiber-of-the-period-map}
The map
\[
\rho \ : \ P_0^{-1}(\period) \ \ \ \longrightarrow \ \ \ \KT(\period),
\]
given by $\rho(X,\eta):=\eta(\K_X)$, 
is a $Mon^2_{Hdg}(\period)$-equivariant bijection.
\end{thm}

%Proposition \ref{prop-fiber} provides a bijection between the fiber
%$P_0^{-1}(\period)$ and K\"{a}hler-type chambers in $\C_\period$.

\begin{rem}
Compare Theorem \ref{thm-modular-description-of-the-fiber-of-the-period-map} 
with the more detailed analogue for $K3$ surfaces, which is provided
in (\cite{looijenga-peters}, Theorem 10.5).
Ideally, one would like to have a description of the set $\KT(\period)$,
depending only on the period $\period$, the deformation type of $X$,
and possibly some additional discrete monodromy invariant of $X$ (see the 
invariant $\iota_X$ introduced in Corollary
\ref{cor-third-characterization-of-Mon-2}). 
Such a description would depend on the
determination of the K\"{a}hler-type chambers in $\C_X$.
In particular, it requires a determination of the
K\"{a}hler cone of an irreducible holomorphic symplectic variety, 
in terms of the Hodge structure of $H^2(X,\Integers)$, the
Beauville-Bogomolov pairing, and the discrete monodromy invariants of $X$.
The determination of the K\"{a}hler cone $\K_X$ in terms of such data  
is a very difficult problem addressed in a sequence of papers 
of Hassett and Tschinkel 
\cite{hassett-tschinkel-conj,hassett-tschinkel,hassett-tschinkel-monodromy,hassett-tschinkel-intersection-numbers}.
Precise conjectures for the determination of the K\"{a}hler cones in the $K3^{[n]}$-type,
for all $n$, 
and for generalized Kummer fourfolds, are provided in 
\cite{hassett-tschinkel-intersection-numbers}, Conjectures 1.2 and 1.4.
The determination of the birational K\"{a}hler cone,
in terms of such data,  is the subject of section 
\ref{sec-generators-for-Mon-2-Hdg}. 
\end{rem}

%*********************************************************************
%
%*********************************************************************
\section{$Mon^2_{Hdg}(X)$ is generated by reflections and $Mon^2_{Bir}(X)$}
\label{sec-generators-for-Mon-2-Hdg}
Throughout this section $X$ denotes a {\em projective} 
irreducible holomorphic symplectic manifold.
Under the projectivity assumption, we can define a subgroup
$W_{Exc}$ of the Hodge-monodromy group $Mon^2_{Hdg}(X)$, which is 
generated by reflections with respect to classes of prime exceptional divisors
(Definition \ref{def-W-Exc} and Theorem \ref{thm-Mon-2-Hdg-is-a-semi-direct-product}
part \ref{thm-item-W-is-generated-by-exceptional-reflections}).
The fundemental exceptional chamber $\FE_X$,
introduced in Definition \ref{def-fundamental-exceptional-chamber}, is the
interior of a fundamental domain for the action of the reflection group $W_{Exc}$
on the positive cone $\C_X$. Significant regularity properties follow from this description of 
$\FE_X$ (Theorem \ref{thm-numerical-characterization-of-BK-X}).
We prove also that $W_{Exc}$ is a normal subgroup of $Mon^2_{Hdg}(X)$ and
the latter is a semi-direct-product of $W_{Exc}$ and $Mon^2_{Bir}(X)$ 
(Theorem \ref{thm-Mon-2-Hdg-is-a-semi-direct-product}).
A weak version of Morrison's movable cone conjecture follows from the
above results in the special case of irreducible holomorphic symplectic manifolds  
(Theorems \ref{thm-finite-number-of-Bir-X-orbits} and 
\ref{thm-movable-conj-conj}).

\subsection{Reflections}
\label{sec-reflections}
Let $X$ be a projective irreducible holomorphic symplectic manifold of dimension $2n$
and $E\subset X$ a prime exceptional divisor (Definition \ref{def-exceptional-divisor}).

\begin{prop}
\label{prop-druel}
(\cite{druel}, Proposition 1.4)
There exists a sequence of flops of $X$,
resulting in a smooth birational model $X'$ of $X$,
such that the strict transform $E'$ of $E$ in $X'$ is contractible
via a projective birational 
morphism $\pi:X' \rightarrow Y$ onto a normal projective variety 
$Y$. The exceptional locus of $\pi$ is equal to the support of $E'$.
\end{prop}

Identify $H^2(X,\RationalNumbers)^*$ with $H_2(X,\RationalNumbers)$. 
Set
\[
%\begin{equation}
%\label{eq-E-vee}
[E]^\vee \ \ \ := \ \ \ \frac{-2([E],\bullet)}{([E],[E])} \ \ \in \ \ H_2(X,\RationalNumbers).
%\end{equation}
\]
\begin{prop}
\label{prop-reflection-by-exceptional-divisor}
(\cite{markman-prime-exceptional}, Corollary 3.6 part 1).
\begin{enumerate}
\item
\label{prop-item-reflection-is-integral}
There exists a Zariski dense open subset $E^0\subset E$ 
and a proper holomorphic fibration $\pi:E^0\rightarrow B$, 
onto a smooth holomorphic symplectic variety of dimension
$2n-2$, with the following property. 
The class $[E]^\vee$ is the class of a generic fiber of $\pi$. 
The generic fiber is either a smooth rational curve, or the union of two
homologous smooth  rational curves meeting at one point non-tangentially.
In particular, the class $[E]^\vee$ is integral, as is the reflection
$R_E:H^2(X,\Integers)\rightarrow H^2(X,\Integers)$, given by
$R_E(x)=x+(x,[E]^\vee)[E]$.
\item
The reflection $R_E$ belongs to $Mon^2_{Hdg}(X)$.
\end{enumerate}
\end{prop}

\begin{rem}
\begin{enumerate}
\item
The proof of Proposition
\ref{prop-reflection-by-exceptional-divisor}
relies heavily on Druel's result stated above in Proposition
\ref{prop-druel}.
The fact that $R_{E'}$ belongs to $Mon^2_{Hdg}(X')$
was proven earlier in (\cite{markman-galois}, Theorem 1.4), using fundamental work of 
Namikawa \cite{namikawa-deformations} (see \cite{namikawa-galois} for an alternative proof).
The author does not know if the 
analogue of Proposition \ref{prop-druel} holds for a non-projective 
irreducible holomorphic symplectic manifold $X$ as well. This is the reason for
the projectivity assumption throughout  section 
\ref{sec-generators-for-Mon-2-Hdg}.
\item
The variety $B$ in part (\ref{prop-item-reflection-is-integral}) 
of the proposition is an \'{e}tale cover of a Zariski open subset of the image of $E'$ in $Y$
(\cite{namikawa-deformations}, section 1.8).
\end{enumerate}
\end{rem}
%*********************************************************************
%
%*********************************************************************
\subsection{Stably prime-exceptional line bundles}
\label{sec-stably-prime-exceptional-line-bundles}
Let $X$ be an irreducible holomorphic symplectic manifold. 
Denote by $Def(X)$ the local Kuranishi deformation space of $X$
and let $0\in Def(X)$ be the special point corresponding to $X$. 
Let $L$ be a line bundle on $X$. 
Set $\Lambda:=H^2(X,\Integers)$.
The period map $P:Def(X)\rightarrow \Omega_\Lambda$ embeds 
$Def(X)$ as an open analytic subset of the period domain $\Omega_\Lambda$
and the intersection $Def(X,L):=Def(X)\cap c_1(L)^\perp$ is the Kuranishi deformation space
of the pair $(X,L)$, i.e., it consists of deformations of the complex structure of $X$
along which $c_1(L)$ remains of type $(1,1)$. 
We assume that both $Def(X)$ and the intersection $Def(X,L)$
are simply connected, possibly after replacing $Def(X)$ by a smaller open neighborhood
of $0$ in the Kuranishi deformation space, which we denote again by $Def(X)$. 

Let $\pi:\X\rightarrow Def(X)$ be the universal family and denote by $X_t$
the fiber of $\pi$ over  $t\in Def(X)$. Denote by $\ell$ the flat section
of the local system $R^2{\pi_*}\Integers$ through $c_1(L)$ and let
$\ell_t\in H^{1,1}(X_t,\Integers)$ be its value at $t\in Def(X,L)$. 
Let $L_t$ be the line bundle on $X_t$ with $c_1(L_t)=\ell_t$. 

\begin{defi}
\label{def-stably-prime-exceptional}
A line bundle $L\in \Pic(X)$ is called 
{\em stably prime-exceptional}, if there exists a closed analytic subset
$Z\subset Def(X,L)$, of positive codimension, such that the linear system 
$\linsys{L_t}$ consists of a prime exceptional divisor $E_t$,
for all $t\in \left[Def(X,L)\setminus Z\right]$.
\end{defi}

Note that a stably prime-exceptional line bundle $L$ is effective, by the 
semi-continuity theorem. Furthermore, if we set $\ell:=c_1(L)$ and define
the reflection $R_\ell(\alpha):=\alpha-2\frac{(\alpha,\ell)}{(\ell,\ell)}\ell$, then
$R_\ell$ belongs to $Mon^2_{Hdg}(X)$.

\begin{rem}
Note that the linear system $\linsys{L}$, of a
stably prime-exceptional line bundle $L$,
may have positive dimension, if the Zariski decomposition of Theorem
\ref{thm-Zariski-decomposition} is non-trivial. Even if $\linsys{L}$ consists of
a single exceptional divisor, it may be reducible or non-reduced, i.e., the special point $0$
may belong to the closed analytic subset $Z$ in Definition
\ref{def-stably-prime-exceptional}.
\end{rem}

\begin{prop}
\label{prop-deformation-of-stably-prime-exceptional}
Let $E$ be a prime exceptional divisor on a 
projective irreducible holomorphic symplectic manifold $X$. 
\begin{enumerate}
\item  (\cite{markman-prime-exceptional}, Proposition 5.2)
The line bundle $\StructureSheaf{X}(E)$ is stably prime-exceptional.
\item
(\cite{markman-prime-exceptional}, Proposition 5.14)
Let $Y$ be an irreducible holomorphic symplectic manifold and
$g:H^2(X,\Integers)\rightarrow H^2(Y,\Integers)$
a parallel-transport operator, which is an isomorphism of Hodge structures.
Set $\alpha:=g([E])\in H^{1,1}(Y,\Integers)$. Then either $\alpha$
or $-\alpha$ is the class of a stably prime-exceptional line bundle.
\end{enumerate}
\end{prop}

\begin{example}
Let $X$ be a $K3$ surface. A line bundle $L$ is stably prime-exceptional,
if and only if $\deg(L)=-2$, and $(c_1(L),\kappa)>0$, for some K\"{a}hler class
$\kappa$ on $X$. 
\end{example}

Denote by $\Spe\subset H^{1,1}(X,\Integers)$ the subset of 
classes of stably prime-exceptional divisors. 

\begin{defi}
\label{def-W-Exc}
Let $W_{Exc}\subset Mon^2_{Hdg}(X)$ be the reflection subgroup generated by
$\{R_\ell \ : \ \ell\in \Spe\}$. 
\end{defi}

Note that $R_\ell=R_{-\ell}$.

\begin{cor}
\label{cor-Mon-Hdg-invariance-of-Spe}
The union $\Spe\cup -\Spe$ is a $Mon^2_{Hdg}(X)$-invariant subset of
$H^{1,1}(X,\Integers)$. 
In particular, $W_{Exc}$ is a normal subgroup of $Mon^2_{Hdg}(X)$
\end{cor}

\begin{prop}
\label{prop-fundamental-exceptional-chamber-is-indeed-an-exceptional-one}
The fundamental exceptional chamber $\FE_X$,
introduced in Definition \ref{def-fundamental-exceptional-chamber}, 
is equal to the subset
\begin{equation}
\label{eq-fundamental-exceptional-chamber-is-indeed-an-exceptional-one}
\{\alpha\in \C_X \ : \ (\alpha,\ell) > 0, \ \mbox{for every} \ \ell\in \Spe\}.
\end{equation}
\end{prop}

\begin{proof}
Denote the exceptional chamber 
(\ref{eq-fundamental-exceptional-chamber-is-indeed-an-exceptional-one}) by
$Ch_0$. Then $Ch_0\subset \FE_X$, since a prime exceptional divisor is
stably prime-exceptional, by Proposition
\ref{prop-deformation-of-stably-prime-exceptional}. 
Let $\alpha$ be a class in $\FE_X$, $\ell\in \Spe$, and 
$
\ell = P(\ell)+N(\ell)
$
its Zariski decomposition of Theorem \ref{thm-Zariski-decomposition}. 
Then $N(\ell)$ is a non-zero exceptional divisor,
since $(\ell,\ell)<0$ and $(P(\ell),P(\ell))\geq 0$.
Furthermore, $(\alpha,P(\ell))\geq 0$, since $\alpha$ and $P(\ell)$ belong to the
closure of the positive cone.
%$\overline{\BK}_X$, by Proposition \ref{prop-closure-of-FE-and-FU-are-equal}. 
Thus, $(\alpha,\ell)\geq (\alpha,N(\ell))>0$. We conclude that $\alpha$ belongs to $Ch_0$ and 
so $\FE_X\subset Ch_0$.
\end{proof}

In section \ref{sec-numerical-determination-of-BK-X} 
we will provide a numerical determination of the 
set $\Spe$, and hence of $\FE_X$, for $X$ of $K3^{[n]}$-type.

%*********************************************************************
%
%*********************************************************************
\subsection{Hyperbolic reflection groups}
\label{sec-hyperbolic-reflection-groups}
Consider the vector space $\RealNumbers^{n+1}$, endowed with the
quadratic form $q(x_0, \dots, x_n)=x_0^2-\sum_{i=1}^n x_i^2$.
We will denote the inner product space $(\RealNumbers^{n+1},q)$ by $V$
and denote by $(v,w)$, $v, w\in V$, the inner product, such that $q(v)=(v,v)$.
Let $v:=(v_0, \dots, v_n)$ be the coordinates of a vector $v$ in $V$. 
The {\em hyperbolic} (or {\em Lobachevsky}) {\em space} is
\[
\HH^n \ \ := \ \ \{v\in V \ : \ q(v)=1 \ \mbox{and} \ v_0>0\}.
\]
$\HH^n$ has two additional descriptions. It is the set of 
$\RealNumbers_{>0}$ orbits (half lines) in one of the two connected component 
of the cone $\C'_V:=\{v\in V \ : \ q(v)>0\}$. We will denote by $\C_V$ the 
chosen connected component of $\C'_V$ and refer to $\C_V$ as the
{\em positive cone}. 
$\HH^n$ also naturally embeds in $\PP^n(\RealNumbers)$ as the image of $\C_V$.
A {\em hyperplane} in $\HH^n$ is a non-empty intersection of $\HH^n$ with 
a hyperplane in $\PP^n(\RealNumbers)$.

The first description of $\HH^n$ above depended on the diagonal form of the quadratic form $q$. The
last two descriptions of $\HH^n$ produce a copy of $\HH^n$ associated more generally to any quadratic form 
$q(x_0, \dots, x_n)=\sum_{i,j=0}^n a_{ij}x_ix_j,$ $a_{ij}\in\RationalNumbers$, of signature $(1,n)$. We
will consider from now on this more general set-up. 

$\HH^n$ admits a metric of constant curvature \cite{vinberg}.
Let $O^+(V)$ be the subgroup of the isometry group of $V$
mapping $\C_V$ to itself. 
Then $O^+(V)$ acts transitively on $\HH^n$ via isometries. 
The stabilizer $Stab_{O^+(V)}(t)$, of every point $t\in \HH^n$,
is compact, since the hyperplane $t^\perp\subset V$ is negative definite.

A subgroup $\Gamma\subset O^+(V)$ is said to be a {\em discrete group of motions of}
$\HH^n$,
if for each point $t\in \HH^n$, the stabilizer $Stab_\Gamma(t)$ is finite
and the orbit $\Gamma\cdot t$ is a discrete subset of $\HH^n$.
The arithmetic group $O^+(V,\Integers)$ 
is a discrete group of motions (\cite{vinberg}, Ch. 1, section 2.2). Furthermore, if
a subgroup $\Gamma\subset O^+(V)$ is commensurable to a discrete group of motions,
then $\Gamma$ is a discrete group of motions as well
(\cite{vinberg}, Ch. 1, Proposition 1.13). 
Given a group homomorphism $\widetilde{\Gamma}\rightarrow O^+(V)$, we say that 
{\em $\widetilde{\Gamma}$ acts on $\HH^n$ via a discrete group of motions,} if its image
$\Gamma\subset O^+(V)$ is a discrete group of motions.

\begin{lem}
Let $X$ be a projective irreducible holomorphic symplectic manifold.
Then $Mon^2_{Hdg}(X)$ acts via a discrete group of motions on 
the hyperbolic space $\HH_X$ associated to $V:=H^{1,1}(X,\Integers)\otimes_\Integers\RealNumbers$
as well as on 
the hyperbolic space $\widetilde{\HH}_X$ associated to $H^{1,1}(X,\RealNumbers)$.
\end{lem}

\begin{proof}
Let $\rho$ be the rank of $\Pic(X)$. The Beauville-Bogomolov pairing restricts 
to $H^{1,1}(X,\Integers)$ as a non-degenerate pairing of signature $(1,\rho-1)$.
The action of $Mon^2_{Hdg}(X)$ on $\HH_X$ factors through the action of $O^+[H^{1,1}(X,\Integers)]$.
The latter acts as a discrete group of motions on $\HH_X$
(see \cite{vinberg}, Ch. 1, section 2.2). The statement of the lemma follows for $\HH_X$.

Let $G$ be the kernel of
the restriction homomorphism $Mon^2_{Hdg}(X)\rightarrow O^+[H^{1,1}(X,\Integers)]$.
We prove next that $G$ is a finite group.
Let $T(X)$ be the subspace of $H^2(X,\RealNumbers)$ orthogonal to $H^{1,1}(X,\Integers)$.
Set $T^{1,1}(X):=T(X)\cap H^{1,1}(X,\RealNumbers)$. 
The Beauville-Bogomolov pairing restricts to $T^{1,1}(X)$ as a negative definite pairing. 
Let $T^+(X)\subset T(X)$ be the orthogonal complement of $T^{1,1}(X)$ in $T(X)$.
Then $T^+(X)$ is the two-dimensional positive definite subspace of $T(X)$, spanned 
by the real and imaginary parts of a holomorphic $2$-form on $X$. 
%We have the $G$-invariant direct sum decomposition $T(X)=T^+(X)\oplus T^{1,1}(X)$.
$G$ acts faithfully on $T(X)$ and it embeds as a discrete subgroup of the 
compact group $O\left(T^+(X)\right)\times O\left(T^{1,1}(X)\right)$. 
We conclude that $G$ is finite. 

The linear subspace $\PP\left(T^{1,1}(X)\right)$ of $\PP\left(H^{1,1}(X,\RealNumbers)\right)$
is disjoint from $\widetilde{\HH}_X$ and so the orthogonal projection 
$H^{1,1}(X,\RealNumbers)\rightarrow V$ induces a well defined $Mon^2_{Hdg}(X)$-equivariant
map $\pi:\widetilde{\HH}_X\rightarrow \HH_X$. Explicitly, a point
$\tilde{v}$ in the positive cone of $H^{1,1}(X,\RealNumbers)$
can be uniquely decomposed as a sum $\tilde{v}=v+t$, with $v\in V$
and $t\in T^{1,1}(X)$, and $\pi$ takes the image of $\tilde{v}$ in $\widetilde{\HH}_X$ to the
image  of $v$ in $\HH_X$. 

We show next that $Mon^2_{Hdg}(X)$ acts on $\widetilde{\HH}_X$ via a discrete group of motions.
%points in $\widetilde{\HH}_X$ have finite stabilizers in $Mon^2_{Hdg}(X)$.
Set $\Gamma:=Mon^2_{Hdg}(X)/G$. 
Let $\tilde{x}$ be a point of $\widetilde{\HH}_X$ and set $x:=\pi(\tilde{x})$.
The stabilizing subgroup $Stab_\Gamma(x)$ is finite, since  $\Gamma$ acts 
on $\HH_X$ as a discrete group of motions. 
The preimage of $Stab_\Gamma(x)$ in $Mon^2_{Hdg}(X)$ is 
finite and contains the stabilizer of $\tilde{x}$ in $Mon^2_{Hdg}(X)$. 
Hence, the latter stabilizer is finite. 
%
%Let $\tilde{v}$ be a point in the positive cone of $H^{1,1}(X,\RealNumbers)$.
%Decompose it as a sum $\tilde{v}=v+t$, with $v\in V$
%and $t\in T^{1,1}(X)$. Then $v$ belongs to the positive cone of 
%$V$ and its image
%$[v]$ in $\HH_X$ has a finite 
%stabilizing subgroup $Stab_\Gamma([v])$ in $\Gamma:=Mon^2_{Hdg}(X)/G$. 
%The preimage of $Stab_\Gamma([v])$ in $Mon^2_{Hdg}(X)$ is finite and contains the stabilizer of 
%$[\tilde{v}]\in \widetilde{\HH}_X$ in $Mon^2_{Hdg}(X)$. 
%
Let $y$ be a point in the orbit $\Gamma\cdot x$ in $\HH_X$. Then 
$\pi^{-1}(y)$ intersects the orbit $Mon^2_{Hdg}(X)\cdot\tilde{x}$ in an orbit of a finite subgroup, namely, an orbit of
the preimage
of $Stab_\Gamma(y)$  in $Mon^2_{Hdg}(X)$. 
The orbit $Mon^2_{Hdg}(X)\cdot\tilde{x}$ is a discrete subset of $\widetilde{\HH}_X$, since $\pi$
restricts to it as a finite map onto the discrete orbit of $x$ in $\HH_X$.
\end{proof}

Given an element $e\in V$, with $q(e)<0$, we get the reflection
$R_e\in O^+(V)$, given by $R_e(w)=w-2\frac{(e,w)}{(e,e)}e$.

\begin{defi}
\label{def-hyperbolic-reflection-group}
A {\em hyperbolic reflection group} is a discrete group of motions of $\HH^n$ generated
by reflections.
\end{defi}

Given a vector $e\in V$, with $q(e)<0$, set 
\[
H_e^+ \ := \ \{v\in \C_V \ : \ (v,e)>0\}/\RealNumbers_{>0}.
\]
Define $H_e^-$ similarly using the inequality $(v,e)<0$.
Set $H_e:=e^\perp\cap \HH^n$, where $e^\perp$ is the hyperplane of 
$\PP(V)$ orthogonal to $e$. 
Then $\HH^n\setminus H_e$ is the disjoint union of its two connected components 
$H_e^+$ and $H_e^-$.
The closures $\overline{H}_e^{\pm}$ are called {\em half-spaces}. 

\begin{defi}
\label{def-generalized-convex-polyhedron}
\begin{enumerate}
\item
A set $\{\Sigma_i \ : i\in I\}$, of subsets of a topological space $X$,
is {\em locally finite}, if each point $x\in X$ has an open neighborhood $U_x$,
such that the intersection $\Sigma_i\cap U_x$ is empty, for all but
finitely many indices $i\in I$.
\item
A {\em decomposition} of $\HH^n$ is a locally finite covering of $\HH^n$
by closures of open connected subsets, no two of which have common interior points.
\item
A closure $D$ of an open subset of $\HH^n$ is said to be a 
{\em fundamental domain} of a discrete group of motions $\Gamma$,
if $\{\gamma(D) \ : \ \gamma\in \Gamma\}$ is a decomposition of $\HH^n$.
\item
(\cite{AVS}, Ch. 1, Definition 3.9)
A {\em convex polyhedron} is an intersection of finitely many half-spaces, having
a non-empty interior.
\item
(\cite{vinberg}, Ch 1, Definition 1.9)
A closed subset $P\subset \HH^n$ is a {\em generalized convex polyhedron}, if
$P$ is the closure of an open subset, and the intersection of $P$ with 
every bounded convex polyhedron, containing at least one common interior point,
is a convex polyhedron. 
\item
A closed cone in $\C_V$ is a {\em generalized convex polyhedron}, if
its image in $\HH^n$ is a generalized convex polyhedron.
\item
\label{def-item-rational-polydedral-cone}
A closed cone $\Pi$ in $\C_V$ is a {\em rational convex polyhedron}, if its 
image in $\HH^n$ is a convex polyhedron,
which is the intersection of finitely many half spaces $H^+_e$
with $e\in \RationalNumbers^{n+1}$.
\end{enumerate}
\end{defi}

\begin{thm} 
\label{thm-fundamental-domain-of-a-discrete-group-of-motions}
\begin{enumerate}
\item
(\cite{vinberg}, Ch. 1 Theorem 1.11)
Any discrete group of motions of $\HH^n$ has a fundamental domain, 
which is a generalized convex polyhedron.
\item
\label{thm-item-polyhedron-fundamental-domain-of-an-arithmetic-subgroup}
(\cite{vinberg}, Ch. 2 Theorem 2.5) The action on $\HH^n$ 
of any arithmetic subgroup of $O^+(V)$ has a fundamental domain, 
which is a convex polyhedron.
\end{enumerate}
\end{thm}

The decomposition of $\HH^n$, induced by translates of the fundamental domain 
in Theorem \ref{thm-fundamental-domain-of-a-discrete-group-of-motions}, is not canonical
in general. A canonical decomposition exists, 
if the discrete group of motions is a reflection group.
The hyperplanes of $n-1$ dimensional faces of a generalized convex polyhedron
are called its {\em walls}.

Let $\Gamma$ be a hyperbolic reflection group and 
$\R_\Gamma\subset \Gamma$ the subset of reflections.
Given a reflection $\rho\in\R_\Gamma$, let $H_\rho\subset \HH^n$ be the hyperplane 
fixed by $\rho$. 
Connected components of 
$\HH^n\setminus \bigcup_{\rho\in \R_\Gamma} H_\rho$ are called {\em chambers}.

\begin{thm}
\label{thm-each-chamber-is-a-generalized-convex-polyhedron}
(\cite{vinberg}, Ch. 5 Theorem 1.2 and Proposition 1.4)
\begin{enumerate}
\item
The closure of each chamber of $\Gamma$ in $\HH^n$ 
%connected component of $\HH^n\setminus \bigcup_{\rho\in \R_\Gamma} H_\rho$ 
is a generalized convex polyhedron,\footnote{This polyhedron is moreover
a generalized Coxeter polyhedron (\cite{vinberg}, Ch. 5 Definition 1.1),
but we will not use this fact.}
which is a fundamental domain for $\Gamma$.
\item
%Let $Pol$ be the closure of a chamber of $\Gamma$. 
$\Gamma$ is generated by reflections in the walls of any of its chambers in $\HH^n$.
\end{enumerate}
\end{thm}

Let $\Gamma$ be any discrete group of motions of $\HH^n$. 
Denote by $\Gamma_r$ the subgroup of $\Gamma$ generated by all 
reflections in $\Gamma$. We call $\Gamma_r$ the
{\em reflection subgroup of $\Gamma$.}
Choose a chamber $D$ of $\Gamma_r$. 
Let $\Gamma_D\subset \Gamma$ be the subgroup 
$\{\gamma\in \Gamma \ : \ \gamma(D)=D\}$.

\begin{thm} 
\label{thm-any-discrete-group-of-motions-is-a-semi-direct-product}
(\cite{vinberg}, Ch. 5 Proposition 1.5)
$\Gamma_r$ is a normal subgroup of $\Gamma$, and  
$\Gamma$ is the semi-direct product of $\Gamma_r$ and $\Gamma_D$.
\end{thm}

We refer the reader to the book \cite{vinberg} and the interesting recent survey
\cite{dolgachev-reflections} for detailed expositions of the subject of
hyperbolic reflection groups.

Let $X$ be a projective irreducible holomorphic symplectic manifold.

\begin{thm}
\label{thm-numerical-characterization-of-BK-X}
The fundamental exceptional chamber $\FE_X$, introduced in
Definition \ref{def-fundamental-exceptional-chamber},
is equal to the connected component of
\begin{equation}
\label{eq-positive-cone-minus-walls}
\C_X \ \setminus \  \bigcup \{\ell^\perp \ : \ \ell\in \Spe\}
\end{equation}
containing the K\"{a}hler cone. In particular, $\FE_X$ is the interior of a  
generalized convex polyhedron 
(Definition \ref{def-generalized-convex-polyhedron}).
\end{thm}

\begin{proof}
The group $W_{Exc}$ is a hyperbolic reflection group and 
the set $U$ in equation (\ref{eq-positive-cone-minus-walls}) is an open subset 
of $\C_X$, which is the union of the interiors of the fundamental chambers of the $W_{Exc}$-action on $\C_X$, 
by Theorem \ref{thm-each-chamber-is-a-generalized-convex-polyhedron}.
The intersection of $\FE_X$ and  $U$ is the union of connected components of $U$,
by the definitions of $\FE_X$ and $W_{Exc}$. 
$\FE_X$ is contained in $U$, by Proposition 
\ref{prop-fundamental-exceptional-chamber-is-indeed-an-exceptional-one}.
$\FE_X$ is convex cone, hence a connected component of $U$. 
$\FE_X$ contains $\K_X$, by the definition of $\FE_X$.
%$\FE_X$ is a fundamental chamber of $W_{Exc}$, by Proposition 
%\ref{prop-fundamental-exceptional-chamber-is-indeed-an-exceptional-one}.
%The group $W_{Exc}$ is a hyperbolic reflection group, hence
%$\FE_X$ is the interior of a generalized convex polyhedron,
%by Theorem \ref{thm-each-chamber-is-a-generalized-convex-polyhedron}. 
\end{proof}

%*********************************************************************
%
%*********************************************************************
\subsection{$Mon^2_{Hdg}(X)$ is a semi-direct product of $W_{Exc}$ and $Mon^2_{Bir}(X)$}
\label{sec-Mon-2-Hdg-is-semi-direct-product}
Denote by $\Pex$ the set of prime exceptional divisors in $X$.
Given  $E\in \Pex$, denote by $R_E$ the corresponding 
reflection (Proposition \ref{prop-reflection-by-exceptional-divisor}).

\begin{thm}
\label{thm-Mon-2-Hdg-is-a-semi-direct-product}
\begin{enumerate}
\item
\label{thm-item-W-acts-simply-transitively}
The group $Mon^2_{Hdg}(X)$ acts transitively on the set of exceptional chambers, introduced in Definition
\ref{def-Kahler-type-chambers},  
and the subgroup $W_{Exc}$ acts simply-transitively on this set.
\item
\label{thm-item-exceptional-chambers-are-connected-components}
The exceptional chambers  are precisely the 
connected component of 
the open set in equation (\ref{eq-positive-cone-minus-walls}), i.e., 
each exceptional chamber is the interior of a fundamental domain of the $W_{Exc}$
action on $\C_X$.
\item
\label{thm-item-W-is-generated-by-exceptional-reflections}
The group $W_{Exc}$ is generated by $\{R_e \ : \ e\in \Pex\}$.
%\item
%\label{thm-item-locally-finite}
%The set of hyperplane sections $\{\C_X\cap e^\perp\}_{e\in \Pex}$
%is locally finite (every compact subset of $\C_X$ intersects only finitely many 
%such hyperplanes).
\item
\label{thm-item-stabilizer-of-FE}
The subgroup of $Mon^2_{Hdg}(X)$ stabilizing the fundamental exceptional chamber 
$\FE_X$ is equal to $Mon^2_{Bir}(X)$.
\item
\label{thm-item-semi-direct-product}
$Mon^2_{Hdg}(X)$ is the semi-direct product of its subgroups 
$W_{Exc}$ and $Mon^2_{Bir}(X)$.
\end{enumerate}
\end{thm}

When $X$ is a $K3$ surface $Mon^2_{Hdg}(X)$ is equal to the group of Hodge isometries of 
$H^2(X,\Integers)$ preserving the spinor norm and 
$Mon^2_{Bir}(X)$ is equal to the group of biregular automorphisms of $X$.
Furthermore, the fundamental exceptional chamber is equal to the K\"{a}hler cone
of the $K3$ surface. 
Theorem \ref{thm-Mon-2-Hdg-is-a-semi-direct-product} 
is well known in the case of $K3$ surfaces
\cite{burns-rapoport,PS}, or (\cite{looijenga-peters}, Proposition 1.9).

\begin{proof}
Parts (\ref{thm-item-W-acts-simply-transitively}) and (\ref{thm-item-exceptional-chambers-are-connected-components}):
$Mon^2_{Hdg}(X)$ acts transitively on the set of exceptional chambers, by their definition.
The subgroup $W_{Exc}$ acts simply-transitively on the set of connected components of
the set $U$ in equation (\ref{eq-positive-cone-minus-walls}), by Theorem 
\ref{thm-each-chamber-is-a-generalized-convex-polyhedron}. 
One of these is $\FE_X$,
by Theorem \ref{thm-numerical-characterization-of-BK-X}.
Hence, every connected component of $U$ is an exceptional chamber. 
$Mon^2_{Hdg}(X)$ acts on the set of   connected component of $U$, by
Corollary \ref{cor-Mon-Hdg-invariance-of-Spe}. Hence, every exceptional chamber is 
a connected component of $U$.

Part (\ref{thm-item-W-is-generated-by-exceptional-reflections}):
The walls in the boundary of the fundamental exceptional chamber
are all of the form $[E]^\perp\cap \C_X$, for some
prime exceptional divisor $E$, by definition. $\FE_X$ is the interior
of a chamber of $W_{Exc}$, by Theorem \ref{thm-numerical-characterization-of-BK-X}.
We conclude that $W_{Exc}$ is generated by $\{R_e \ : \ e\in\Pex\}$, by 
Theorem \ref{thm-each-chamber-is-a-generalized-convex-polyhedron}.
%Let $W_1\subset W_{Exc}$ be the subgroup generated by $\{R_e \ : \ e\in\Pex\}$. Then 
%$W_1=W_{Exc}$ and $W_{Exc}$ 
%acts simply transitively on the set of exceptional chambers, by
%Theorem \ref{thm-each-chamber-is-a-generalized-convex-polyhedron}.

%Part (\ref{thm-item-exceptional-chambers-are-connected-components}) 
%The fundamental exceptional chamber $\FE_X$ is an open subset of $\C_X$,
%by Theorem \ref{thm-numerical-characterization-of-BK-X}.  
%$\FE_X$  is a closed subset of $\C_X  \setminus \ \bigcup \{\ell^\perp \ : \ \ell\in \Spe\}$,
%by definition. 
%$\FE_X$ is a convex cone, hence a connected component
%of $\C_X  \setminus \ \bigcup \{\ell^\perp \ : \ \ell\in \Spe\}$.
%It suffices to show that $W_{Exc}$ acts transitively both on the set of exceptional chambers and 
%on the set $\Sigma$ of  connected components of $\C_X  \setminus \ \bigcup \{\ell^\perp \ : \ \ell\in \Spe\}$.
%The set $\Sigma$ 
%is the set of inverse images of the set of chambers of $W_{Exc}$ in the hyperbolic space 
%associated to $\C_X$. $W_{Exc}$ acts transitively on $\Sigma$, 
%by Theorem \ref{thm-each-chamber-is-a-generalized-convex-polyhedron}. 
%$W_{Exc}$ acts transitively on the set of exceptional chambers, by part
%(\ref{thm-item-W-acts-simply-transitively}). 
Part (\ref{thm-item-stabilizer-of-FE}): 
$Mon^2_{Bir}(X)$ is the subgroup of $Mon^2_{Hdg}(X)$ leaving $\FE_X$ invariant,
by Lemma \ref{lemma-on-semi-chambers}
part \ref{lemma-item-stabilizer-of-fundamental-uniruled-chamber}.

Part (\ref{thm-item-semi-direct-product}):
$Mon^2_{Hdg}(X)$ is generated by $W_{Exc}$ and $Mon^2_{Bir}(X)$,
by parts (\ref{thm-item-W-acts-simply-transitively}) and  (\ref{thm-item-stabilizer-of-FE}).
The intersection $W_{Exc}\cap Mon^2_{Bir}(X)$ is trivial, since the action of 
$W_{Exc}$ on the set of exceptional chambers is free. $W_{Exc}$ is a normal subgroup of 
$Mon^2_{Hdg}(X)$, by Corollary \ref{cor-Mon-Hdg-invariance-of-Spe}.
\end{proof}

\begin{caution}
When $X$ is a $K3$ surface, then $W_{Exc}$ is the reflection subgroup of 
$Mon^2_{Hdg}(X)$, i.e., every reflection $g\in Mon^2_{Hdg}(X)$ is of the form
$R_\ell$, for a class $\ell$ satisfying $(\ell,\ell)=-2$. This follows easily from the fact that
$H^2(X,\Integers)$ is a unimodular lattice.
$W_{Exc}$ may be strictly smaller than the reflection subgroup of $Mon^2_{Hdg}(X)$,
for a higher dimensional irreducible holomorphic symplectic manifold $X$.
In other words, there are examples of elements $\alpha\in H^{1,1}(X,\Integers)$, 
with $(\alpha,\alpha)<0$, such that $R_\alpha$ belongs to $Mon^2_{Hdg}(X)$,
but neither $\alpha$, nor $-\alpha$, belongs to $\Spe$. Instead,
$R_\alpha$ is induced by a bimeromorphic map from $X$ to itself 
(see Example \ref{example-monodromy-reflective-but-not-spe} below, 
and section 11 of \cite{markman-prime-exceptional} for additional examples).
\end{caution}

Let $L$ be a stably prime-exceptional line bundle and set $\ell:=c_1(L)$. 
The hyperplane $\ell^\perp$ intersects $\overline{\FE}_X$ in a top dimensional cone
in $\ell^\perp$, only if $L=\StructureSheaf{X}(E)$ 
for some prime exceptional divisor $E$,
by Proposition \ref{prop-fundamental-exceptional-chamber-is-indeed-an-exceptional-one}.
We show next that the condition is also sufficient.

\begin{lem}
\label{lemma-faces-of-FE}
Let $E$ be a prime exceptional divisor on $X$. 
Then $E^\perp\cap \overline{\FE}_X$ is a top dimensional cone in 
the hyperplane $E^\perp$.
Consequently, $W_{Exc}$ can not be generated by 
any proper subset of $\{R_e \ : \ e \in\Pex\}$.
\end{lem}

\begin{proof}
Let $e$ be an element of $\Pex$.
It suffices to show that $e^\perp\cap\overline{\FE}_X\cap\C_X$ contains elements, which are not orthogonal
to any other $e'\in \Pex$. Choose $x\in\FE_X$ and set 
$y:=x-\frac{(x,e)}{(e,e)}e$. Then $(y,e)=0$. 
Given $e'\in\Pex$, $e'\neq e$, then $(e,e')\geq 0$ and $(x,e')>0$. 
Now $(e,e)<0$. 
We get the following inequalities.
\begin{eqnarray*}
(e',y)&=&(e',x)-\frac{(x,e)}{(e,e)}(e',e) \ > \ 0.
\\
(y,y) & = & (x,x)-\frac{(x,e)^2}{(e,e)} \ > \ 0.
\end{eqnarray*}
%The first equality above implies that $(e',y)>0$, 
%since $(e,e)<0$. The second equality implies that $(y,y)>0$, provided we
%choose $x$, such that $(x,e)$ is sufficiently small. 
We conclude that $y$ belongs to $e^\perp\cap \overline{\FE}_X\cap\C_X$, and $y$ does not belong to
$(e')^\perp$, for any $e'\in \Pex\setminus\{e\}$.
\end{proof}

%*******************
% Hide
%*******************
\hide{
\begin{proof}
We use the notation of Proposition \ref{prop-druel}.
Given a birational morphism $f:X\rightarrow X'$ onto an irreducible holomorphic
symplectic manifold $X'$,
we have $f_*(\FE_X)=\FE_{X'}$ and the strict transform $E'$ of $E$ 
is a prime exceptional divisor. Hence,
it suffices to prove the lemma in the case $X=X'$, where $X'$ is the one in  
Proposition \ref{prop-druel}, and $\pi:X\rightarrow Y$ is a regular morphism. 
Set $NS(X)_\RealNumbers:=H^{1,1}(X,\Integers)\otimes_\Integers\RealNumbers$.
Set $\T\C_X:=\C_X\cap \left[NS(X)_\RealNumbers^\perp\right]$.
$\T\C_X$ is the positive cone associated to the transcendental lattice.
Set $\A_X:=\overline{\FE}_X\cap NS(X)_\RealNumbers$. 
Then the sum  of the two cones $\T\C_X+\A_X$ 
is a subcone of $\overline{\FE}_X$, and
$\dim\left(\T\C_X+\A_X\right)=\dim(\T\C_X)+\dim(\A_X)$. 
It suffices to show that the intersection $\A_X\cap E^\perp$ is a top dimensional cone in
$NS(X)_\RealNumbers\cap E^\perp$. 

The {\em dissident locus} $\Sigma_0$ in $Y$ is the locus in $Y_{sing}$,
where the singularities of $Y$ are not of ADE type.
$\Sigma_0$ is a closed subscheme of co-dimension $\geq 4$ in $Y$
(\cite{namikawa-deformations}, Proposition 1.6).
Set $U:=Y\setminus\Sigma_0$ and $\widetilde{U}:=\pi^*{U}\subset X$.
We get the injective homomorphism
\[
H^2(U,\Integers)\LongRightArrowOf{\pi^*} H^2(\widetilde{U},\Integers)
\cong H^2(X,\Integers).
\]
The image of $H^2(U,\Integers)$ in $H^2(X,\Integers)$ is equal to $E^\perp$,
by (\cite{markman-galois}, Lemma 4.2). 
Now $\A_X$ contains the image of the ample cone of $U$, and 
the latter is a top dimensional cone in 
$\pi^*NS(U)_\RealNumbers=E^\perp\cap NS(X)_\RealNumbers$.
%since the class $[E]^\vee$, given in (\ref{eq-E-vee}), is equal to the class of a generic fiber of
%$E\rightarrow \pi(E)\subset Y$
%by Proposition \ref{prop-reflection-by-exceptional-divisor}.
\end{proof}
%**************
% End hide
%**************
}

\noindent
{\bf Proof of Theorem \ref{thm-decomposition-of-a-Hodge-isometry-as-wg}:}
$W_{Exc}$ is a normal subgroup of $Mon^2_{Hdg}(X)$, 
by Corollary \ref{cor-Mon-Hdg-invariance-of-Spe}.
There exists a bimeromorphic map $h:X_1\rightarrow X_2$, 
by Theorem \ref{thm-Hodge-theoretic-Torelli}, and $h^*$ is a parallel transport operator,
by Theorem \ref{thm-bimeromorphic-implies-inseparable}. 
The composition $f\circ h^*$ belongs to 
$Mon^2_{Hdg}(X_2)$. There exists an element $w$ of $W_{Exc}(X_2)$, such that 
$w^{-1}f\circ h^*$ belongs to $Mon^2_{Bir}(X_2)$, 
by Theorem \ref{thm-Mon-2-Hdg-is-a-semi-direct-product}.
Let $\varphi:X_2\rightarrow X_2$ be a bimeromorphic map, such that
$\varphi_*=w^{-1}f\circ h^*$. Then $f=w(\varphi h)_*$.
Set $g:=\varphi h$ to obtain the desired decomposition $f=w\circ g_*$. 
%The element $w^{-1}$ can be written as a product 
%$R_{E_1}\cdots R_{E_k}$, for some prime exceptional divisors $E_i$,
%by Theorem \ref{thm-Mon-2-Hdg-is-a-semi-direct-product}.

Assume that $\tilde{g}:X_1\rightarrow X_2$ is a birational map and $\tilde{w}$
is an element of $W_{Exc}(X_2)$, such that $f=\tilde{w}\tilde{g}_*$.
Then $w^{-1}\tilde{w}=(\tilde{g}^{-1}g)_*$ belongs to the intersection of
$W_{Exc}(X_2)$ and $Mon^2_{Bir}(X_2)$, which is trivial, by Theorem 
\ref{thm-Mon-2-Hdg-is-a-semi-direct-product}. Thus, 
$w=\tilde{w}$ and $g_*=\tilde{g}_*$. Now,
$\tilde{g}=g(g^{-1}\tilde{g})$, and $g^{-1}\tilde{g}$ is a birational map 
inducing the identity on $H^2(X_1,\Integers)$. In particular, $g^{-1}\tilde{g}$
maps $\K_{X_1}$ to itself, and hence is a biregular automorphism.
\EndProof

%*****************
% Hide
%*****************
\hide{
%*********************************************************************
%
%*********************************************************************
\subsection{Moduli space of bimeromorphic classes?}

Let $X_0$ be an irreducible holomorphic symplectic manifold satisfying the following
property. 

\begin{assumption}
\label{assumption-druel-prop-holds}
For every irreducible holomorphic symplectic manifold $X$, deformation equivalent to $X_0$,
but not necessarily projective, and for every prime exceptional divisor $E\subset X$,
the analogue of Proposition \ref{prop-druel} holds.
\end{assumption}

Theorem \ref{thm-Mon-2-Hdg-is-a-semi-direct-product}
holds if we replace the projectivity assumption by assumption 
\ref{assumption-druel-prop-holds} on $X$. The Torelli theorem 
for K\"{a}hler $K3$ surfaces admits a modular formulation
in terms of the moduli space of isomorphism classes of pairs $(X,\kappa)$,
where $X$ is a K\"{a}hler $K3$ surface and $\kappa\in\K_X$ is a 
K\"{a}hler class (\cite{looijenga-peters}, Theorem 10.5) and
(\cite{BHPV}, Ch. VIII, Theorem (12.3)). 
We work out a higher dimensional analogue for deformation types satisfying 
the additional assumption \ref{assumption-druel-prop-holds}.

Let $\C_{\Omega_\Lambda}\subset \Omega_\Lambda\times \Lambda_\RealNumbers$
be the subset
\[
\{(\sigma,\lambda) \ : \  (\sigma, \lambda)=0, \  \mbox{and} \ (\lambda,\lambda)>0\}.
\]
Fix a connected component  $\FM^0_\Lambda$ of the moduli space of marked pairs
and denote by $Mon^2(\FM^0_\Lambda)$ its monodromy group
(Definition \ref{def-polarized-monodromy-subgroup-of-O-Lambda}).
Let $\Spe_\Lambda\subset\Lambda$ be the set of all $Mon^2(\FM^0_\Lambda)$-orbits of 
$\eta(\ell)$, where $(X,\eta)$ is a marked pair in $\FM^0_\Lambda$,
and $\ell$ is a stably prime-exceptional class in $H^{1,1}(X,\Integers)$.
Set $\ell^\perp_{\ComplexNumbers}:=\Omega_\Lambda\cap \ell^\perp$
and $\ell^\perp_\RealNumbers\subset \Lambda_\RealNumbers$ the 
hyperplane orthogonal to $\ell$. Set
\[
\C_{\Omega_\Lambda}^0 \ := \ 
\C_{\Omega_\Lambda}\setminus\bigcup_{\ell\in\Spe_\Lambda}
\left[\ell^\perp_{\ComplexNumbers}\times \ell^\perp_\RealNumbers\right].
\]

The group $Mon^2(\FM^0_\Lambda)$ acts diagonally on 
$\Omega_\Lambda\times \Lambda_\RealNumbers$, and the subset 
$\C_{\Omega_\Lambda}^0$ is $Mon^2(\FM^0_\Lambda)$-invariant.
The stabilizer $Mon^2_{Hdg}(\sigma)$
in $Mon^2(\FM^0_\Lambda)$, of a point $\sigma\in \Omega_\Lambda$,  
is equal to 
$\eta\circ Mon^2_{Hdg}(X)\circ \eta^{-1}$, for any marked pair $(X,\eta)\in \FM^0_\Lambda$
with period $\sigma$. Let $\C^0_\sigma$
be the fiber of $\C^0_{\Omega_\Lambda}$ over $\sigma$. 
Let $\C^0_X$ be the union of all exceptional chambers in $\C_X$.
Then 
\[
\C^0_\sigma/Mon^2_{Hdg}(\sigma)\cong
\C^0_X/Mon^2_{Hdg}(X)\cong
\FE_X/Mon^2_{Bir}(X).
\]
%******************
% End Hide
%******************
}
%*********************************************************************
%
%*********************************************************************
\subsection{Morrison's movable cone conjecture}
\label{sec-finiteness}
Let $X$ be a projective irreducible holomorphic symplectic manifold.
We describe first an analogy between results on the ample cone of a projective $K3$ surface
and results on the movable cone of $X$. 
Set $\NS:=H^{1,1}(X,\Integers)$,
$\NS_\RealNumbers:=\NS\otimes_\Integers\RealNumbers$, and 
$\NS_\RationalNumbers:=\NS\otimes_\Integers\RationalNumbers$. 
Let $\C_{\NS}$ be the intersection $\C_X\cap \NS_\RealNumbers$.

\begin{defi}
\label{def-movable-cone}
\begin{enumerate}
\item
A line bundle $L$ on $X$ is {\em movable}, if the base locus of the linear system $\linsys{L}$
has codimension $\geq 2$.
\item
The {\em movable cone} $\MV_X$ is the convex hull in $\NS_\RealNumbers$ of all classes of
movable line bundles. 
\end{enumerate}
\end{defi}

Let $\MV_X^0$ be the interior of $\MV_X$ and $\overline{\MV}_X$ the closure of $\MV_X$ in
$\NS_\RealNumbers$. 

\begin{lem}
\label{lemma-movable-cone-is-a-fundamental-domain}
The equality $\MV_X^0=\FE_X\cap \NS_\RealNumbers$ holds. 
$W_{Exc}$ acts faithfully on $\C_{\NS}$ and 
the map $Ch\mapsto Ch\cap \NS_\RealNumbers$ induces a one-to-one correspondence between 
the set of exceptional chambers and the chambers in $\C_{\NS}$ of the $W_{Exc}$ action.
In particular, the closure of $\MV_X$ in $\C_{\NS}$ is a fundamental domain for the action of $W_{Exc}$ on
$\C_{\NS}$.
\end{lem}

\begin{proof}
The equality $\MV_X^0=\FE_X\cap \NS_\RealNumbers$  follows immediately from the
Zariski decomposition (Theorem \ref{thm-Zariski-decomposition}).
%and the equality $\FE_X=\FU_X$ (Proposition \ref{prop-closure-of-FE-and-FU-are-equal}). 
The set $\Spe$ is contained in $\NS$, hence the $W_{Exc}$ action on $\C_{\NS}$ is faithful
and the map $Ch\mapsto Ch\cap \NS_\RealNumbers$ induces a bijection.
\end{proof}

Let $\rho:Mon^2_{Hdg}(X)\rightarrow O(\NS)$ be the restriction homomorphism. We denote 
$\rho(W_{Exc})$ by $W_{Exc}$ as well.

\begin{lem}
\label{lemma-kernel-and-image-of-restriction-homo}
\begin{enumerate}
\item
\label{lemma-item-image-of-rho-has-finite-index}
The image $\Gamma$ 
%$\overline{M}on^2_{Hdg}(X)$ 
of $\rho$ is a finite index subgroup of $O^+(\NS)$.
\item
\label{lemma-item-ker-rho-in-Mon-2-Bir}
The kernel of $\rho$ is a subgroup of $Mon^2_{Bir}(X)$.
\item
\label{lemma-item-overline-Mon-Hdg-a-semi-direct-product}
$\Gamma$ is a semi-direct product of its normal subgroup $W_{Exc}$ and the quotient group
$\Gamma_{Bir}:=Mon^2_{Bir}(X)/\ker(\rho)$.
\end{enumerate}
\end{lem}

\begin{proof}
(\ref{lemma-item-image-of-rho-has-finite-index}) 
The positive cone $\C_X$ is $Mon^2_{Hdg}(X)$-invariant and $\C_\NS=\C_X\cap \NS$ is thus
$\Gamma$-invariant. Hence, $\Gamma$ is a subgroup of $O^+(\NS)$.
Let $O^+_{Hdg}\left(H^2(X,\Integers)\right)$ be the subgroup of $O^+\left(H^2(X,\Integers)\right)$
preserving the Hodge structure. Then $O^+_{Hdg}\left(H^2(X,\Integers)\right)$ maps onto
a finite index subgroup of $O^+(\NS)$. 
The index of $Mon^2(X)$ in $O^+[H^2(X,\Integers)]$ is finite, by a result of Sullivan \cite{sullivan}
(see also \cite{verbitsky}, Theorem 3.4). Hence, $Mon^2_{Hdg}(X)$ is a finite index subgroup of 
$O^+_{Hdg}\left(H^2(X,\Integers)\right)$. Part (\ref{lemma-item-image-of-rho-has-finite-index}) follows.

(\ref{lemma-item-ker-rho-in-Mon-2-Bir}) 
Let $g$ be an element of $\ker(\rho)$. Then $g$ acts trivially on $\Spe$. Hence, $g$ maps
$\FE_X$ to itself. It follows that $g$ belongs to $Mon^2_{Bir}(X)$, by 
Theorem \ref{thm-Mon-2-Hdg-is-a-semi-direct-product} part \ref{thm-item-stabilizer-of-FE}.

Part (\ref{lemma-item-overline-Mon-Hdg-a-semi-direct-product})
is an immediate consequence of part (\ref{lemma-item-ker-rho-in-Mon-2-Bir}) and 
Theorem \ref{thm-Mon-2-Hdg-is-a-semi-direct-product} part \ref{thm-item-semi-direct-product}.
\end{proof}

%Let $T(X)$ be the subspace  of $H^2(X,\RealNumbers)$ orthogonal to $\NS$, 
%$T^{1,1}(X):=T(X)\cap H^{1,1}(X,\RealNumbers)$, and $T^+(X)\subset T(X)$ the orthogonal complement of 
%$T^{1,1}(X)$. Any element of $Mon^2_{Hdg}(X)$ preserves the orientation of $\widetilde{\C}_X$ and 

Let $\Eff_X\subset \NS_\RealNumbers$ be the convex cone generated by classes of effective divisors on $X$.
Set $\MV_X^e:=\overline{\MV}_X\cap \Eff_X$. 
Following is Morrison's movable cone conjecture. 

\begin{conj} \cite{morrison-compactifications,morison-cone-conj,kawamata-cone-conj}
There exists a rational convex polyhedral cone (Definition \ref{def-generalized-convex-polyhedron} 
part \ref{def-item-rational-polydedral-cone}) 
$\Pi$, which is a fundamental domain for the action of
$\Bir(X)$ on $\MV_X^e$.
\end{conj}

Morrison formulated a version of the conjecture for the ample cone as well. 
The two versions coincide in dimension $2$ and for abelian varieties.
The  $K3$ surface case of the conjecture is proven by Looijenga and Sterk (\cite{sterk}, Lemma 2.4),
the Enriques surfaces case by Namikawa (\cite{namikawa-yukihiko}, Theorem 1.4), the case of abelian and 
hyperelliptic surfaces by Kawamata (\cite{kawamata-cone-conj}, Theorem 2.1),  
the case of two-dimensional Calabi-Yau pairs by Totaro \cite{totaro}, and 
the case of abelian varieties by Prendergast-Smith \cite{prendergast}.
A version of the conjectures 
for fiber spaces was formulated by Kawamata and proven in dimension $3$ in \cite{kawamata-cone-conj}.

The following theorem is a weaker version of  Morrison's movable cone conjecture, in the special case of 
projective irreducible holomorphic symplectic manifolds.
Let $\MV_X^+$ be the convex hull of $\overline{\MV}_X\cap \NS_\RationalNumbers$. 
Clearly, $\MV_X^0$ is equal to the interior of both 
$\MV_X^+$ and $\MV_X^e$.
When $X$ is a $K3$ surface the equality
$\MV_X^+=\MV_X^e$ holds. In the $K3$ case the inclusion $\MV_X^+\subset \MV_X^e$ follows from 
(\cite{BHPV}, Proposition 3.6 part i) and the inclusion $\MV_X^+\supset \MV_X^e$ is proven 
in (\cite{kawamata-cone-conj}, Proposition 2.4). 

\begin{thm}
\label{thm-movable-conj-conj}
There exists a rational convex polyhedral cone  $\Pi$ in $\MV_X^+$, such that 
$\Pi$ is a fundamental domain for the action of $\Gamma_{Bir}$ on $\MV_X^+$.
\end{thm}

\begin{proof}
The proof is identical to that of Lemma 2.4 in \cite{sterk}, which proves the $K3$-surface case of the Theorem.
When $X$ is a $K3$ surface, $\MV_X^0$ is the ample cone and $\Pex$ is the set of nodal $-2$ classes. 
The proof is lattice theoretic. Following is the dictionary translating our notation to that of Sterk.

\smallskip
\noindent
\begin{tabular}{|c|c|c|c|c|c|c|c|c|}  \hline 
\hspace{1ex}
Our notation & $\MV_X^0$ & $\C_\NS$ &  $\MV_X^+$ & $\Pex$ & $\Spe$ & $\Gamma$ & $\Gamma_{Bir}$ & $W_{Exc}$
\\
\hline
Sterk's notation & $K$ & $\C$ & $\overline{K}\cap \C_+$ &  $B$ & $\Delta^+$ & $\Gamma$ & $\Gamma_B$ & $W$
\\
\hline
\end{tabular}\\
One slight inaccuracy in the above dictionary 
is that Sterk chose $\Gamma$ to be the subgroup of $O^+\left(H^2(X,\Integers)\right)$
acting trivially on the transcendental lattice $\NS^\perp$, while we consider (in case $X$ is a $K3$ surface)
the image of $O^+_{Hdg}\left(H^2(X,\Integers)\right)$ in $O^+(\NS)$. So Sterk's $\Gamma$ is the finite index 
subgroup of our $\Gamma$ acting trivially on the finite discriminant group $\NS^*/\NS$.
Both choices satisfy the 
following complete list of assertions needed for the Looijenga-Sterk argument (in Sterk's notation).
\begin{enumerate}
\item
\label{input-item-arthmetic-group}
$\NS$ is a lattice of signature $(1,*)$ and 
$\Gamma$ is an arithmetic subgroup of $O^+(\NS)$.
\item
\label{input-item-W-is-generated-by-reflections-in-B}
$W\subset O^+(\NS)$ is the reflection group generated by reflections in elements of $B\subset \NS$.
\item
\label{input-item-Gamma-B-is-the-stabilizer-of-B}
$\Gamma_B$ is equal to the subgroup $\{g\in\Gamma \ : \ g(B)=B\}$.
\item
\label{input-item-Gamma-is-a-semi-direct-product}
$W$ is a normal subgroup of $\Gamma$ and 
$\Gamma=\Gamma_B\cdot W$ is a semi-direct product decomposition.
\item
\label{input-item-K-is-fundamental-domain}
$\overline{K}\cap \C$ is a 
fundamental domain for the action of $W$ on $\C$, cut out by closed half-spaces associated to elements of $B$. 
\end{enumerate}

Assertion (\ref{input-item-arthmetic-group}) is verified in our case in
Lemma \ref{lemma-kernel-and-image-of-restriction-homo} part \ref{lemma-item-image-of-rho-has-finite-index}.
Assertion (\ref{input-item-W-is-generated-by-reflections-in-B}) is verified in 
Theorem \ref{thm-Mon-2-Hdg-is-a-semi-direct-product} part \ref{thm-item-W-is-generated-by-exceptional-reflections}.
$Mon^2_{Bir}(X)=\{g\in Mon^2_{Hdg}(X) \ : \ g(\Pex)=\Pex\}$, by 
Theorem \ref{thm-Mon-2-Hdg-is-a-semi-direct-product} part 
\ref{thm-item-stabilizer-of-FE} and Lemma \ref{lemma-faces-of-FE}.
Assertion (\ref{input-item-Gamma-B-is-the-stabilizer-of-B}) follows from the latter equality by 
Lemma \ref{lemma-kernel-and-image-of-restriction-homo} part \ref{lemma-item-ker-rho-in-Mon-2-Bir}.
Assertion (\ref{input-item-Gamma-is-a-semi-direct-product})
is verified in Lemma \ref{lemma-kernel-and-image-of-restriction-homo} part 
\ref{lemma-item-overline-Mon-Hdg-a-semi-direct-product}.
Assertion (\ref{input-item-K-is-fundamental-domain}) is verified in 
Lemma \ref{lemma-movable-cone-is-a-fundamental-domain}.

The argument proceeds roughly as follows. 
Choose a rational element $x_0\in \MV_X$ which is not fixed by any element of $\Gamma$. 
Let $\C_{+}$ be the convex hull of $\overline{\C}_\NS\cap \NS_\RationalNumbers$
in $\NS_\RealNumbers$.
Set 
\[
\Pi:= \{x\in \C_+ \ : \ (x_0,x) \leq (x_0,\gamma(x)), \ \mbox{for all} \ \gamma\in\Gamma\}.
\]
Then $\Pi$ is a fundamental domain for the $\Gamma$ action on $\C_+$, known as the
{\em Dirichlet domain with center $x_0$} (compare\footnote{The bilinear pairing 
$(x_0,x)$ in the above
definition of the Dirichlet domain is replaced with the hyperbolic distance $\rho(x_0,x)$ in 
Definition 1.8 in Ch. 1 of \cite{vinberg}. However, the two definitions are equivalent, by 
the relation $cosh (\rho(x_0,x))=(x_0,x)$ (see Ch. 1 section 4.2 in \cite{AVS}).
} 
with \cite{vinberg}, Ch. 1 Proposition 1.10). 
$\Pi$ is shown to be a rational convex polyhedron (\cite{sterk}, Lemma 2.3, 
see also Theorem \ref{thm-fundamental-domain-of-a-discrete-group-of-motions} 
part (\ref{thm-item-polyhedron-fundamental-domain-of-an-arithmetic-subgroup}) above). 
The above depends only on 
Assertion (\ref{input-item-arthmetic-group}). The interior of 
any fundamental domain for $\Gamma$ can not intersect any hyperplane $e^\perp$,  $e\in \Pex$.
Hence, $\Pi$ is contained in $\MV_X^+$, by Assertions (\ref{input-item-W-is-generated-by-reflections-in-B})
and (\ref{input-item-K-is-fundamental-domain}). 
$\MV_X^+$ is a fundamental domain for the $W_{Exc}$ action on $\C_+$, by Assertion
(\ref{input-item-K-is-fundamental-domain}).
Hence, any fundamental domain for the $\Gamma$-action on $\C_+$ which is contained
in $\MV_X^+$, is a fundamental domain for the $\Gamma_{Bir}$ action on 
$\MV_X^+$, by Assertions (\ref{input-item-Gamma-B-is-the-stabilizer-of-B}) and  
(\ref{input-item-Gamma-is-a-semi-direct-product}).
\end{proof}

\begin{proof} (Of Theorem \ref{thm-finite-number-of-Bir-X-orbits})
Assume that $D$ is an irreducible divisor on $X$. Then $D$ is either prime exceptional, or the class 
%some positive integer multiple of 
$[D]$ belongs to $\overline{\MV}_X$, by Theorem \ref{thm-Zariski-decomposition}.
If $D$ is prime exceptional, the statement follows by the same argument used in the $K3$ surface 
case (\cite{sterk}, Proposition 2.5). Otherwise, $[D]$ belongs to $\MV_X^+$, and there exists 
$g\in \Gamma_{Bir}$, such that $g([D])$ belongs to the rational convex polyhedron $\Pi$ in Theorem
\ref{thm-movable-conj-conj}. 
The intersection $\Pi\cap \NS$ is a finitely generated semi-group. Choose generators $\{x_1, \dots, x_m\}$.
Then $(x_i,x_i)\geq 0$, and $(x_i,x_j)>0$, if $x_i$ and $x_j$ are linearly independent. It follows that
$\Pi\cap \NS$ contains at most finitely many elements of any given positive Beauville-Bogomolov 
degree, and at most finitely many primitive isotropic classes. 
\end{proof}

%**************
% Hide
%**************
\hide{
%*********************************************************************
%
%*********************************************************************
\subsection{Monodromy-reflective chambers and $Mon^2_{Bir}$ as a semi-direct
product}
\label{sec-monodromy-reflective-chambers}
In addition to the K\"{a}hler-type and exceptional chamber decompositions,
there is yet a third natural 
chamber decomposition of the positive cone $\C_X$, which exists by general results
on hyperbolic reflection groups.
%One chamber decomposition consists of fundamental domains for the 
%$Mon^2_{Hdg}(X)$-action 
%(Theorem \ref{thm-fundamental-domain-of-a-discrete-group-of-motions}), 
%the other 
This chamber decomposition consists of  fundamental domains for
the reflection subgroup of $Mon^2_{Hdg}$
(Theorem \ref{thm-any-discrete-group-of-motions-is-a-semi-direct-product}).
%We discuss the latter in this section.

Let $R_\ell$ be a reflection in $Mon^2_{Hdg}(X)$, 
given by $R_\ell(\lambda)=\lambda-\frac{2(\ell,\lambda)}{(\ell,\ell)}\ell$,
where $\ell\in H^2(X,\RealNumbers)$, with $(\ell,\ell)\neq 0$.
Then the kernel of $id+R_\ell$ is a rank $1$ integral Hodge substructure,
which is contained in the line ${\rm span}_\RealNumbers(\ell)$. 
We may thus assume that $\ell$ is integral, primitive, and of Hodge-type $(1,1)$,
since $R_\ell=R_{k\ell}$, for any non-zero scalar $k$.
Furthermore, $(\ell,\ell)<0$, since for $R_\ell$ to be a monodromy operator,
it must be orientation preserving.

\begin{defi}
\label{def-monodromy-reflective}
\begin{enumerate}
\item
A class $\ell\in H^{1,1}(X,\Integers)$ is called 
{\em monodromy-reflective}, if $\ell$ is a primitive class, $(\ell,\ell)<0$, and 
$R_\ell$ is a monodromy operator.
A holomorphic line bundle $L\in \Pic(X)$ is called 
{\em monodromy-reflective}, if the class  $c_1(L)$ is {\em monodromy-reflective}. 
\item
A {\em monodromy-reflective chamber} is a connected component 
of 
\[
\C_X \setminus \bigcup\left\{\ell^\perp \ : \ \ell  \ 
\mbox{is a monodromy reflective class}\right\}.
\] 
\end{enumerate}
\end{defi}

The set $\Spe$, of stably prime-exceptional classes, is contained in the set of all
monodromy-reflective classes.  Hence, the monodromy-reflective chamber
decomposition is a refinement of the exceptional chamber decomposition.
Let $W_{Hdg}$ be the subgroup of $Mon^2_{Hdg}(X)$ generated by
all the reflections in $Mon^2_{Hdg}(X)$. 
Then $W_{Hdg}$ acts simply transitively on the set of monodromy-reflective
chambers, by Theorem \ref{thm-each-chamber-is-a-generalized-convex-polyhedron}.
$W_{Exc}$ is a normal subgroup of $W_{Hdg}$, and 
$W_{Hdg}$ is the semi-direct product of 
$W_{Exc}$ and $W_{Bir}:=W_{Hdg}\cap Mon^2_{Bir}(X)$,
by Theorem \ref{thm-Mon-2-Hdg-is-a-semi-direct-product} part \ref{thm-item-semi-direct-product}.

\begin{lem}
$W_{Bir}$ is generated by reflections. 
\end{lem}

\begin{proof}
Choose a fundamental chamber $ch$ of $W_{Hdg}$, which is contained in
$\FE_X$. Then $W_{Bir}$ consists of $g\in W_{Hdg}$,
such that $g(ch)$ is contained in $\FE_X$, by 
Theorems \ref{thm-Mon-2-Hdg-is-a-semi-direct-product}. 
Furthermore, $W_{Bir}$ acts simply-transitively on the set 
of monodromy-reflective chambers contained in $\FE_X$.
%and \ref{thm-numerical-characterization-of-BK-X}.
Now $\FE_X$ is connected, and given $g\in W_{Bir}$, we can
find a sequence of adjacent monodromy-reflective 
chambers in $\FE_X$ (sharing a wall), 
connecting $ch$ to $g(ch)$.
\end{proof}

Let $\ell$ be a monodromy-reflective class and assume that 
$R_\ell$ belongs to $W_{Bir}$. Then $R_\ell(\FE_X)=\FE_X$ and 
there exists a bimeromorphic map $g:X\rightarrow X$, such that 
$g^*=R_\ell$, by Theorem \ref{thm-Mon-2-Hdg-is-a-semi-direct-product}.
Assume that $X$ satisfies the following.

\begin{assumption}
\label{assumption-reflection-is-not-induced-by-regular-automorphism}
There does not exist\footnote{A weaker version of this assumption, namely
the non-existence of a fixed-point free such automorphism $g$, 
is always true. Indeed, if $g^*=R_\ell$, and $g$ is a fixed-point-free
(necessarily symplectic) automorphism, then $g^2$ acts trivially on $H^2(X,\Integers)$.
Hence, $g^2$ is an isometry with respect to a K\"{a}hler metric.
It follows that $g$ has finite order, since it generates a discrete
subgroup of the compact isometry group. Thus,
$X/\langle g\rangle$ is a non simply connected holomorphic symplectic K\"{a}hler manifold,
with $h^{k,0}(X)=1$, for even $k$ in the range $0\leq k\leq \dim_\ComplexNumbers(X)$,
and $h^{k,0}(X)=0$, otherwise.
Such $X$ does not exist, by \cite{huybrechts-nieper}, Proposition A.1.} 
any regular automorphism $g:X\rightarrow X$, 
such that $g^*=R_\ell$, for some monodromy-reflective class $\ell$. 
\end{assumption}

It follows that $R_\ell(\K_X)\neq \K_X$, for every reflection
$R_\ell\in W_{Bir}$. 
In particular, there exists a
unique monodromy-reflective chamber, which contains the K\"{a}hler cone $\K_X$.
We call it the {\em fundamental monodromy-reflective chamber} and 
denote it by
\[
\FR_X.
\]
We then have $\K_X\subset \FR_X\subset \FE_X\subset \C_X$.
Let $Mon^2_{\FR}(X)$ be the subgroup of $Mon^2(X)$, consisting of 
monodromy-operators $g$, such that $g(\FR_X)=\FR_X$. Then $Mon^2_{Bir}(X)$
is the semi-direct product of $W_{Bir}$ and $Mon^2_{\FR}(X)$, by Theorem
\ref{thm-any-discrete-group-of-motions-is-a-semi-direct-product}.

Examples of $X$ of type $K3^{[n]}$, with non-trivial $W_{Bir}$, are provided in
(\cite{markman-prime-exceptional}, section 11). Thus, $\FR_X$ may be strictly 
smaller than $\FE_X$.
The group $Mon^2_{\FR}(X)$ contains the image 
$Mon^2_{\Aut}(X)\subset Mon^2_{Hdg}(X)$ of the
automorphism group of $X$. There are examples of $X$
of $K3^{[n]}$-type, where $Mon^2_{\FR}(X)$
is strictly larger than $Mon^2_{\Aut}(X)$, containing an element of the 
form\footnote{The operator $-R_\ell$ is given by
$-R_\ell(\lambda)=-\lambda+\frac{2(\ell,\lambda)}{(\ell,\ell)}\ell$.}
$-R_\ell$, where $(\ell,\ell)=2$,
induced by a birational involution, which is not
biregular \cite{markman-brill-noether}. 
%****************
% End Hide
%****************
}
%****************
% Hide
%****************
\hide{
It seems however reasonable to expect 
every $g\in Mon^2_{\FR}(X)$ to be of stably biregular type in the following sense.

\begin{defi}
Let $X$ be an irreducible holomorphic symplectic manifold and $f$
a bimeromorphic map from $X$ to itself.
We call $f$  {\em stably biregular},
if a generic flat\footnote{The deformation 
$(X_t,f_t)$ is said to be flat, if the closures of the graph of $f_t$ vary
in a flat family over the base of the deformation.} 
deformation $(X_t,f_t)$,
of the pair $(X,f)$, consists of a biregular automorphism $f_t$ of $X_t$. 
An element $g$ of $Mon^2_{Hdg}(X)$ is said to be 
of {\em stably biregular type}, if $g$ is induced by a
stably biregular bimeromorphic self-map.
\end{defi}
%****************
% End Hide
%****************
}

%****************************************************************
% 
%****************************************************************
\section{The monodromy and polarized monodromy groups}
\label{sec-polarized-monodromy}
In section \ref{sec-polarized-parallel-transport-operators}
we prove Proposition \ref{prop-intro-Mon-X-H-is-a-stabilizer},
stating that the polarized monodromy group $Mon^2(X,H)$ is the stabilizer of $c_1(H)$
in $Mon^2(X)$.
In section \ref{sec-deformation-types-of-polarized-marked-pairs} 
we fix a lattice $\Lambda$ and define the coarse moduli space of polarized 
$\Lambda$-marked pairs of a given deformation type. 

%****************************************************************
% 
%****************************************************************
\subsection{Polarized parallel transport operators}
\label{sec-polarized-parallel-transport-operators}
Let $\Omega_\Lambda$ be a period domain as in equation
(\ref{eq-period-domain}).
Choose a connected component $\FM^0_\Lambda$ of the moduli space of
marked pairs, a class $h\in\Lambda$ with $(h,h)>0$, and let 
$\Omega_{h^\perp}^+$ be the period domain given in equation
(\ref{eq-Omega-h-perp-plus}). Let $P_0:\FM^0_\Lambda\rightarrow \Omega_\Lambda$
be the period map.
Denote the inverse image $P_0^{-1}\left(\Omega_{h^\perp}^+\right)$ by 
%\label{eq-M-h-perp-plus}
$\FM_{h^\perp}^+.$
The discussion in section \ref{sec-orientation} provides the following modular 
description of $\FM_{h^\perp}^+$.
A marked pair $(X,\eta)$ belongs to $\FM_{h^\perp}^+$, if and only if
$(X,\eta)$ belongs to $\FM^0_\Lambda$, 
%$P(X,\eta)$ belongs to $\Omega_{h^\perp}$, so that 
the class $\eta^{-1}(h)$ is of Hodge type $(1,1)$, 
and $\eta^{-1}(h)$ belongs to the positive cone $\C_X$. 

\begin{prop}
\label{prop-path-connectedness}
$\FM_{h^\perp}^+$ is path-connected.
\end{prop}

\begin{proof}
The proof is similar to that of Proposition 5.11 in 
\cite{markman-prime-exceptional}. The proof relies on the 
Global Torelli Theorem \ref{thm-global-torelli} and the connectedness of $\Omega_{h^\perp}^+$.
\end{proof}

\begin{defi}
\label{def-polarized-monodromy-subgroup-of-O-Lambda}
Let $Mon^2\left(\FM_\Lambda^0\right)$ be the subgroup
$\eta\circ Mon^2(X)\circ\eta^{-1}\subset O(\Lambda),$ for some marked pair $(X,\eta)\in \FM^0_\Lambda$.
Let $Mon^2\left(\FM_\Lambda^0\right)_h$ be the subgroup of $Mon^2\left(\FM_\Lambda^0\right)$
stabilizing $h$.
%$Mon^2\left(\FM_{h^\perp}^+\right)\subset O(\Lambda)$ be the subgroup 
%\[
%Mon^2\left(\FM_{h^\perp}^+\right) \ \ \ := \ \ \ \eta\circ Mon^2(X,H)\circ\eta^{-1}, 
%\]
%for some marked pair $(X,\eta)\in \FM_{h^\perp}^+$. 
\end{defi}

The subgroup $Mon^2\left(\FM_\Lambda^0\right)$ is independent of 
the choice of $(X,\eta)$, since $\FM_\Lambda^0$ is connected, by definition.
%The subgroup $Mon^2\left(\FM_{h^\perp}^+\right)$
%is independent of the choice of $(X,\eta)$, 
%by Proposition \ref{prop-path-connectedness}. Furthermore, the class
%$h$ is $Mon^2\left(\FM_{h^\perp}^+\right)$-invariant.
%The subgroup $Mon^2\left(\FM_{h^\perp}^+\right)$ is the largest 
%subgroup of $O(\Lambda)$ leaving the subset $\FM_{h^\perp}^+$
%of $\FM_\Lambda$ invariant. In particular, 
$Mon^2\left(\FM_\Lambda^0\right)_h$
naturally acts on $\FM_{h^\perp}^+$.

Let 
\begin{equation}
\label{eq-connected-component-ofmoduli-of-polarized-marked-pairs}
\FM_{h^\perp}^a
\end{equation} 
be the subset of $\FM_{h^\perp}^+$,
consisting of isomorphism classes of pairs $(X,\eta)$, such that $\eta^{-1}(h)$
is an ample class of $X$. The stability of K\"{a}hler manifolds 
implies that $\FM_{h^\perp}^a$ is an open subset of $\FM_{h^\perp}^+$
(\cite{voisin-book-vol1}, Theorem 9.3.3). 
We refer to $\FM_{h^\perp}^a$ as a {\em connected component of the 
moduli space of polarized marked pairs.}

\begin{cor}
\label{cor-moduli-of-polarized-marked-pairs-is-path-connected}
$\FM_{h^\perp}^a$ is a $Mon^2\left(\FM_\Lambda^0\right)_h$-invariant
path-connected open Hausdorff subset of $\FM_{h^\perp}^+$.
The period map restricts as an injective open 
$Mon^2\left(\FM_\Lambda^0\right)_h$-equivariant morphism from $\FM_{h^\perp}^a$
onto an open dense subset of $\Omega_{h^\perp}^+$.
\end{cor}

\begin{proof}
Let us check first that 
$\FM_{h^\perp}^a$ is  $Mon^2\left(\FM_\Lambda^0\right)_h$-invariant.
Indeed, let $(X,\eta)$ belong to $\FM_{h^\perp}^a$ and let 
$g$ be an element of $Mon^2\left(\FM_\Lambda^0\right)_h$.
Denote by $H$ the line bundle with $c_1(H)=\eta^{-1}(h)$.
Then $g=\eta f \eta^{-1}$, for some $f\in Mon^2(X)$ stabilizing $c_1(H)$, 
by definition of $Mon^2\left(\FM_\Lambda^0\right)_h$. The pair $(X,g\eta)=(X,\eta f)$
belongs to $\FM_\Lambda^0$, since $f$ is a monodromy-operator. 
We have
\[
(g\eta)^{-1}(h)=
f^{-1}(\eta^{-1}(h))=f^{-1}(c_1(H))=c_1(H).
\] 
%where the third equality follows from the definition of $Mon^2(X,H)$. 
Hence, $(g\eta)^{-1}(h)$ is an ample class in $H^{1,1}(X,\Integers)$.

Let $(X,\eta)$ and $(Y,\psi)$ be two inseparable points of $\FM_{h^\perp}^a$.
Then $\psi^{-1}\eta$ is a parallel-transport operator, preserving the Hodge structure,
by Theorem \ref{thm-non-separated-implies-bimeromorphic}. Furthermore, 
$\psi^{-1}\eta$ maps the ample class $\eta^{-1}(h)$ to the
ample class $\psi^{-1}(h)$, by definition. 
Hence, there exists an isomorphism $f:X\rightarrow Y$, such that $f_*=\psi^{-1}\eta$, 
by Theorem \ref{thm-Hodge-theoretic-Torelli} part \ref{thm-item-isomorphic}.
The two pairs $(X,\eta)$ and $(Y,\psi)$ are thus isomorphic. Hence, 
$\FM_{h^\perp}^a$ is a Hausdorff subset of $\FM_{h^\perp}^+$.
%Write $x\sim y$, 
%if $x$ and $y$ are inseparable points of $\FM_{h^\perp}^+$.
%Verbitsky proves that the relation $\sim$ is an equivalence relation, and the quotient
%$\FM_{h^\perp}^+/\sim$, with its quotient topology, 
%is a Hausdorff complex manifold, called 
%the Hausdorff reduction of $\FM_{h^\perp}^+$ (\cite{verbitsky}, section 4).
%$\FM_{h^\perp}^a$ maps isomorphically onto an open subset of $\FM_{h^\perp}^+/\sim$,

$\FM_{h^\perp}^a$  is the complement of a countable union of 
closed complex analytic subsets of $\FM_{h^\perp}^+$.
Hence, $\FM_{h^\perp}^a$ is path-connected
(see, for example, \cite{verbitsky}, Lemma 4.10).

The period map restricts to an injective map on any Hausdorff subset of a
connected component of the moduli space of marked pairs, by Theorem
\ref{thm-global-torelli}. The image of $\FM_{h^\perp}^a$ contains the 
subset of $\Omega_{h^\perp}^+$, consisting of points $\period$,
such that $\Lambda^{1,1}(\period)={\rm span}_\Integers\{h\}$, by Huybrechts' projectivity criterion
\cite{huybrechts-basic-results}, and Theorem \ref{thm-global-torelli}. Hence, the image of 
$\FM_{h^\perp}^a$ is dense in $\Omega_{h^\perp}^+$. The image is open, since 
$\FM_{h^\perp}^a$ is an open subset and the period map is open, being a local
homeomorphism.
\end{proof}

Let $(X_i,H_i)$, $i=1,2$, be two pairs, each  consisting of a 
projective irreducible holomorphic symplectic manifold $X_i$,
and an ample line bundle $H_i$. Set $h_i:=c_1(H_i)$. 

\begin{cor}
\label{cor-Mon-X-H-is-a-stabilizer}
A parallel transport operator $f:H^2(X_1,\Integers)\rightarrow H^2(X_2,\Integers)$
is a polarized parallel transport operator from $(X_1,H_1)$ to $(X_2,H_2)$
(Definition \ref{def-monodromy}),
if and only if $f(h_1)=h_2$.
\end{cor}

\begin{proof}
The `only if' part is clear. We prove the `if' part.
%The proof uses the Global Torelli Theorem \ref{thm-global-torelli}.
Assume that $f(h_1)=h_2$.
Choose a marking $\eta_2:H^2(X_2,\Integers)\rightarrow \Lambda$,
and set $\eta_1:=\eta_2\circ f$. Then $\eta_1(h_1)=\eta_2(h_2)$.
Denote both $\eta_i(h_i)$ by $h$. 
Let $\FM_\Lambda^0$ be the connected component of $(X_1,\eta_1)$.
Then $(X_2,\eta_2)$ belongs to $\FM_\Lambda^0$, by the assumption that $f$ is
a parallel transport operator. Consequently, $P_0(X_i,\eta_i)$, $i=1,2$,
both belong to the same connected component of
$\Omega_{h^\perp}$.
We may choose $\eta_2$, so that this connected component is $\Omega_{h^\perp}^+$.
Then $(X_1,\eta_1)$ and $(X_2,\eta_2)$ both belong to $\FM_{h^\perp}^a$. 

Choose a path $\gamma:[0,1]\rightarrow \FM_{h^\perp}^a$,
with $\gamma(0)=(X_1,\eta_1)$ and $\gamma(1)=(X_2,\eta_2)$. 
This is possible, by Corollary \ref{cor-moduli-of-polarized-marked-pairs-is-path-connected}.
For each $t\in [0,1]$, there exists a simply-connected  open neighborhood $U_t$
of $\gamma(t)$ in $\FM_{h^\perp}^a$ and a semi-universal family
$\pi_t:\X_t\rightarrow U_t$.  The covering $\{U_t\}_{t\in [0,1]}$ of $\gamma([0,1])$
has a finite sub-covering $\{V_j\}_{j=1}^k$, for some 
%Footnote:
integer\footnote{We could take $k=1$, if there exists a universal family over 
$\FM_{h^\perp}^a$, but such a family need not exist.} 
%End Footnote
$k\geq 1$,  with the property that $\gamma\left(\left[\frac{j-1}{k},\frac{j}{k}\right]\right)$
is contained in $V_j$.
Consider the analytic space $B$, obtained from the disjoint union 
of $V_j$, $1\leq j\leq k$, by gluing $V_j$ to $V_{j+1}$
at the single point $\gamma\left(\frac{j}{k}\right)$ with transversal Zariski tangent spaces.
%Choose an isomorphism of $\X_{t_i,u_{i,i+1}}$ with $\X_{t_{i+1},u_{i,i+1}}$, 
Let $\pi_j:\X_j\rightarrow V_j$ be the universal family and denote its fiber over $v\in V_j$ by 
$\X_{j,v}$. 
Endow each fiber $\X_{j,v}$, of $\pi_j$ over
$v\in V_j$, with the marking $H^2(\X_{j,v},\Integers)\rightarrow \Lambda$
corresponding to the point $v$. 
For $1\leq j\leq k$, choose an isomorphism of
$\X_{j,\gamma\left(\frac{j}{k}\right)}$  with $\X_{j+1,\gamma\left(\frac{j}{k}\right)}$
compatible with the marking chosen, and use it to glue the family
$\pi_j$ to the family $\pi_{j+1}$. We get a family $\pi:\X\rightarrow B$. 
The paths $\gamma:\left[\frac{j-1}{k},\frac{j}{k}\right]\rightarrow V_j$ 
can now be reglued to a path $\tilde{\gamma}:[0,1]\rightarrow B$.
Parallel transport along $\tilde{\gamma}$ induces the isomorphism 
$\eta_{\tilde{\gamma}(1)}^{-1}\circ\eta_{\tilde{\gamma}(0)}=
\eta_{\gamma(1)}^{-1}\circ\eta_{\gamma(0)}=\eta_2^{-1}\circ\eta_1=f$.
Hence, $f$ is a polarized parallel transport operator from 
$(X_1,H_1)$ to $(X_2,H_2)$.
\end{proof}

%****************************************************************
% 
%****************************************************************
\subsection{Deformation types of polarized marked pairs}
\label{sec-deformation-types-of-polarized-marked-pairs}
Fix an irreducible holomorphic symplectic manifold $X_0$ and let 
$\Lambda$ be the lattice $H^2(X_0,\Integers)$, endowed with the Beauville-Bogomolov 
pairing. Let $\tau$ be the set of connected components of $\FM_\Lambda$,
consisting of pairs $(X,\eta)$, such that $X$ is deformation equivalent to $X_0$.

\begin{lem}
\label{lemma-tau-is-finite}
The set $\tau$ is finite. The group $O(\Lambda)$ acts transitively on $\tau$
and the stabilizer of a connected component $\FM_\Lambda^0\in \tau$ is  the subgroup $Mon^2(\FM_\Lambda^0)$, introduced in Definition 
\ref{def-polarized-monodromy-subgroup-of-O-Lambda}.
\end{lem}

\begin{proof}
The set $O[H^2(X,\Integers)]/Mon^2(X)$ is finite, by a result of Sullivan
\cite{sullivan}
(see also \cite{verbitsky}, Theorem 3.4). The rest of the statement is clear.
\end{proof}

Denote by $\FM_\Lambda^\tau$
the disjoint union of connected components parametrized by the set $\tau$. 
We refer to $\FM_\Lambda^\tau$ as {\em the moduli space of marked pairs of 
deformation type} $\tau$. An example would be the moduli space of
marked pairs of $K3^{[n]}$-type. 
Given a point $t\in\tau$, denote by $\FM^t_\Lambda$ the corresponding 
connected component of $\FM^\tau_\Lambda$.

\begin{rem}
If  $Mon^2(X)$ is a normal subgroup of $O[H^2(X,\Integers)]$, 
then the subgroup $Mon^2(\FM_\Lambda^t)$ of $O(\Lambda)$ 
is equal to a fixed subgroup $Mon^2(\tau,\Lambda)\subset O(\Lambda)$, 
for all $t\in\tau$. This is the case 
when $X$ is of $K3^{[n]}$-type (Theorem \ref{thm-mon-2-is-a-reflection-group}).
The set $\tau$ is an $O(\Lambda)/Mon^2(\tau,\Lambda)$-torsor,
by Lemma \ref{lemma-tau-is-finite}. We will identify the torsor $\tau$ 
with an explicit lattice theoretic $O(\Lambda)/Mon^2(\tau,\Lambda)$-torsor
in Corollary \ref{cor-enumeration-of-connected-components}.
\end{rem}

We get the {\em refined period map} 
\begin{equation}
\label{eq-refined-period-map}
\widetilde{P} \ : \ \FM_\Lambda^\tau \ \ \ \longrightarrow \ \ \ 
\Omega_\Lambda\times \tau,
\end{equation}
sending a marked pair $(X,\eta)$ to the pair $\left(P(X,\eta),t\right)$,
where $\FM^t_\Lambda$ is the connected component containing $(X,\eta)$. 
Then $\widetilde{P}$ is $O(\Lambda)$-equivariant with respect to the diagonal action
of $O(\Lambda)$ on $\Omega_\Lambda\times \tau$.

Given $h\in \Lambda$, with $(h,h)>0$, denote by 
$\Omega_{h^\perp}^{t,+}$ the period domain associated to $\FM^t_\Lambda$ in 
equation (\ref{eq-Omega-h-perp-plus}).
Set $\FM_{h^\perp}^{t,+}:=\widetilde{P}^{-1}(\Omega_{h^\perp}^{t,+})$. 
Let $\FM_{h^\perp}^{t,a}\subset \FM_{h^\perp}^{t,+}$ be the open subset of 
polarized pairs introduced in equation 
(\ref{eq-connected-component-ofmoduli-of-polarized-marked-pairs}).

%*******************
% Hide
%*******************
\hide{
Given a class $h\in \Lambda$, with $(h,h)>0$, 
we can now define the refined period space
\[
\Omega_{h^\perp}^{\tau,+} \ \ \subset \ \  \Omega_\Lambda\times \tau,
\]
analogous to the one introduced in equation (\ref{eq-Omega-h-perp-plus}), 
and its inverse image
\[
\FM_{h^\perp}^{\tau,+} \ \ := \ \ \widetilde{P}^{-1}(\Omega_{h^\perp}^{\tau,+}).
\]
Denote by $O(\Lambda,h)$ the subgroup of $O(\Lambda)$ stabilizing $h$.
Then the refined period map restricts to an $O(\Lambda,h)$-equivariant 
surjective morphism
\begin{equation}
\label{eq-refined-period-map-polarized-type-tau}
\widetilde{P} \ : \ \FM_{h^\perp}^{\tau,+} \ \ \ \longrightarrow \ \ \ \Omega_{h^\perp}^{\tau,+}.
\end{equation}
The subset $\FM_{h^\perp}^{\tau,a}\subset \FM_{h^\perp}^{\tau,+}$,
consisting of pairs $(X,\eta)$, such that $\eta^{-1}(h)$ is an ample class, 
is again an open, $O(\Lambda,h)$-invariant subset.
%*******************
% End Hide
%*******************
}

We construct next a polarized analogue of the refined period map.
Given an $O(\Lambda)$-orbit $\bar{h}\subset \Lambda\times \tau$, of 
pairs  $(h,t)$ with $(h,h)>0$, consider the disjoint unions
\begin{eqnarray*}
\FM_{\bar{h}}^+ & := &  \bigcup_{(h,t)\in \bar{h}} \FM_{h^\perp}^{t,+},
\\
\Omega_{\bar{h}}^+ & := & \bigcup_{(h,t)\in \bar{h}} \Omega_{h^\perp}^{t,+},
\end{eqnarray*}
and let 
\begin{equation}
\label{eq-universal-period-map-for-polarized-marked-pairs}
\widetilde{P} \ : \  \FM_{\bar{h}}^+ \ \ \ \longrightarrow \ \ \ \ \Omega_{\bar{h}}^+
\end{equation}
be the map induced by the refined  period map on each connected component.
Then $\widetilde{P}$ is $O(\Lambda)$-equivariant and
surjective. The disjoint union 

\begin{equation}
\label{eq-moduli-space-of-polarized-pairs-of-type-tau-bar-h}
\FM_{\bar{h}}^a:=\bigcup_{(h,t)\in \bar{h}} \FM_{h^\perp}^{t,a}
\end{equation}

\noindent
is an $O(\Lambda)$-invariant open subset of $\FM_{\bar{h}}^+$. 
This open subset will be called {\em the moduli space of polarized marked pairs 
of deformation type $\bar{h}$.} 
Indeed, $\FM_{\bar{h}}^a$ coarsely represents a functor
from the category of analytic spaces to sets, associating to
a complex analytic  space $T$ the set of all equivalence classes of families of 
marked polarized triples $(X,L,\eta)$, where $X$ is of deformation type 
$\tau$, $L$ is an ample line bundle, and 
$\eta:H^2(X,\Integers)\rightarrow\Lambda$ is an isometry, such that the pair
$\left[\eta(c_1(L)),t\right]$ belongs to the $O(\Lambda)$-orbit $\bar{h}$,
where $\FM^t_\Lambda$ is the connected component of $(X,\eta)$.
A family $(\pi:\X\rightarrow T, \LB, \tilde{\eta})$ consists of a family $\pi$, an element $\LB$ of $\Pic(\X/T)$
and a trivialization $\tilde{\eta}:R^2\pi_*\Integers\rightarrow (\Lambda)_T$, via isometries.
Two families $(\X\rightarrow T, \LB, \tilde{\eta})$ and $(\X'\rightarrow T, \LB', \tilde{\eta}')$
are equivalent, if there exists a $T$-isomorphism $f:\X\rightarrow \X'$, such that
$f^*\LB'\cong\LB$ and $\tilde{\eta}'=\tilde{\eta}\circ f^*$.
We omit the detailed definition of this functor, as well as the proof that 
$\FM_{\bar{h}}^a$ coarsely represents it, as we will not use 
the latter fact below.

%****************************************************************
% 
%****************************************************************
\section{Monodromy quotients of type IV period domains}
\label{sec-monodromy-quotients-of-period-domains}
Fix a connected component $\FM_{h^\perp}^a$ of the moduli space 
$\FM_{\bar{h}}^a$ of polarized marked pairs
of polarized deformation type  $\bar{h}$. In the notation of section 
\ref{sec-deformation-types-of-polarized-marked-pairs}, 
$\FM_{h^\perp}^a:=\FM_{h^\perp}^{t,a}$, for some $(h,t)\in\bar{h}$. 
Let $\FM_\Lambda^0$ be the connected component of $\FM_\Lambda$ containing $\FM_{h^\perp}^a$.
Set $\Gamma:=Mon^2\left(\FM_\Lambda^0\right)_h$
(Definition \ref{def-polarized-monodromy-subgroup-of-O-Lambda}).
The period domain $\Omega_{h^\perp}^+$ is a homogeneous domain of type IV
(\cite{satake}, Appendix, section 6). $\Gamma$ is an arithmetic group, by
(\cite{verbitsky}, Theorem 3.5).
The quotient $\Omega_{h^\perp}^+/\Gamma$  
is thus a normal quasi-projective variety \cite{baily-borel}. 
%The following statement is tautological.

\begin{lem}
\label{lemma-open-immersion}
There exist natural isomorphisms of complex analytic spaces
\begin{eqnarray*}
\FM_{\bar{h}}^+/O(\Lambda) & \longrightarrow & \FM_{h^\perp}^+/\Gamma,
\\
\FM_{\bar{h}}^a/O(\Lambda) & \longrightarrow & \FM_{h^\perp}^a/\Gamma,
\\
\Omega^+_{\bar{h}}/O(\Lambda) & \longrightarrow & \Omega^+_{h^\perp}/\Gamma.
\end{eqnarray*}
Furthermore, the period map descends to an open embedding
\begin{equation}
\label{eq-open-embedding}
\overline{P} \ : \ 
\FM_{\bar{h}}^a/O(\Lambda) \ \ \hookrightarrow \ \ \Omega_{h^\perp}^+/\Gamma.
\end{equation}
\end{lem}

\begin{proof}
We have the following commutative equivariant diagram
\[
\begin{array}{ccccc}
\FM_{\bar{h}}^a & \LongRightArrowOf{\widetilde{P}} & \Omega_{\bar{h}}^+ &
\longrightarrow & \Omega_{\bar{h}}^+/O(\Lambda)
\\
\uparrow & & \uparrow & & \uparrow
\\
\FM^a_{h^\perp} & \LongRightArrowOf{P_0} & \Omega_{h^\perp}^+ & \longrightarrow
& \Omega_{h^\perp}^+/\Gamma,
\end{array}
\]
with respect to the $O(\Lambda)$ action on the top row, the
$\Gamma$-action on the bottom,
and the inclusion homomorphism $\Gamma\hookrightarrow O(\Lambda)$.
$O(\Lambda)$ acts transitively on its orbit $\bar{h}$, and the stabilizer of the pair
$(h,\ \FM^+_{h^\perp})\in \bar{h}$ is precisely $\Gamma$,
by Lemma \ref{lemma-tau-is-finite} and 
Proposition \ref{prop-intro-Mon-X-H-is-a-stabilizer}.

The morphism (\ref{eq-open-embedding}) is an open embedding,
since the $\Gamma$-equivariant open morphism
$\FM_{h^\perp}^a\rightarrow \Omega^+_{h^\perp}$ is  injective, 
by Corollary
\ref{cor-moduli-of-polarized-marked-pairs-is-path-connected}.
\end{proof}

A {\em polarized irreducible holomorphic symplectic manifold}
is a pair $(X,L)$, consisting of a smooth projective irreducible holomorphic symplectic variety
$X$ and an ample line bundle $L$. 
Consider the contravariant functor $F'$ from the category of schemes over 
$\ComplexNumbers$ to the category of sets,
which associates to a scheme $T$ the set of isomorphism classes of flat families of  
polarized irreducible holomorphic symplectic manifolds $(X,L)$ over $T$, 
with a fixed Hilbert polynomial $p(n):=\chi(L^n)$. 
The coarse moduli space representing the functor $F'$ was constructed by 
Viehweg as a quasi-projective scheme with quotient singularities \cite{viehweg}.
Fix a connected component $\V$ of this moduli space. 
Then $\V$ is a quasi-projective variety. Denote by $F$ the functor 
represented by the connected component $\V$. The universal property of
a coarse moduli space asserts that 
there is a natural transformation 
$\theta:F\rightarrow \Hom(\bullet,\V)$, satisfying the following properties.
\begin{enumerate}
\item
\label{item-universal-property-bijection}
$\theta({\rm Spec}(\ComplexNumbers)):
F({\rm Spec}(\ComplexNumbers))\rightarrow \V$ is bijective.
\item
\label{item-unique-morphism-psi}
Given a scheme $B$ and a natural transformation $\chi:F\rightarrow \Hom(\bullet,B)$,
there is a unique morphism $\psi:\V\rightarrow B$, hence a
natural transformation $\psi_*:Hom(\bullet,\V)\rightarrow \Hom(\bullet,B)$,
with $\chi=(\psi_*)\circ\theta$.
\end{enumerate}

\begin{rem}
\label{remark-unique-morphism-psi-versus-universal-family}
Property (\ref{item-unique-morphism-psi}) replaces the data of a universal family
over $\V$, which may not exist when $\V$ fails to be a fine moduli space.
When a universal family $\U\in F(\V)$ exists, then the morphism $\psi$ is the 
image of $\U$ via $\chi:F(\V)\rightarrow \Hom(\V,B)$.
\end{rem}

Denote by $\bar{h}$ the deformation type of a polarized pair $(X,L)$ in $\V$.
We regard $\bar{h}$ both as a point in $[\Lambda\times \tau]/O(\Lambda)$
and as a subset of $\Lambda\times \tau$.
Choose a point $(h,t)\in \bar{h}$ and set $\Omega^+_{h^\perp}:=\Omega^{t,+}_{h^\perp}$.

%********************
% Hide
%********************
\hide{
Let $G$ be the contravariant functor from the category of analytic spaces to 
the category of sets, analogous to the functor $F$ above.
Then the quotient $\FM^a_{\bar{h}}/O(\Lambda)$ is an analytic space, which coarsely
represents the functor $G$ in the  sense that 
there is a natural transformation 
\begin{equation}
\label{eq-natural-transformation-gamma}
\gamma:G\rightarrow Hom_{an}\left(\bullet,\FM_{\bar{h}}^a/O(\Lambda)\right)
\end{equation}
with the following properties.
\begin{enumerate}
\item
$\gamma({\rm Spec(\ComplexNumbers)}):
G({\rm Spec}(\ComplexNumbers))\rightarrow \FM_{\bar{h}}^a/O(\Lambda)$
is bijective.
\item
Given an analytic space $B$ and a natural transformation 
$\chi:F\rightarrow \Hom_{an}(\bullet,B)$,
there is a unique morphism $\phi:\FM_{\bar{h}}^a/O(\Lambda)\rightarrow B$, hence a
natural transformation 
$\phi_*:Hom_{an}(\bullet,\FM_{\bar{h}}^a/O(\Lambda))\rightarrow \Hom_{an}(\bullet,B)$,
with $\chi=(\phi_*)\circ\gamma$.
\end{enumerate}
%********************
% End Hide
%********************
}

\begin{lem}
There exists a natural injective and surjective  
morphism $\varphi:\V\rightarrow \FM_{\bar{h}}^a/O(\Lambda)$
in the category of analytic spaces. 
\end{lem}

\begin{proof}
The morphism 
$\Phi:\V\rightarrow \Omega_{\bar{h}}^+/O(\Lambda)\cong \Omega^+_{h^\perp}/\Gamma$,
sending an isomorphism class of a polarized pair $(X,L)$
to its period, is constructed in the proof of (\cite{GHS}, Theorem 1.5). 
The morphism $\Phi$ is set-theoretically injective, by the
Hodge theoretic Torelli Theorem \ref{thm-Hodge-theoretic-Torelli}.
The image $\Phi(\V)$ is the same subset as the image 
$P\left(\FM^a_{\bar{h}}\right)$, by definition of the two moduli spaces.
The latter is the image also of the
open immersion $\overline{P}:\FM^a_{\bar{h}}/O(\Lambda)\hookrightarrow 
\Omega_{\bar{h}}^+/O(\Lambda)$, by Lemma \ref{lemma-open-immersion}.
Hence, the composition 
$\overline{P}^{-1}\circ \Phi:\V\rightarrow \FM^a_{\bar{h}}/O(\Lambda)$
is well defined and we denote it by $\varphi$. 
%********************
% Hide
%********************
\hide{
The existence of the 
morphism $\varphi:\V\rightarrow \FM_{\bar{h}}^a/O(\Lambda)$
follows\footnote{For an explicit proof, see the direct construction 
of the composite morphism 
$\overline{P}\circ \varphi:\V\rightarrow \Omega_{\bar{h}}^+/O(\Lambda)$
in the proof of (\cite{GHS}, Theorem 1.5). By construction, their morphism factors through
$\FM_{\bar{h}}^a/O(\Lambda)$.}
from the universal properties of the coarse moduli spaces
$\V$ and $\FM_{\bar{h}}^a/O(\Lambda)$ as follows.

Let $H$ be the natural functor from the category of schemes over $\ComplexNumbers$
to the category of analytic spaces. Let $\chi_1:F\rightarrow GH$
be the natural transformation, corresponding to the 
fact that any complex algebraic family is in particular complex analytic. 
Let 
$\gamma:G\rightarrow Hom_{an}\left(\bullet,\FM_{\bar{h}}^a/O(\Lambda)\right)$
be the natural transformation (\ref{eq-natural-transformation-gamma}).
We get the natural transformation
\[
\gamma(H):GH\rightarrow 
Hom_{an}\left(\bullet,\FM_{\bar{h}}^a/O(\Lambda)\right).
\]
Let $\chi:=\gamma(H)\circ \chi_1:
F\rightarrow Hom_{an}\left(\bullet,\FM_{\bar{h}}^a/O(\Lambda)\right)$
be the composite natural transformation. 
The universal property (\ref{item-unique-morphism-psi}) of $\V$ does not suffice to
associate to the natural transformation $\chi$ the desired morphism 
$\psi$ from $\V^{an}:=H(\V)$ to $\FM_{\bar{h}}^a/O(\Lambda)$. The difficulty arises 
since the target functor of $\chi$ involves the category of analytic spaces.
Nevertheless, the universal property (\ref{item-unique-morphism-psi}) of $\V$ 
can be extended to yield such an analytic morphism $\psi$
(??? reference ???), as is clear from 
Remark \ref{remark-unique-morphism-psi-versus-universal-family}
and the existence of a universal family over the Hilbert scheme, of which $\V$ is a quotient.

The morphism $\psi$ is bijective, since it induces the set-theoretic map 
\[
F({\rm Spec}(\ComplexNumbers))\rightarrow G({\rm Spec}(\ComplexNumbers)),
\]
sending an isomorphism class of a polarized pair $(X,L)$, as an 
algebraic family over a point, to its isomorphism class as an analytic family over a point.
The latter map is clearly a bijection. Hence so is $\psi$, by the universal property
(\ref{item-universal-property-bijection}) of the two coarse moduli spaces.
%********************
% End Hide
%********************
}
\end{proof}

\begin{thm}
\label{thm-viehweg-moduli-is-open-in-monodromy-quotient-of-period-space}
The composition $\Phi$ of 
\[
\V \ \ \LongRightArrowOf{\varphi} \ \
\FM^a_{\bar{h}}/O(\Lambda) \ \ \cong \ \ 
\FM^a_{h^\perp}/\Gamma 
\ \ \LongRightArrowOf{\overline{P}} \ \ 
\Omega_{h^\perp}^+/\Gamma
\]
is an open immersion
%isomorphism, onto a Zariski open subset, 
in the category of algebraic varieties.
\end{thm}

\begin{proof}
The proof is similar to that of 
Theorem 1.5 in \cite{GHS} and 
Claim 5.4 in \cite{ogrady-msri-lecture-notes}. If $\Gamma$ happens to be torsion free, then
any complex analytic morphism, from a complex algebraic variety
to $\Omega^+_{h^\perp}/\Gamma$, is an algebraic morphism, 
as a consequence of Borel's extension theorem 
\cite{borel-extension-theorem}. 
%Borel's theorem requires that 
$\Gamma$ need not be torsion free, but for sufficiently large positive integer $N$, 
the subgroup $\Gamma(N)\subset \Gamma$, acting trivially on $\Lambda/N\Lambda$, 
is torsion free, as a consequence of  (\cite{satake}, IV, Lemma 7.2). 
In our situation, where the domain $\V$ of $\Phi$ is a moduli space, 
one can apply Borel's extension theorem after passage 
to a finite cover $\widetilde{\V}\rightarrow \V$, where  $\widetilde{\V}$ is a connected component of the 
moduli space of polarized irreducible holomorphic symplectic manifolds with a  
level-$N$  structure,
as done in the proofs of (\cite{hassett-special-cubics}, Proposition 2.2.2) 
and (\cite{GHS}, Theorem 1.5).
The morphism $\Phi$ lifts to a morphism $\widetilde{\Phi}:\widetilde{\V}\rightarrow \Omega^+_{h^\perp}/\Gamma(N)$.
$\widetilde{\Phi}$ is algebraic, by Borel's extension theorem,
and a descent argument implies that so is $\Phi$. 
%We conclude that 
%$\Phi:=\overline{P}\circ \varphi:\V\rightarrow \Omega^+_{h^\perp}/\Gamma$
%is an algebraic morphism.

The morphism $P:\FM^a_{\bar{h}}\rightarrow \Omega^+_{\bar{h}}/O(\Lambda)$
is open. Hence, the image $\overline{P}(\FM^a_{h^\perp}/\Gamma)$ of  
$\Phi$ is an open subset of $\Omega^+_{h^\perp}/\Gamma$
in the analytic topology. The image of $\Phi$ 
is also a constructibe set, in the Zariski topology. The image is 
thus a Zariski dense open subset. $\Phi$ is thus 
an algebraic open immersion, by Zariski's Main Theorem.
\end{proof}

\begin{rem}
Theorem \ref{thm-viehweg-moduli-is-open-in-monodromy-quotient-of-period-space}
answers Question 2.6 in the paper \cite{GHS}, concerning the polarized 
$K3^{[n]}$-type moduli spaces. The map $\Phi$ in 
Theorem \ref{thm-viehweg-moduli-is-open-in-monodromy-quotient-of-period-space}
is denoted by $\tilde{\varphi}$ in (\cite{GHS}, Question 2.6) and
is defined in (\cite{GHS}, Theorem 2.3). There is a typo in the definition of
$\tilde{\varphi}$ in \cite{GHS}; its target 
$\widetilde{O}^+(L_{2n-2},h)\backslash{\mathcal D}_h$ should be replaced by 
$\widehat{O}^+(L_{2n-2},h)\backslash{\mathcal D}_h$. 
When $n=2$, these two quotients are the same, but for $n\geq 3$, the former is 
a branched double cover of the latter. Modulo this minor change, Theorem
\ref{thm-viehweg-moduli-is-open-in-monodromy-quotient-of-period-space}
provides an affirmative answer to (\cite{GHS}, Question 2.6).
\end{rem}

%****************************************************************
% 
%****************************************************************
\section{The $K3^{[n]}$ deformation type}
\label{sec-K3-type}
In section \ref{sec-computation-of-Mon-2}
we review results about parallel-transport operators of $K3^{[n]}$-type.
In section \ref{sec-numerical-determination-of-BK-X}
we explicitly calculate the fundamental exceptional chamber $\FE_X$ 
of a projective manifold $X$ of $K3^{[n]}$-type.

%****************************************************************
% 
%****************************************************************
\subsection{Characterization of parallel-transport operators  of $K3^{[n]}$-type}
\label{sec-computation-of-Mon-2}

In sections 
\ref{sec-first-two-characterizations-of-Mon-2-K3-n}, 
\ref{sec-third-characterization-of-Mon-2}, and \ref{sec-integral-generators}, 
we provide three useful  characterizations
of the monodromy group
$Mon^2(X)$ of an irreducible holomorphic symplectic manifold
of $K3^{[n]}$-type.
Given $X_1$ and $X_2$ of $K3^{[n]}$-type, we state 
in section \ref{sec-parallel-transport-operators-of-K3-n-type} 
a necessary and sufficient condition for an isometry
$g:H^2(X_1,\Integers)\rightarrow H^2(X_2,\Integers)$ to be a parallel-transport operator.

%****************************************************************
% 
%****************************************************************
\subsubsection{First two characterizations of $Mon^2(K3^{[n]})$}
\label{sec-first-two-characterizations-of-Mon-2-K3-n}
Let $X$ be an irreducible holomorphic symplectic manifold of $K3^{[n]}$-type.
If $n=1$, then $X$ is a $K3$ surface. In that case it is well known that
$Mon^2(X)=O^+[H^2(X,\Integers)]$ (see \cite{borcea}).
%An easy derivation of this equality from the Torelli theorem can be found,
%for example,  in \cite{markman-monodromy-I}, Corollary 6.7.
From now on we assume that $n\geq 2$.

Given a class $u\in H^2(X,\Integers)$, with $(u,u)\neq 0$, let
$R_u:H^2(X,\RationalNumbers)\rightarrow H^2(X,\RationalNumbers)$ be the reflection
$R_u(\lambda)=\lambda-\frac{2(u,\lambda)}{(u,u)}u$.
Set 
$\rho_u:=\left\{\begin{array}{ccc}
R_u & \mbox{if} & (u,u)<0,
\\
-R_u & \mbox{if} & (u,u)>0.
\end{array}\right.
$
Then $\rho_u$ belongs to $O^+[H^2(X,\RationalNumbers)]$.
Note that $\rho_u$ is an integral isometry, if $(u,u)=2$ or $-2$.
Let $\N\nolinebreak \subset \nolinebreak O^+[H^2(X,\Integers)]$ be the subgroup generated by such $\rho_u$.
\begin{equation}
\label{eq-N}
\N \ \ := \ \ 
\langle\rho_u \ : \ u\in H^2(X,\Integers) \ \mbox{and} \ (u,u)=2 \ \mbox{or} \ 
(u,u)=-2\rangle.
\end{equation}
Clearly, $\N$ is a normal subgroup.

\begin{thm}
\label{thm-mon-2-is-a-reflection-group}
(\cite{markman-constraints}, Theorem 1.2)
$Mon^2(X)=\N.$
\end{thm}

A second useful description of $Mon^2(X)$ depends on the fact that the lattice 
$H^2(X,\Integers)$ is isometric to the orthogonal direct sum
\[
\Lambda \ \ := \ \ E_8(-1)\oplus E_8(-1) \oplus U\oplus U \oplus U \oplus \Integers\delta,
\]
where $E_8(-1)$ is the negative definite (unimodular) $E_8$ root lattice, 
$U$ is the rank $2$ unimodular lattice of signature $(1,1)$, and $(\delta,\delta)=2-2n$.
See \cite{beauville} for a proof of this fact. 

Set $\Lambda^*:=\Hom(\Lambda,\Integers)$. Then $\Lambda^*/\Lambda$
is a cyclic group of order $2n-2$. Let $O(\Lambda^*/\Lambda)$ be
the subgroup of $\Aut(\Lambda^*/\Lambda)$ consisting of 
multiplication by all elements of $t\in \Integers/(2n-2)\Integers$,
such that $t^2=1$. Then $O(\Lambda^*/\Lambda)$ is isomorphic to 
$(\Integers/2\Integers)^{r}$, where $r$ is the number of distinct primes 
in the prime factorization $n-1=p_1^{d_1}\cdots p_r^{d_r}$ 
of $n-1$ (see \cite{oguiso}).
The isometry group $O(\Lambda)$ 
acts on $\Lambda^*/\Lambda$ and the
image of $O^+(\Lambda)$ in $\Aut(\Lambda^*/\Lambda)$
is equal to $O(\Lambda^*/\Lambda)$ (\cite{nikulin}, Theorem 1.14.2). 

Let $\pi:O^+(\Lambda)\rightarrow O(\Lambda^*/\Lambda)$ be the natural homomorphism.
The following characterization of the monodromy group follows from
Theorem \ref{thm-mon-2-is-a-reflection-group} via 
lattice theoretic arguments.
\begin{lem}
\label{lemma-mon-2-maps-to-pm-1}
(\cite{markman-constraints}, Lemma 4.2)
$Mon^2(X)$ is equal to the inverse image via $\pi$
of the subgroup $\{1,-1\}\subset O(\Lambda^*/\Lambda)$.
\end{lem}

We conclude that the index of $Mon^2(X)$ in $O^+[H^2(X,\Integers)]$
is $2^{r-1}$, and $Mon^2(X)=O^+[H^2(X,\Integers)]$, if and only if
$n=2$ or $n-1$ is a prime power. 
If $n=7$, for example, then $Mon^2(X)$ has index two in $O^+[H^2(X,\Integers)]$.

%****************************************************************
% 
%****************************************************************
\subsubsection{A third characterization of $Mon^2(K3^{[n]})$}
\label{sec-third-characterization-of-Mon-2}
The third characterization of $Mon^2(X)$ is more subtle, as it 
depends also on $H^4(X,\Integers)$. It is however this third 
characterization that will generalize to the case of parallel transport operators.

Given a $K3$ surface $S$, denote by $K(S)$ the integral $K$-ring
generated by the classes of complex topological vector bundles over $S$.
Let $\chi:K(S)\rightarrow \Integers$ be the Euler characteristic
$\chi(x)=\int_Sch(x)td_S$. 
Given classes $x,y\in K(S)$, let $x^\vee$ be the dual class and set
\begin{equation}
\label{eq-Mukai-pairing}
(x,y) \ \ := \ \ -\chi(x^\vee\otimes y).
\end{equation}
The above yields a unimodular symmetric bilinear pairing on $K(S)$,
called the {\em Mukai pairing} \cite{mukai-hodge}.
The lattice $K(S)$, endowed with the Mukai pairing, is isometric to the
orthogonal direct sum
\[
\widetilde{\Lambda} \ \ := \ \ E_8(-1)\oplus E_8(-1) \oplus U\oplus U \oplus U \oplus U
\]
and is called the {\em Mukai lattice}.

Let $Q^4(X,\Integers)$ be the quotient of $H^4(X,\Integers)$ by the image of the
cup product homomorphism
$
\cup :H^2(X,\Integers)\otimes H^2(X,\Integers) \rightarrow H^4(X,\Integers).
$
Clearly, $Q^4(X,\Integers)$ is a $Mon(X)$-module, and it comes with a pure
integral Hodge structure of weight $4$.
Let $q:H^4(X,\Integers)\rightarrow Q^4(X,\Integers)$ be the natural homomorphism and set
$\bar{c}_2(X):=q(c_2(TX))$.

\begin{thm}
\label{thm-Q-4}
(\cite{markman-constraints}, Theorem 1.10)
Let $X$ be of $K3^{[n]}$-type, $n\geq 4$.
\begin{enumerate}
\item
\label{thm-item-Q-4-is-torsion-free}
$Q^4(X,\Integers)$ is a free abelian group of rank $24$.
\item
The element $\frac{1}{2}\bar{c}_2(X)$ is an integral and primitive class in 
$Q^4(X,\Integers)$.
\item
\label{thm-item-Mukai-pairing-on-Q-4}
There exists a unique symmetric, even, integral, unimodular, $Mon(X)$-invariant
bilinear pairing $(\bullet,\bullet)$ on 
$Q^4(X,\Integers)$, such that 
$
\left(\frac{\bar{c}_2(X)}{2},\frac{\bar{c}_2(X)}{2}\right)=2n-2.
$
The resulting lattice $\left[Q^4(X,\Integers), (\bullet,\bullet)\right]$
is isometric to the Mukai lattice $\widetilde{\Lambda}$.
\item
\label{thm-item-the-Mon-equivariant-isometric-embedding-e}
The $Mon(X)$-module $\Hom\left[H^2(X,\Integers),Q^4(X,\Integers)\right]$
contains a unique integral rank $1$ saturated $Mon(X)$-submodule
\[
\EE(X),
\]
which is a sub-Hodge structure of type $(1,1)$.
A generator $e\in \EE(X)$ induces a Hodge-isometry
\[
e \ : \ H^2(X,\Integers) \ \ \longrightarrow \ \ \bar{c}_2(X)^\perp
\]
onto the co-rank $1$ sublattice of $Q^4(X,\Integers)$ orthogonal to $\bar{c}_2(X)$.
\end{enumerate}
\end{thm}

Parts (\ref{thm-item-Q-4-is-torsion-free}), (\ref{thm-item-Mukai-pairing-on-Q-4}), and
(\ref{thm-item-the-Mon-equivariant-isometric-embedding-e}) 
of the Theorem are explained in the following section
\ref{sec-integral-generators}.

Denote by $O(\Lambda,\widetilde{\Lambda})$ the set of primitive isometric embeddings of the
$K3^{[n]}$-lattice $\Lambda$ into the Mukai lattice $\widetilde{\Lambda}$.
The isometry groups $O(\Lambda)$ and $O(\widetilde{\Lambda})$
act on $O(\Lambda,\widetilde{\Lambda})$. The action
on $\iota\in O(\Lambda,\widetilde{\Lambda})$, of elements
$g\in O(\Lambda)$, and $f\in O(\widetilde{\Lambda})$, is given by
$
(g,f)\iota=f\circ\iota\circ g^{-1}.
$
\begin{lem}
\label{lemma-orbits-of-isometric-embeddings}
(\cite{markman-constraints}, Lemma 4.3)
$O^+(\Lambda)\times O(\widetilde{\Lambda})$ acts transitively on 
$O(\Lambda,\widetilde{\Lambda})$.
The subgroup $\N\subset O^+(\Lambda)$, given in
(\ref{eq-N}),
is equal to 
the stabilizer in $O^+(\Lambda)$
of every point in the orbit space
$O(\Lambda,\widetilde{\Lambda})/O(\widetilde{\Lambda})$.
\end{lem}

The lemma implies that $O(\Lambda,\widetilde{\Lambda})$ is a finite set of order
$[\N:O^+(\Lambda)]$.
The following is our third characterization of $Mon^2(X)$.

\begin{cor}
\label{cor-third-characterization-of-Mon-2}
\begin{enumerate}
\item
\label{cor-item-natural-orbit-iota-X}
An irreducible holomorphic symplectic manifold
$X$ of $K3^{[n]}$-type, $n\geq 2$, comes with a 
natural choice of an $O(\widetilde{\Lambda})$-orbit
$\iota_X$ of primitive isometric embeddings of 
$H^2(X,\Integers)$ in the Mukai lattice $\widetilde{\Lambda}$.
\item
\label{cor-item-Mon-2-is-stabilizer-of-iota-X}
The subgroup $Mon^2(X)$ of $O^+[H^2(X,\Integers)]$ is equal to 
the stabilizer of $\iota_X$ as an element of the orbit space
$O\left(H^2(X,\Integers),\widetilde{\Lambda}\right)/O(\widetilde{\Lambda})$.
\end{enumerate}
\end{cor}

\begin{proof}
Part (\ref{cor-item-natural-orbit-iota-X}):
If $n=2$, or $n=3$, then 
$O(\Lambda,\widetilde{\Lambda})$ is a singleton, and there is nothing to prove.
Assume that $n\geq 4$. Let $e:H^2(X,\Integers)\rightarrow Q^4(X,\Integers)$ 
be one of the two generators of $\EE(X)$.
Choose an isometry $g:Q^4(X,\Integers)\rightarrow\widetilde{\Lambda}$.
This is possible by Theorem \ref{thm-Q-4}.
Set $\iota:=g\circ e:H^2(X,\Integers)\rightarrow \widetilde{\Lambda}$
and let $\iota_X$ be the orbit $O(\widetilde{\Lambda})\iota$.
Then $\iota_X$ is independent of the choice of $g$. 
If we choose $-e$ instead we get the same orbit, since $-1$ belongs to
$O(\widetilde{\Lambda})$.

Part (\ref{cor-item-Mon-2-is-stabilizer-of-iota-X}):
Follows immediately from Theorem \ref{thm-mon-2-is-a-reflection-group} and 
Lemma \ref{lemma-orbits-of-isometric-embeddings}.
\end{proof}

\begin{example}
\label{example-iota-is-inverse-of-Mukai-isometry}
Let $S$ be a projective $K3$ surface, $H$ an ample line bundle on $S$, and
$v\in K(S)$ a class in the $K$-group. Denote by $M_H(v)$
the moduli space of Gieseker-Maruyama-Simpson $H$-stable coherent 
sheaves on $S$ of class $v$.
A good reference about these moduli spaces is the
book \cite{huybrechts-lehn-book}.
Assume that $M_H(v)$ is smooth and projective (i.e., we assume that 
every $H$-semi-stable sheaf is automatically also $H$-stable). 
Then $M_H(v)$ is known to be connected and of $K3^{[n]}$-type,
by a theorem due to Mukai, Huybrechts, O'Grady, and Yoshioka.
It can be found in its final form in \cite{yoshioka-abelian-surface}.

Let $\pi_i$ be the projection from $S\times M_H(v)$ onto the $i$-th factor, $i=1,2$.
Denote by $\pi_{2_!}:K[S\times M_H(v)]\rightarrow K[M_H(v)]$ the Gysin map
and by $\pi_1^!:K(S)\rightarrow K[S\times M_H(v)]$ the pull-back homomorphism.
Assume, further, that there exists a universal sheaf $\E$ over $S\times M_H(v)$.
Let $[\E]\in K[S\times M_H(v)]$ be the class of the universal sheaf in
the $K$-group. We get the natural homomorphism
\begin{equation}
\label{eq-u}
u \ : \ K(S) \ \ \rightarrow \ \ K(M_H(v)),
\end{equation}
given by $u(x):=\pi_{2_!}\left\{\pi_1^!(x^\vee)\otimes [\E]\right\}$.
Let $v^\perp\subset K(S)$ be the co-rank $1$ sub-lattice of $K(S)$
orthogonal to the class $v$ and consider Mukai's homomorphism
\begin{equation}
\label{eq-Mukai-isometry}
\theta \ : \ v^\perp \ \ \longrightarrow \ \ H^2(M_H(v),\Integers),
\end{equation}
given by $\theta(x)=c_1\left[u(x)\right]$.
Then $\theta$ is an isometry, with respect to the
Mukai pairing on $v^\perp$, and the Beauville-Bogomolov pairing on 
$H^2(\M_H(v),\Integers)$, by the work of 
Mukai, Huybrechts, O'Grady, and Yoshioka \cite{yoshioka-abelian-surface}.
Furthermore, the orbit $\iota_{M_H(v)}$ of Corollary
\ref{cor-third-characterization-of-Mon-2} is represented by the inverse of $\theta$
\begin{equation}
\label{eq-iota-is-theta-inverse}
\iota_{M_H(v)} = O[K(S)]\cdot \theta^{-1},
\end{equation}
by (\cite{markman-constraints}, Theorem 1.14).
\end{example}

%*******************************************************************************
%
%*******************************************************************************
\subsubsection{Generators for the cohomology ring $H^*(X,\Integers)$}
\label{sec-integral-generators}
Part (\ref{thm-item-Q-4-is-torsion-free}) of Theorem \ref{thm-Q-4}
is a simple consequence of the following result. Consider the case, where $X$ is a moduli
space $M$ of $H$-stable sheaves on a $K3$ surface $S$ 
and $M$ is of $K3^{[n]}$-type, as in Example \ref{example-iota-is-inverse-of-Mukai-isometry}.  
%Keep the notation of Example \ref{example-iota-is-inverse-of-Mukai-isometry}. 
Choose a basis $\{x_1, x_2, \dots, x_{24}\}$ of $K(S)$.
Let $u:K(S)\rightarrow K(M)$ be the homomorphism given in equation (\ref{eq-u}).

\begin{thm}
\label{thm-integral-generators}
(\cite{markman-integral-generators}, Theorem 1)
The cohomology ring $H^*(M,\Integers)$ is generated by the Chern classes $c_j(u(x_i))$, 
for $1\leq i\leq 24$,
and for $j$  an integer in the range $0\leq j\leq 2n$.
\end{thm}

The map $\tilde{\varphi}:K(S)\rightarrow H^4(M,\Integers)$, given by 
$\tilde{\varphi}(x)=c_2(u(x))$, is not a group homomorphism. Nevertheless, 
the composition $\varphi:=q\circ \tilde{\varphi}:K(S)\rightarrow Q^4(M,\Integers)$,
%$K(S)\LongRightArrowOf{\tilde{\varphi}} H^4(M,\Integers) \rLongRghtArrowOf{q} Q^4(M,\Integers),$
of $\tilde{\varphi}$ with the projection 
$q:H^4(M,\Integers)\rightarrow Q^4(M,\Integers)$, is a homomorphism of abelian groups 
%$\varphi:K(S)\rightarrow Q^4(M,\Integers)$
(\cite{markman-integral-generators}, Proposition 2.6). 
We note here only that $2\varphi$ is clearly a group homomorphism,
since $2c_2(y)=c_1^2(y)-2ch_2(y)$, 
the map $2ch_2:K(M)\rightarrow H^4(M,\Integers)$  is known to be 
a group homomorphism, 
and the term $c_1^2(y)$ is annihilated by the projection to $Q^4(M,\Integers)$.

Part (\ref{thm-item-Q-4-is-torsion-free}) of Theorem \ref{thm-Q-4} follows from the fact that $\varphi$
is an isomorphism. 
The homomorphism $\varphi$ is surjective, by Theorem \ref{thm-integral-generators}.
It remains to prove that $\varphi$ is injective. Injectivity would follow, once we show
that $Q^4(M,\Integers)$ has rank $24$. Now cup product 
induces an injective homomorphism $\Sym^2H^2(M,\RationalNumbers)\rightarrow H^4(M,\RationalNumbers)$,
for any irreducible holomorphic symplectic manifold of dimension $\geq 4$, by a general result
of Verbitsky \cite{verbitsky-cohomology}.
When $n\geq 4$, i.e., $\dim_\ComplexNumbers(M)\geq 8$, then 
$\dim H^4(M,\RationalNumbers)-\dim \Sym^2H^2(M,\RationalNumbers)=24$,
by G\"{o}ttsche's formula for the Betti numbers of $S^{[n]}$ \cite{gottsche}.
Hence, the rank of $Q^4(M,\Integers)$ is $24$. 

The bilinear pairing on $Q^4(M,\Integers)$, constructed in part (\ref{thm-item-Mukai-pairing-on-Q-4}) of
Theorem \ref{thm-Q-4}, is simply the push-forward via the isomorphism $\varphi$ of the Mukai pairing on 
$K(S)$. We then show that this bilinear pairing is monodromy invariant, hence it defines 
a bilinear pairing on $Q^4(X,\Integers)$, for any $X$ of $K3^{[n]}$-type.

The isometric embedding $e:H^2(M,\Integers)\rightarrow Q^4(M,\Integers)$, constructed in 
part (\ref{thm-item-the-Mon-equivariant-isometric-embedding-e}) of
Theorem \ref{thm-Q-4}, is simply the composition 
$\varphi\circ\theta^{-1}$, where $\theta$ is given in equation (\ref{eq-Mukai-isometry}).
We show that the composition is $Mon(M)$-equivariant, up to sign, hence defines 
the $Mon(X)$-submodule $\EE(X)$ in part (\ref{thm-item-the-Mon-equivariant-isometric-embedding-e}) of
Theorem \ref{thm-Q-4}, for any $X$ of $K3^{[n]}$-type.
%****************************************************************
% 
%****************************************************************
\subsubsection{Parallel transport operators of $K3^{[n]}$-type}
\label{sec-parallel-transport-operators-of-K3-n-type}
Let $X_1$ and $X_2$ be irreducible holomorphic symplectic manifolds of
$K3^{[n]}$-type. Denote by $\iota_{X_i}$ the natural $O(\widetilde{\Lambda})$-orbit 
of primitive isometric embedding of 
$H^2(X_i,\Integers)$ into the Mukai lattice $\widetilde{\Lambda}$, 
given in Corollary \ref{cor-third-characterization-of-Mon-2}.

\begin{thm}
\label{thm-necessary-and-sufficient-condition-to-be-a-parallel-transport-operator}
An isometry $g:H^2(X_1,\Integers)\rightarrow H^2(X_2,\Integers)$
is a parallel-transport operator, if and only if $g$ is orientation preserving and
\begin{equation}
\label{eq-g-pulls-back-orbit-to-orbit}
\iota_{X_1}=\iota_{X_2}\circ g.
\end{equation}
\end{thm}

\begin{proof}
Assume first that $g$ is a parallel-transport operator. 
Then $g$ lifts to a parallel-transport operator
$\tilde{g}:H^*(X_1,\Integers)\rightarrow H^*(X_2,\Integers)$. Now
$\tilde{g}$ induces a parallel-transport operators 
$\tilde{g}_4:Q^4(X_1,\Integers)\rightarrow Q^4(X_2,\Integers)$, 
as well as 
\[
Ad_{\tilde{g}}:
\Hom\left[H^2(X_1,\Integers),Q^4(X_1,\Integers)\right] \ \ \longrightarrow \ \ 
\Hom\left[H^2(X_2,\Integers),Q^4(X_2,\Integers)\right],
\]
given by $f\mapsto \tilde{g}_4\circ f\circ g^{-1}$. 
We have the equality $Ad_{\tilde{g}}(\EE_{X_1})=\EE_{X_2}$,
by the characterization of the
$Mon(X_i)$-module $\EE(X_i)$ provided in Theorem
\ref{thm-Q-4}. Hence, the equality 
(\ref{eq-g-pulls-back-orbit-to-orbit}) holds, by construction of $\iota_{X_i}$.

Conversely, assume that the isometry $g$ satisfies the equality 
(\ref{eq-g-pulls-back-orbit-to-orbit}). There
exists a parallel-transport operator 
$f:H^2(X_1,\Integers)\rightarrow H^2(X_2,\Integers)$, since $X_1$ and $X_2$
are deformation equivalent. Hence, the equality
$\iota_{X_1}=\iota_{X_2}\circ f$ holds, as well. We get the equality
$\iota_{X_1}=\iota_{X_1}\circ f^{-1}g$. We conclude that 
$f^{-1}g$ belongs to $Mon^2(X_1)$, by Corollary
\ref{cor-third-characterization-of-Mon-2}.
The equality $g=f(f^{-1}g)$ represents $g$ as a composition of two
parallel-transport operators. Hence, $g$ is a parallel-transport operator.
\end{proof}

The following statement is an immediate corollary of Theorems
\ref{thm-Hodge-theoretic-Torelli} and 
\ref{thm-necessary-and-sufficient-condition-to-be-a-parallel-transport-operator}.

\begin{cor}
\label{cor-Hodge-theoretic-Torelli-for-K3-n-type}
Let $X$ and $Y$ be two manifolds of $K3^{[n]}$-type.
\begin{enumerate}
\item
\label{thm-item-bimeromorphic-iff}
$X$ and $Y$ are bimeromorphic, if and only if there exists a Hodge-isometry
$f:H^2(X,\Integers)\rightarrow H^2(Y,\Integers)$, 
satisfying $\iota_X=\iota_Y\circ f$.
\item
$X$ and $Y$ are isomorphic, if and only if there exists a Hodge-isometry $f$ as in part
(\ref{thm-item-bimeromorphic-iff}), which maps some K\"{a}hler class of $X$ to a 
K\"{a}hler class of $Y$.
\end{enumerate}
\end{cor}

We do not require $f$ in part (\ref{thm-item-bimeromorphic-iff}) to be orientation preserving,
since if it is not then $-f$ is,  
and the orbits $\iota_Y\circ f$ and $\iota_Y\circ (-f)$ are equal. 

Let $\tau$ be the set of connected components of the moduli space 
of marked pairs $(X,\eta)$, where $X$ is of $K3^{[n]}$-type, and
$\eta:H^2(X,\Integers)\rightarrow \Lambda$ is an isometry.
Denote by $\FM^\tau_\Lambda$ the moduli space of isomorphism classes
of marked pairs $(X,\eta)$, where $X$ is of $K3^{[n]}$-type.
The group $O(\Lambda)$ acts on the set $\tau$ and the stabilizer of a
connected component $\FM^t_\Lambda$, $t\in \tau$, is the monodromy
group $Mon^2(\FM^t_\Lambda)\subset O(\Lambda)$ 
(Definition \ref{def-polarized-monodromy-subgroup-of-O-Lambda}).
Let 
\[
{\rm orb} \ : \ \FM^\tau_\Lambda \ \ \rightarrow \ \ 
O(\Lambda,\widetilde{\Lambda})/O(\widetilde{\Lambda})
\]
be the map given by $(X,\eta)\mapsto \iota_X\circ \eta^{-1}$.
Let $\orient:\FM^\tau_\Lambda\rightarrow\Orient(\Lambda)$ be the map
given in equation (\ref{eq-orient}).
The characterization of the monodromy group in Corollary 
\ref{cor-third-characterization-of-Mon-2} 
yields the following enumeration of $\tau$.

\begin{cor}
\label{cor-enumeration-of-connected-components}
The map $({\rm orb},\orient):\FM^\tau_\Lambda\rightarrow 
O(\Lambda,\widetilde{\Lambda})/O(\widetilde{\Lambda})\times \Orient(\Lambda)$
factors through a bijection 
\[
\tau \ \ \ \rightarrow \ \ \ 
O(\Lambda,\widetilde{\Lambda})/O(\widetilde{\Lambda})\times \Orient(\Lambda).
\]
\end{cor}

%****************************************************************
% 
%****************************************************************
\subsection{A numerical determination of the fundamental exceptional chamber}
\label{sec-numerical-determination-of-BK-X}

\begin{defi}
\label{def-monodromy-reflective}
A class $\ell\in H^{1,1}(X,\Integers)$ is called 
{\em monodromy-reflective}, if $\ell$ is a primitive class, $(\ell,\ell)<0$, and 
$R_\ell$ is a monodromy operator.
A holomorphic line bundle $L\in \Pic(X)$ is called 
{\em monodromy-reflective}, if the class  $c_1(L)$ is {\em monodromy-reflective}. 
\end{defi}

Let $X$ be a manifold of $K3^{[n]}$-type, $n\geq 2$.
%Recall that a class $\ell\in H^{1,1}(X,\Integers)$ is monodromy-reflective,
%if $\ell$ is a primitive class of negative Beauville-Bogomolov degree
%and the reflection $R_\ell$ belongs to $Mon^2(X)$
%(Definition \ref{def-monodromy-reflective}).
In section \ref{sec-Monodromy-reflective-classes-K3-n-type}
%and \ref{sec-monodromy-invariants}
we classify monodromy-orbits of monodromy-reflective classes.
This is done in terms of  explicitly computable monodromy invariants.
In section \ref{sec-stably-prime-exceptional-classes-K3-n-type}
we describe the values of the monodromy invariants,
for which the monodromy-reflective class is stably prime-exceptional
(Theorem \ref{thm-stably-prime-exceptional-K3-n}).
%(\cite{markman-prime-exceptional}, Theorem 1.12).
When $X$ is projective,  
Theorems \ref{thm-numerical-characterization-of-BK-X} and
\ref{thm-stably-prime-exceptional-K3-n} combine 
to provide a  determination of the 
closure $\overline{\BK}_X$ of the birational K\"{a}hler cone in $\C_X$
in terms of explicitly computable invariants.

%****************************************************************
% 
%****************************************************************
\subsubsection{Monodromy-reflective classes of $K3^{[n]}$-type}
\label{sec-Monodromy-reflective-classes-K3-n-type}
Set $\Lambda:=H^2(X,\Integers)$. 
Recall that if $\ell\in\Lambda$ is monodromy-reflective, then $R_\ell$
acts on $\Lambda^*/\Lambda$ via multiplication by $\pm1$
(Lemma \ref{lemma-mon-2-maps-to-pm-1}).
The set of monodromy-reflective classes is determined by the following 
statement.

\begin{prop}
\label{prop-monodromy-reflective-of-K3-n-type}
(\cite{markman-prime-exceptional}, Proposition 1.5)
Let $\ell\in H^2(X,\Integers)$ be a primitive class of negative degree 
$(\ell,\ell)<0$. Then $R_\ell$ belongs to $Mon^2(X)$, if and only if $\ell$ has one of
the following two properties.
\begin{enumerate}
\item
\label{prop-item-minus-2-is-monodromy-reflective}
$(\ell,\ell)=-2$.
\item
\label{prop-item-2-minus-2n-may-be-monodromy-reflective}
$(\ell,\ell)=2-2n$, and $(n-1)$ divides the class $(\ell,\bullet)\in H^2(X,\Integers)^*.$
\end{enumerate}
$R_\ell$ acts on $\Lambda/\Lambda^*$ as the identity 
in case (\ref{prop-item-minus-2-is-monodromy-reflective}), 
and via multiplication by $-1$ 
in case (\ref{prop-item-2-minus-2n-may-be-monodromy-reflective}).
\end{prop}

%****************************************************************
% 
%****************************************************************
%\subsubsection{Monodromy-invariants of monodromy-reflective classes}
%\label{sec-monodromy-invariants}
%Let $I(X)\subset H^2(X,\Integers)$ be a $Mon^2(X)$-invariant subset
%and $\Sigma$ a set. 

%\begin{defi}
%A function $f:I(X)\rightarrow \Sigma$ is called a {\em monodromy invariant}, 
%if $f$ is constant on each $Mon^2(X)$-orbit in $I(X)$. 
%If, furthermore, $f$ factors through an injective function from
%$I(X)/Mon^2(X)$ into $\Sigma$, we say that $f$ is a 
%{\em faithful} monodromy invariant.
%\end{defi}

%Clearly, the degree $(e,e)$ and the divisibility $(e,\bullet)$
%are monodromy invariants of the class $e$. 
%We describe next a faithful monodromy invariant for monodromy-reflective
%classes.

Given a primitive class $e\in H^2(X,\Integers)$, we denote by $\div(e,\bullet)$
the largest positive integer dividing the class $(e,\bullet)\in H^2(X,\Integers)^*$.
Let $\R_n(X)\subset H^2(X,\Integers)$ be the subset
of primitive classes of degree $2-2n$, such that $n-1$ divides $\div(e,\bullet)$. 
Let $\ell\in\R_n(X)$ and 
%be a monodromy reflective class in $H^2(X,\Integers)$
%with $(\ell,\ell)=2-2n$ and such that $n-1$ divides $\div(\ell,\bullet)$.
choose an embedding $\iota:H^2(X,\Integers)\hookrightarrow \widetilde{\Lambda}$
in the natural 
orbit $\iota_X$ provided by Corollary \ref{cor-third-characterization-of-Mon-2}.
Choose a generator $v\in \widetilde{\Lambda}$ of the rank $1$ sublattice 
orthogonal to the image of $\iota$. 
Set $e:=\iota(\ell)$ and let 
\begin{equation}
\label{eq-L}
L\subset \widetilde{\Lambda}
\end{equation}
be the saturation of the rank $2$ sublattice spanned by $e$ and $v$.

\begin{defi}
Two pairs $(L_i,e_i)$, $i=1,2$, 
each consisting of a lattice $L_i$ and a class $e_i\in L_i$, are said to be {\em isometric},
if there exists an isometry $g:L_1\rightarrow L_2$, such that $g(e_1)=e_2$.
\end{defi}

Given a rank $2$ lattice $L$, let $I_n(L)\subset L$ be the subset
of primitive classes $e$ with $(e,e)=2-2n$.

\begin{lem}
There exists a natural one-to-one correspondence between the orbit set $I_n(L)/O(L)$ 
and the set of isometry classes
of pairs $(L',e')$, such that $L'$ is isometric to $L$
and $e'$ is a primitive class in $L'$ with $(e',e')=2-2n$.
\end{lem}

\begin{proof}
Let $\P(L,n)$ be the set of isometry classes
of pairs $(L',e')$ as above. Define the map $f:\P(L,n)\rightarrow I_n(L)/O(L)$
as follows. Given a pair $(L',e')$ representing a class in $\P(L,n)$, 
choose an isometry $g:L'\rightarrow L$ and set $f(L',e'):=O(L)g(e')$. The map $f$ is well defined,
since the orbit $O(L)g(e')$ is clearly independent of the choice of $g$. 
The map $f$ is surjective, since given $e\in I_n(L)$, $f(L,e)=O(L)e$. 
If $f(L_1,e_1)=f(L_2,e_2)$, then there exist isometries $g_i:L_i\rightarrow L$
and an element $h\in O(L)$, such that $g_2(e_2)=hg_1(e_1)$. 
Then $g_2^{-1}hg_1$ is an isometry from $(L_1,e_1)$ to $(L_2,e_2)$.
Hence, the map $f$ is injective.
\end{proof}

Let $U$ be the unimodular hyperbolic plane. Let $U(2)$ be the
rank $2$ lattice with Gram matrix
$\left(\begin{array}{cc}
0 & -2 \\ -2 & 0
\end{array}
\right)$
and let $D$ be the rank $2$ lattice with Gram matrix
$\left(\begin{array}{cc}
-2 & 0 \\ 0 & -2
\end{array}
\right)$.

\begin{prop}
(\cite{markman-prime-exceptional}, Propositions 1.8 and  6.2)
\begin{enumerate}
\item
If $(\ell,\ell)=-2$ then the $Mon^2(X)$-orbit of $\ell$ is determined by 
$\div(\ell,\bullet)$.
\item
Let $\ell\in \R_n(X)$. 
\begin{enumerate}
\item
The lattice $L$, given in (\ref{eq-L}), 
is isometric to one of the lattices $U$, $U(2)$, or $D$.
\item
Let 
$
f\ : \ \R_n(X) \ \ \ \longrightarrow \ \ \ I_n(U)/O(U) \ \cup \ 
I_n(U(2))/O(U(2)) \ \cup \ I_n(D)/O(D)
$
be the function, sending a class $\ell$ to the isometry class of the pair
$(L,\iota(\ell))$.
Then the values $\div(\ell,\bullet)$ and $f(\ell)$ determine the $Mon^2(X)$-orbit
of $\ell$.
\end{enumerate}
\end{enumerate}
\end{prop}

The values of the function $f$ can be conveniently enumerated and 
calculated as follows. Set $e:=\iota(\ell)\in L$. 
Let $\rho$ be the largest integer, such that $(e+v)/\rho$ is an integral class of $L$.
Let $\sigma$ be the largest integer, such that $(e-v)/\sigma$ 
is an integral class of $L$. 
If $\div(\ell,\bullet)=n-1$ and $n$ is even, set $\{r,s\}(\ell)=\{\rho,\sigma\}$.
Otherwise, set $\{r,s\}(\ell)=\{\frac{\rho}{2},\frac{\sigma}{2}\}$.
The unordered pair $\{r,s\}:=\{r,s\}(\ell)$ has the following properties.

\begin{prop}
\label{prop-r-s}
(\cite{markman-prime-exceptional}, Lemma 6.4)
\begin{enumerate}
\item
The isometry class of the lattice $L$ and the product $rs$ are determined in
terms of $(\ell,\ell)$, $\div(\ell,\bullet)$, $n$, and $\{\rho,\sigma\}$
by the following table.

\begin{tabular}{|c|c|c|c|c|c|c|c|}  \hline 
\hspace{1ex}
& $(\ell,\ell)$ & $\div(\ell,\bullet)$ & $n$ & $\rho\sigma$ & $L$ & 
$\{r,s\}$ & $r\cdot s$ 
\\
\hline
1) &$2-2n$ & $2n-2$ & $\geq 2$ & $4n-4$ & $U$ & 
$\{\frac{\rho}{2},\frac{\sigma}{2}\}$ & $n-1$
\\
\hline
2) & $2-2n$ & $n-1$ & even & $n-1$ & $D$ & $\{\rho,\sigma\}$ & $n-1$
\\
\hline
3) & $2-2n$ & $n-1$ & odd & $2n-2$ & $U(2)$ & 
$\{\frac{\rho}{2},\frac{\sigma}{2}\}$ & $(n-1)/2$
%\\
%\hline
%4) & $2-2n$ & $n-1$ & $\equiv 1$ modulo $8$ & $2n-2$ & $U(2)$ & 
%$\{\frac{\rho}{2},\frac{\sigma}{2}\}$ & $(n-1)/2$
\\
\hline
4) & $2-2n$ & $n-1$ & $\equiv 1$ modulo $8$ & $n-1$ & $D$ & 
$\{\frac{\rho}{2},\frac{\sigma}{2}\}$ & $(n-1)/4$
\\
\hline
\end{tabular}
\item
The pair $\{r,s\}$ consists of relatively prime positive integers.
All four rows in the above table do occur, 
and every relatively prime decomposition $\{r,s\}$ of the
integer in the rightmost column occurs, for some $\ell\in \R_n(X)$.
\item
If $\ell\in\R_n(X)$, then 
$\div(\ell,\bullet)$ and $\{r,s\}(\ell)$ determine the $Mon^2(X)$-orbit of $\ell$.
\end{enumerate}
\end{prop}

%****************************************************************
% 
%****************************************************************
\subsubsection{Stably prime-exceptional classes of $K3^{[n]}$-type}
\label{sec-stably-prime-exceptional-classes-K3-n-type}
\begin{thm}
\label{thm-stably-prime-exceptional-K3-n}
(\cite{markman-prime-exceptional}, Theorem 1.12).
Let $\kappa\in H^{1,1}(X,\RealNumbers)$ be a K\"{a}hler class
and $L$ a monodromy reflective line bundle.
Set $\ell:=c_1(L)$. Assume that $(\kappa,\ell)>0$.
\begin{enumerate}
\item
If $(\ell,\ell)=-2$, then $L^k$ is stably prime-exceptional, where
\[
k=\left\{\begin{array}{ccccc}
2, & \mbox{if} & \div(\ell,\bullet)=2 & \mbox{and} & n=2,
\\
1, & \mbox{if} & \div(\ell,\bullet)=2 & \mbox{and} & n>2,
\\
1 & \mbox{if}  & \hspace{1ex}\div(\ell,\bullet)=1.
\end{array}\right.
\]
\item
If $\div(\ell,\bullet)=2n-2$ and $\{r,s\}(\ell)=\{1,n-1\}$, 
then $L^2$ is stably prime-exceptional.
\item
If $\div(\ell,\bullet)=2n-2$ and $\{r,s\}(\ell)=\{2,(n-1)/2\}$, 
then $L$ is stably prime-exceptional.
\item
If $\div(\ell,\bullet)=n-1$, $n$ is even, 
and $\{r,s\}(\ell)=\{1,n-1\}$, then $L$ is stably prime-exceptional.
\item
If $\div(\ell,\bullet)=n-1$, $n$ is odd, 
and $\{r,s\}(\ell)=\{1,(n-1)/2\}$, then $L$ is stably prime-exceptional.
\item
In all other cases, $H^0(L^k)$ vanishes, and so $L^k$ is not
stably prime-exceptional,
for every non-zero integer $k$.
\end{enumerate}
\end{thm}

When $X$ is projective, 
Proposition \ref{prop-monodromy-reflective-of-K3-n-type} and
Theorem  \ref{thm-stably-prime-exceptional-K3-n}
determine the set $\Spe\subset H^{1,1}(X,\Integers)$, 
of stably prime-exceptional classes, and hence
also the fundamental exceptional chamber $\FE_X$, by Proposition 
\ref{prop-fundamental-exceptional-chamber-is-indeed-an-exceptional-one}.

The proof of Theorem \ref{thm-stably-prime-exceptional-K3-n} 
has two ingredients. First we deform the pair $(X,L)$ to a pair 
$(M,L_1)$,  where $M$ is a moduli space of sheaves on a projective $K3$ surface,
and $L_1$ is a monodromy-reflective line bundle
with the same monodromy invariants. Then $L$ is stably prime-exceptional,
if and only if $L_1$ is, by Proposition 
\ref{prop-deformation-of-stably-prime-exceptional}.
We then laboriously check an example, one for each value of the 
monodromy invariants $n$, $(\ell,\ell)$, $\div(\ell,\bullet)$, 
and $\{r,s\}(\ell)$,
and show that either $R_\ell$ is induced by a birational map $f:M\rightarrow M$,
such that $f^*(L_1)=L_1^{-1}$, or that the linear system 
$\linsys{L_1^k}$ consists of
a single prime exceptional divisor, for
the power $k$ prescribed by Theorem \ref{thm-stably-prime-exceptional-K3-n}.

The two possible values of the degree $-2$ or $2-2n$, of a prime exceptional
divisor, correspond to two types of well known constructions
in the theory of moduli spaces of sheaves on a $K3$ surface $S$.
We briefly describe these constructions below.

Pairs $(M,\StructureSheaf{M}(E))$, 
where $M:=M_H(v)$ is a moduli space of 
$H$-stable coherent sheaves 
of class $v\in K(S)$, and $E$ is a prime exceptional 
divisor of Beauville-Bogomolov  degree $-2$, arise as follows.
The Mukai isometry (\ref{eq-Mukai-isometry}) associates to the line bundle 
$\StructureSheaf{M}(E)$ a class $e\in v^\perp$, with $(e,e)=-2$.
In the examples considered in \cite{markman-prime-exceptional}, 
$e$ is the class of an $H$-stable 
sheaf $F$ on $S$. Such a sheaf is necessarily rigid,
i.e., $\Ext^1(F,F)=0$. Indeed, 
\[
\dim\Ext^1(F,F)=\dim\Hom(F,F)+\dim\Ext^2(F,F)-\chi(F^\vee\otimes F)=1+1-2=0.
\]
Furthermore, the moduli space $M_H(e)$ is connected, by a theorem of Mukai,
and consists of the single point $\{F\}$ (see \cite{mukai-hodge}).
The prime exceptional divisor $E$ is the Brill-Noether locus 
\[
\{V\in M_H(v) \ : \  \dim\Ext^1(F,V)>0\}.
\]

Specific examples are easier to describe using Mukai's notation.
Recall Mukai's isomorphism  
\begin{equation}
\label{eq-Mukai-vector}
ch(\bullet)\sqrt{td_S} \ : \ K(S) \ \ \longrightarrow \ \ H^*(S,\Integers),
\end{equation}
sending a class $v\in K(S)$ to the integral singular cohomology group. 
Let $D:H^*(S,\Integers)\rightarrow H^*(S,\Integers)$ be the automorphism 
acting by $(-1)^i$ on $H^{2i}(S,\Integers)$.
The homomorphism (\ref{eq-Mukai-vector}) 
is an isometry once we endow $H^*(S,\Integers)$ with the pairing
\[
(x,y) \ \ \ := \ \ \  
%-\chi(x^\vee\otimes y) \ \ \ = \ \ \  
 -\int_S D(x)\cup y,
\]
by the Hirzebruch-Riemann-Roch theorem and the definition
of the Mukai pairing in equation (\ref{eq-Mukai-pairing}).
We have $ch(v)\sqrt{td_S}=(r,c_1(v),s)$, where $r=\rank(v)$,
$s=\chi(v)-r$, and we identify $H^0(S,\Integers)$ and $H^4(S,\Integers)$ with $\Integers$,
using the classes Poincar\'{e}-dual to $S$ and to a point.
Given two classes $v_i\in K(S)$, with $\rank(v_i)=r_i$, $c_1(v_i)=\alpha_i$,
and $s_i:=\chi(v_i)-r_i$, then
\[
(v_1,v_2) \ \ = \ \ \left(\int_S\alpha_1\alpha_2\right)-r_1s_2-r_2s_1.
\]

\begin{example}
Following is a simple example in which a prime exceptional divisor $E$
of degree $-2$ and divisibility $\div([E],\bullet)=1$ is realized as a Brill-Noether locus. 
Consider a $K3$ surface $S$,
containing a smooth rational curve $C$.
Consider the Hilbert scheme $M:=S^{[n]}$ as the moduli
space of ideal sheaves of length $n$ subschemes.
Let $F$ be the torsion sheaf $\StructureSheaf{C}(-1)$, supported 
on $C$ as a line bundle of degree $-1$. 
Let $v\in K(S)$ be the class of an ideal sheaf in $S^{[n]}$ and $e$ the class of $F$.
The Mukai vector of $v$ is $(1,0,1-n)$, that of $e$ is $(0,[C],0)$, and 
$(v,e)=0$. 
Let $E\subset M$ be the divisor of ideal sheaves $I_Z$ of subscheme 
$Z$ with non-empty intersection $Z\cap C$. 
The space $\Hom(F,I_Z)$ vanishes for all $I_Z\in M$, and so
$\dim\Ext^1(F,I_Z)=\dim\Ext^2(F,I_Z)$, for all $I_Z\in M$. Now, 
$\Ext^2(F,I_Z)\cong\Hom(I_Z,F)^*$ vanishes, if and only if $Z\cap C=\emptyset$.
Hence, $\Ext^1(F,I_Z)\neq 0$, if and only if $I_Z$ belongs to $E$.
See \cite{markman-brill-noether,yoshioka-brill-noether} for many more examples
of prime exceptional divisors $E$ of degree $-2$ and $\div([E],\bullet)=1$.
See \cite{markman-prime-exceptional}, Lemma 10.7 for the case $(e,e)=-2$,
$\div(e,\bullet)=2$, and $n\equiv 2$ modulo $4$.
\end{example}

Jun Li constructed a birational morphism from 
the moduli space of Gieseker-Maruyama $H$-stable sheaves on a $K3$ surface
to the Uhlenbeck-Yau compactification of the moduli space of 
$H$-slope-stable locally-free sheaves \cite{jun-li}.
The examples of prime exceptional  divisors of degree $2-2n$ on a moduli space of sheaves,
provided in \cite{markman-prime-exceptional},
were all constructed as exceptional divisors for Jun Li's morphism.

\begin{example}
The simplest example is the
Hilbert-Chow morphism, from the Hilbert scheme $S^{[n]}$, $n\geq 2$, to the 
symmetric product $S^{(n)}$ of a $K3$ surface $S$,
where the exceptional divisor $E$ is the big diagonal. 
The Mukai vector of the ideal sheaf is $v=(1,0,1-n)$.
In this case $[E]=2\delta$, where $\delta=(1,0,n-1)$. Note that $(\delta,\delta)=2-2n$.
The second cohomology of $S^{[n]}$ is an orthogonal direct sum 
$H^2(S,\Integers)\oplus \Integers\delta$, by \cite{beauville} or by Mukai's isometry
(\ref{eq-Mukai-isometry}). Hence, 
$\div(\delta,\bullet)=2n-2$. The  largest integer $\rho$ dividing $\delta+v=(2,0,0)$
is $2$ and the largest integer $\sigma$ dividing $\delta-v=(0,0,2n-2)$ is $2n-2$.
Hence, 
$\{r,s\}(\delta)=\{1,n-1\}$, by Proposition \ref{prop-r-s} and Equation
(\ref{eq-iota-is-theta-inverse}).
\end{example}

\begin{example}
\label{example-monodromy-reflective-but-not-spe}
Consider, more generally, the moduli space $M_H(r,0,-s)$ of $H$-stable 
sheaves with Mukai vector $v=(r,0,-s)$, satisfying $s>r\geq 1$ and
$\gcd(r,s)=1$. Then $M_H(r,0,-s)$ is of $K3^{[n]}$-type, $n=rs+1$.
The Mukai vector $e:=(r,0,s)\in v^\perp$ maps to
a monodromy-reflective class  $\ell\in H^2(M_H(v),\Integers)$ 
of degree $(\ell,\ell)=2-2n$, divisibility $\div(\ell,\bullet)=2n-2$, 
and $\{r,s\}(\ell)=\{r,s\}$, 
by Proposition \ref{prop-r-s} and Equation (\ref{eq-iota-is-theta-inverse}).
When  $r=2$, $\ell$ is the class of the exceptional divisor $E$ of Jun Li's morphism.
$E$ is the locus of sheaves, which are not locally-free or not $H$-slope-stable
(\cite{markman-prime-exceptional}, 
Lemma 10.16).
When $r>2$, the exceptional locus has co-dimension $\geq 2$, and
no multiple of the class $\ell$ is effective.
Instead, the reflection $R_\ell$ is induced by the birational map
$f:M_H(r,0,-s)\rightarrow M_H(r,0,-s)$, sending a locally-free $H$-slope stable sheaf $F$
of class $(r,0,-s)$ to the dual sheaf $F^*$
(\cite{markman-prime-exceptional}, Proposition 11.1).
\end{example}

\begin{rem}
\label{rem-BK-depends-on-orbit-iota-X}
Fix an integer $n>0$, such that $n-1$ is not a prime power, and  consider 
all possible factorizations $n-1=rs$, with $s>r\geq 1$ and
$\gcd(r,s)=1$. 
The sub-lattice $(r,0,-s)^\perp$ of the Mukai lattice of a $K3$ surface $S$
is the orthogonal direct sum $H^2(S,\Integers)\oplus \Integers (r,0,s)$.
We get the isometry 
\[
\theta \ : \ H^2(S,\Integers)\oplus \Integers (r,0,s)\ \ \longrightarrow \ \ 
H^2\left(M_H(r,0,-s),\Integers\right),
\]
using Mukai's isometry given in equation (\ref{eq-Mukai-isometry}). 
Let $n-1=r_1s_1=r_2s_2$ be two different such factorizations.
Then the two moduli spaces $M_H(r_1,0,-s_1)$ and  $M_H(r_2,0,-s_2)$,
considered in Example \ref{example-monodromy-reflective-but-not-spe}, come
with a natural Hodge isometry
\[
g \ : \ H^2(M_H(r_1,0,-s_1),\Integers) \ \ \longrightarrow \ \ 
H^2(M_H(r_2,0,-s_2),\Integers),
\]
which restricts as the identity on the direct summand 
$\theta\left(H^2(S,\Integers)\right)$ and maps
the class $\ell_1:=\theta(r_1,0,s_1)\in H^2(M_H(r_1,0,-s_1),\Integers)$
to the class $\ell_2:=\theta(r_2,0,s_2)\in H^2(M_H(r_2,0,-s_2),\Integers)$.
The Hodge isometry $g$ is not a parallel-transport operator, 
since the monodromy-invariants
$\{r,s\}(\ell_i)=\{r_i,s_i\}$ are distinct.
Indeed, these moduli spaces
are not birational in general (\cite{markman-constraints}, Proposition 4.10). 
Furthermore, if $n-1=rs$ is such a factorization with $r>2$, then 
%the Hodge isometry 
%$g:H^2(S^{[n]},\Integers)\rightarrow H^2(M_H(r,0,-s),\Integers)$, compatible with 
%the above orthogonal direct sum decomposition, maps 
the birational K\"{a}hler cones $\BK_{S^{[n]}}$ and $\BK_{M_H(r,0,-s)}$
are not isometric in general. 
Indeed, $S^{[n]}$ admits a stably prime-exceptional class, 
while $M_H(r,0,-s)$ does not, for a $K3$ surface with a suitably chosen Picard lattice.
\end{rem}

%****************************************************************
% 
%****************************************************************
\section{Open problems}
\label{sec-open-problems}

Following is a very brief list of central open problems closely related to this survey.
See \cite{beauville-list} for a more complete recent survey of open problems
in the subject of irreducible holomorphic symplectic manifolds.

\begin{question}
Let $X$ be one of the known examples of irreducible holomorphic symplectic manifolds,
i.e., of $K3^{[n]}$-type, a generalized Kummer variety, or one of the two exceptional 
examples of O'Grady \cite{ogrady-ten,ogrady-six}.
%Let $\Lambda$ be a lattice isometric to $H^2(X,\Integers)$.
Let $Y$ be an irreducible holomorphic symplectic manifold, with $H^2(Y,\Integers)$
isometric to $H^2(X,\Integers)$. Is $Y$ necessarily deformation equivalent to $X$?
\end{question}

Let $\Lambda$ be a lattice isometric to $H^2(X,\Integers)$.
At present it is only known that the number of deformation types
of irreducible holomorphic symplectic manifolds of
a given dimension $2n$, and with second cohomology lattice  isometric to $\Lambda$,
is finite \cite{huybrechts-finiteness}. 
The moduli space $\FM_\Lambda$, of isomorphism classes of marked 
pairs $(X,\eta)$, with $X$ of dimension $2n$ and $\eta:H^2(X,\Integers)\rightarrow \Lambda$
an isometry, has finitely many connected components, by Huybrechts' result and 
Lemma \ref{lemma-tau-is-finite}. O'Grady has made substantial progress towards 
the proof of uniqueness of the deformation type in case the dimension is $4$
and the lattice $\Lambda$ is of $K3^{[2]}$-type
\cite{ogrady-msri-lecture-notes}.

\begin{problem}
\label{problem-Kahler-cone}
Let $X$ be an irreducible holomorphic symplectic manifold of $K3^{[n]}$-type, $n\geq 2$.
Determine the K\"{a}hler-type chamber (Definition \ref{def-Kahler-type-chambers})
in the fundamental exceptional chamber $\FE_X$ of $X$, 
containing a given very general class $\alpha\in \FE_X$,
in terms of the weight $2$ integral Hodge structure $H^2(X,\Integers)$, 
the Beauville-Bogomolov pairing, 
and the orbit $\iota_X$ of isometric embeddings of $H^2(X,\Integers)$
in the Mukai lattice, given in Corollary \ref{cor-third-characterization-of-Mon-2}. 
\end{problem}

Note that the data specified in Problem \ref{problem-Kahler-cone}
determines the isomorphism class of an irreducible holomorphic symplectic manifold
$Y$, bimeromorphic to $X$, 
%an element $w\in W_{Exc}(X)$, 
and an $\Aut(X)\times \Aut(Y)$-orbit\footnote{The orbit of $f$
is the set $\{g_1fg_2^{-1} \ : \ g_1\in \Aut(X), \ g_2\in \Aut(Y)\}$.} 
of a bimeromorphic map $f:Y\rightarrow X$, such that $f^*(\alpha)$ is a 
K\"{a}hler class on $Y$, by 
%Theorem \ref{thm-decomposition-of-a-Hodge-isometry-as-wg} 
Corollaries \ref{cor-stabilizer-of-birational-Kahler-cone} 
and \ref{cor-Hodge-theoretic-Torelli-for-K3-n-type}. 
The homomorphism $f^*$ takes the K\"{a}hler-type chamber in Problem 
\ref{problem-Kahler-cone} to $\K_Y$.
%the Hodge-theoretic Torelli Theorem
%\ref{thm-Hodge-theoretic-Torelli}, Lemma \ref{lemma-on-semi-chambers}, 
%and the determination of parallel-transport operators of
%$K3^{[n]}$-type provided in Theorem
%\ref{thm-necessary-and-sufficient-condition-to-be-a-parallel-transport-operator}. 
Hassett and Tschinkel formulated a precise conjectural 
solution to problem \ref{problem-Kahler-cone}
%the determination of the K\"{a}hler cone 
\cite{hassett-tschinkel-intersection-numbers}.
The K\"{a}hler cone, according to 
their conjecture, does not depend on the orbit $\iota_X$.
The birational K\"{a}hler cone does, as we saw in Remark
\ref{rem-BK-depends-on-orbit-iota-X}.

%Let $Torelli$ be the category, whose objects are irreducible
%holomorphic symplectic manifolds and whose morphisms 
%$g\in \Hom_{Torelli}(X,Y)$ are parallel transport operators 
%$g:H^2(X,\Integers)\rightarrow H^2(Y,\Integers)$. 
%Fix a deformation type $\tau$ of some 
%irreducible holomorphic symplectic manifold $X_0$.
%Denote by $Torelli_\tau$ the full sub-category with objects
%of deformation type $\tau$.

%\begin{defi}
%Let $X$ be an irreducible holomorphic symplectic manifold. 
%A {\em Torelli data} for a variety $Y$ deformation equivalent to $X$
%consists of the following data.
%\begin{enumerate}
%\item
%The dimension of $Y$.
%\item
%The weight $2$ integral Hodge structure $H^2(Y,\Integers)$.
%\item
%The K\"{a}hler cone $\K_X\subset H^{1,1}(X,\RealNumbers)$.
%\item
%The Beauville-Bogomolov pairing on $H^2(Y,\Integers)$.
%\item
%If $Mon^2(X)$ is strictly smaller than $O^+[H^2(X,\Integers)]$, then 
%specify a functor

%some additional  monodromy invariant function 
%$\iota$, continuous over the base of flat families, 
%from the set of isomorphism classes of 
%manifolds deformation equivalent to $X$, o a finite set $\Sigma$.
%\end{enumerate}
%The invariant function $\iota$ 
%\end{defi}

\begin{problem}
\label{problem-parallel-transport-operators}
Find an explicit necessary and sufficient condition for a Hodge isometry 
$g:H^2(X,\Integers)\rightarrow H^2(Y,\Integers)$
to be a parallel-transport operator,
in the case $X$ and $Y$ are deformation equivalent to
generalized Kummer varieties, or 
to O'Grady's two exceptional examples.
\end{problem}

\begin{problem}
\label{problem-spe-in-other-examples}
Let $X$ be deformation equivalent to a generalized Kummer variety, or 
to one of O'Grady's two exceptional examples.
Find an explicit necessary and sufficient condition for a class $\ell\in H^{1,1}(X,\Integers)$
to be stably prime-exceptional (Definition \ref{def-stably-prime-exceptional}).
\end{problem}

Problem \ref{problem-spe-in-other-examples} is solved in the
$K3^{[n]}$-type case (Proposition \ref{prop-monodromy-reflective-of-K3-n-type} and
Theorem \ref{thm-stably-prime-exceptional-K3-n}). 
A solution to problem \ref{problem-spe-in-other-examples} yields a determination of
the fundamental exceptional chamber $\FE_X$, by Proposition
\ref{prop-intro-fundamental-exceptional-chamber-in-indeed-such}, and of
the closure of the birational K\"{a}hler cone, by
Proposition \ref{prop-closure-of-FE-and-FU-are-equal}.
%Theorem \ref{thm-numerical-characterization-of-BK-X}.
Once solutions to Problems \ref{problem-parallel-transport-operators}
and \ref{problem-spe-in-other-examples} 
are provided, the analogue of Problem \ref{problem-Kahler-cone} 
may be formulated as well.

\begin{question}
Is the monodromy group $Mon^2(X)$, 
of an irreducible holomorphic symplectic manifold $X$, 
necessarily a normal subgroup of the isometry group of
$H^2(X,\Integers)$?
\end{question}

Let $X$ be deformation equivalent to a generalized Kummer variety of dimension $2n$, $n\geq 2$.
Then $H^2(X,\Integers)$ is isometric to the lattice 
$\Lambda:=U\oplus U\oplus U\oplus \Integers\delta$,
where $U$ is the unimodular rank $2$ lattice
of signature $(1,1)$, and 
$(\delta,\delta)=-2-2n$ (see \cite{beauville,yoshioka-abelian-surface}). 
Given a class $u\in H^2(X,\Integers)$ with $(u,u)=2$, set $\rho_u=-R_u$. 
If $(u,u)=-2$, set $\rho_u:=R_u$, as in equaton (\ref{eq-N}).

\begin{conj}
$Mon^2(X)$ is equal to the subgroup $\N(X)$ of the signed isometry group 
$O^+[H^2(X,\Integers)]$, generated by products 
$\rho_{u_1}\cdots \rho_{u_{k}}$, where 
$(u_i,u_i)=-2$, for an even number of indices $i$,
and $(u_i,u_i)=2$ for the rest of the indices $i$.
\end{conj}

The inclusion $\N(X)\subset Mon^2(X)$ was proven by the author in an
unpublished work. When $n=2$, the equality $\N(X)= Mon^2(X)$
follows from the Global Torelli Theorem \ref{thm-global-torelli}
and Namikawa's counter example to the naive Hodge theoretic Torelli
statement \cite{namikawa-torelli}. If $H^{1,1}(X,\Integers)$
is cyclic, generated by a class $u$ with $(u,u)=2$, then $\rho_u$ is induced 
by a regular anti-symplectic involution
$f:X\rightarrow X$, by the inclusion $\N(X)\subset Mon^2(X)$, the Hodge theoretic Torelli Theorem
\ref{thm-Hodge-theoretic-Torelli}, and the equality of the K\"{a}hler and positive cones of $X$ 
(Theorem \ref{thm-global-torelli} part \ref{thm-item-cyclic-picard-and-projective-imply-K-X=C-X}).
The analogous statement in the $K3^{[n]}$ case is treated in \cite{ogrady-involutions,ogrady-epw}.

Let $X$ be an irreducible holomorphic symplectic manifold deformation equivalent to 
O'Grady's $10$-dimensional exceptional  example \cite{ogrady-ten}. 
Then $H^2(X,\Integers)$ is isometric to the orthogonal direct sum
of $H^2(S,\Integers)\oplus G_2$, where
$S$ is a $K3$ surface, and $G_2$ is the negative definite root lattice 
of type $G_2$, with Gram matrix
$\left(\begin{array}{cc}
-2 & 3\\
3 & -6
\end{array}
\right)$ (see \cite{rapagnetta}).
The isometry group $O(G_2)$ is equal to the Weyl group of $G_2$ and its extension to $H^2(X,\Integers)$,
via the trivial action on  $H^2(S,\Integers)$, is contained in $Mon^2(X)$, 
by (\cite{markman-galois}, Lemma 5.1).

\begin{conj}
$Mon^2(X)=O^+[H^2(X,\Integers)]$.
\end{conj}

There are many examples of non-isomorphic $K3$ surfaces with
equivalent bounded derived categories of coherent sheaves 
\cite{orlov}.

\begin{question}
\label{question-derived-equivalent}
Let $X$ and $Y$ be projective irreducible holomorphic symplectic manifolds,
such that $H^2(X,\Integers)$ and $H^2(Y,\Integers)$ are Hodge isometric.
Are their bounded derived categories of coherent sheaves 
necessarily equivalent?
\end{question}

When $X=S_1^{[n]}$ and $Y=S_2^{[n]}$, where $S_1$ and $S_2$ are
$K3$ surfaces, the answer to Question \ref{question-derived-equivalent} 
is affirmative (see  the proof of \cite{ploog}, Proposition 10).
%See \cite{ploog}, Remark 11 (3), for such an example
%with non-birational $X$ and $Y$. 
See \cite{huybrechts-seatle} for a survey on the topic of question 
\ref{question-derived-equivalent}.

%Recall that a line bundle $L$ is monodromy-reflective, 
%if $c_1(L)$ is a primitive class, and
%the reflection $R_{c_1(L)}$ is a monodromy operator 
%(Definition \ref{def-monodromy-reflective}).

%\begin{problem}
%\label{problem-stably-prime-exceptional-classes}
%Find an algorithm to determine if some power of 
%a monodromy-reflective line bundle
%is stably prime-exceptional, for the known deformation types of
%irreducible holomorphic symplectic manifolds.
%\end{problem}

Recall that a class $\ell\in H^{1,1}(X,\Integers)$ is monodromy-reflective, 
if it is a primitive class, and
the reflection $R_{\ell}$ is a monodromy operator 
(Definition \ref{def-monodromy-reflective}).

\begin{question}
Let $\ell\in H^{1,1}(X,\Integers)$ be a monodromy-reflective class. 
Is there always some non-zero integer $\lambda$,
such that the class $\lambda(\ell,\bullet)\in H^2(X,\Integers)^*\cong H_2(X,\Integers)$
corresponds to an effective one-cycle?
\end{question}

An affirmative answer to the above question implies that the reflection $R_\ell$ 
can not be induced by a regular automorphism\footnote{A weaker version of this assertion, namely
the non-existence of a fixed-point free such automorphism $g$, 
is always true. Indeed, if $g^*=R_\ell$, and $g$ is a fixed-point-free
(necessarily symplectic) automorphism, then $g^2$ acts trivially on $H^2(X,\Integers)$.
Hence, $g^2$ is an isometry with respect to a K\"{a}hler metric.
It follows that $g$ has finite order, since it generates a discrete
subgroup of the compact isometry group. Thus,
$X/\langle g\rangle$ is a non simply connected holomorphic symplectic K\"{a}hler manifold,
with $h^{k,0}(X)=1$, for even $k$ in the range $0\leq k\leq \dim_\ComplexNumbers(X)$,
and $h^{k,0}(X)=0$, otherwise.
Such $X$ does not exist, by \cite{huybrechts-nieper}, Proposition A.1.} 
of $X$.
%(compare with Assumption 
%\ref{assumption-reflection-is-not-induced-by-regular-automorphism}).
It follows that the K\"{a}hler cone is contained in a unique 
chamber of the subgroup of $Mon^2_{Hdg}(X)$ generated by all reflections in $Mon^2_{Hdg}(X)$
(see Theorem \ref{thm-each-chamber-is-a-generalized-convex-polyhedron}).
%monodromy-reflective chamber (Definition \ref{def-monodromy-reflective}).

\begin{problem}
Prove an analogue 
of Proposition \ref{prop-druel}, about birational contractibility
of a prime exceptional divisor, for non-projective 
irreducible holomorphic symplectic manifolds. 
\end{problem}

Druel's proof of Proposition \ref{prop-druel} relies on results in the minimal model program,
which are currently not available in the K\"{a}hler category \cite{druel}.

\begin{question}
\label{question-semigroup-of-effective-divisors}
Let $X$ be a projective irreducible holomorphic symplectic manifold.
Is the semi-group $\Sigma$,  of effective divisor classes on  $X$, equal to the semi-group
$\Sigma'$ generated by the prime exceptional classes and integral points on the closure
$\overline{\BK}_X$ of the birational K\"{a}hler cone in $H^{1,1}(X,\RealNumbers)$?
\end{question}

The answer is affirmative for any $K3$ surface, even without the 
projectivity assumption (\cite{BHPV}, Ch. IIIV, Proposition 3.7). 
Stronger results hold true for projective $K3$ surfaces \cite{kovacs}.
The inclusion $\Sigma\subset \Sigma'$ is known in general,
%The effective semi-group is contained in the cone in question
%\ref{question-semigroup-of-effective-divisors}, 
by the divisorial Zariski decomposition (Theorem \ref{thm-Zariski-decomposition}). 
The integral points of $\C_X\cap\overline{\BK}_X$ are known to be contained in $\Sigma$. 
This is seen as follows.
The integral points of the positive cone are known to correspond to 
big line bundles, by (\cite{huybrechts-basic-results}, Corollary 3.10). 
Each integral point of $\C_X\cap\overline{\BK}_X$ thus coresponds to a big and nef
line bundle $L$ on some birational irreducible
holomorphic symplectic manifold $Y$, by Theorems
\ref{thm-kahler-cone} and \ref{thm-numerical-characterization-of-BK-X},
and so the cohomology groups $H^i(Y,L)$ vanish, for $i>0$,
by the Kawamata-Viehweg vanishing theorem. 
Set $\ell:=c_1(L)$. 
If $X$ is of $K3^{[n]}$-type or deformation equivalent to a generalized Kummer variety, then an
explicit formula is known for the  Euler characteristic $\chi(L)$ of a line bundle $L$,
in terms of its Beauville-Bogomolov degree $(\ell,\ell)$ 
(\cite{huybrechts-norway}, Examples 7 and 8).
One sees, in particular, that $\chi(L)>0$, if $(\ell,\ell)\geq 0$, and so $L$ is effective. 

An affirmative answer to  Question \ref{question-semigroup-of-effective-divisors} would thus follow, 
%It remains to settle Question \ref{question-semigroup-of-effective-divisors} in the case 
if one could prove that 
nef line bundles with $(\ell,\ell)=0$ are effective. Some experts conjectured that 
such line bundles are related to 
Lagrangian fibrations (\cite{markushevich}, Conjecture 1.7; \cite{sawon}, Conjecture \nolinebreak 1,
\cite{verbitsky-syz}, Conjecture 1.7).
We refer the reader also to the important work of Matsushita on Lagrangian fibrations
\cite{matsushita,matsushita-higher-direct-images} and to the survey
(\cite{beauville-list}, section 1.6).

\begin{question}
Which components, of the moduli spaces of polarized 
projective irreducible holomorphic symplectic manifolds, 
are unirational? Which are of general type?
\end{question}

Gritsenko, Hulek, and Sankaran had studied this question for 
fourfolds $X$ of $K3^{[2]}$-type, and for primitive 
polarizations $h\in H^2(X,\Integers)$, with $\div(h,\bullet)=2$.
%such that the orthogonal complement $h^\perp$ is an orthogonal direct sum
%$2U\oplus 2E_8(-1)\oplus $
Let $(h,h)=2d$. They show that for $d\geq 12$, the moduli space is of general type
(\cite{GHS}, Theorem 4.1). They use the theory of modular forms to show 
that the quotient of the period domain $\Omega_{h^\perp}^+$, 
given in equation (\ref{eq-Omega-h-perp-plus}), by 
the polarized monodromy group $Mon^2(X,h)$, is of general type.

On the other hand, unirational components are those likely to admit
explicit and very beautiful geometric descriptions  
\cite{beauville-donagi,debbare-voisin,iliev-ranestad,mukai-unirational,ogrady-epw}.

%****************************************************************
% Bibliography:
%****************************************************************

\end{document}